\documentclass[10pt]{amsart}
\usepackage{amsmath,amsfonts,amsthm}
\usepackage[dvips]{graphicx}
\usepackage{psfrag}

\theoremstyle{plain}
\newtheorem{maintheorem}{Theorem}

\newtheorem{theorem}{Theorem}[section]
\newtheorem{lemma}[theorem]{Lemma}
\newtheorem{corollary}[theorem]{Corollary}
\newtheorem{proposition}[theorem]{Proposition}
\theoremstyle{remark}
\newtheorem{definition}[theorem]{Definition}
\newtheorem{remark}[theorem]{Remark}

\newcommand{\RR}{{\mathbb R}}

\newcommand{\ZZ}{{\mathbb Z}}
\newcommand{\NN}{{\mathbb N}}
\newcommand{\CC}{{\mathbb C}}

\newcommand{\PC}{\mathbb{P}(\mathbb{C}^d)}
\newcommand{\Grass}{\operatorname{Grass}(\ell,d)}

\newcommand{\cB}{{\mathcal B}}
\newcommand{\cD}{{\mathcal D}}
\newcommand{\cE}{{\mathcal E}}

\newcommand{\cG}{{\mathcal G}}
\newcommand{\cH}{{\mathcal H}}
\newcommand{\cI}{{\mathcal I}}

\newcommand{\cM}{{\mathcal M}}
\newcommand{\cN}{{\mathcal N}}

\newcommand{\cP}{{\mathcal P}}
\newcommand{\cQ}{{\mathcal Q}}
\newcommand{\cR}{{\mathcal R}}
\newcommand{\cU}{{\mathcal U}}

\newcommand{\Chi}{{\mathcal X}}
\newcommand{\cZ}{{\mathcal Z}}

\newcommand{\hA}{{\hat A}}
\newcommand{\hB}{{\hat B}}
\newcommand{\hI}{{\hat I}}
\newcommand{\hx}{{\hat x}}
\newcommand{\hy}{{\hat y}}
\newcommand{\hz}{{\hat z}}
\newcommand{\hw}{{\hat w}}
\newcommand{\hm}{{\hat m}}
\newcommand{\hf}{{\hat f}}
\newcommand{\hg}{{\hat g}}
\newcommand{\hp}{{\hat p}}

\newcommand{\hpi}{\hat\pi}
\newcommand{\hmu}{\hat\mu}
\newcommand{\hnu}{\hat\nu}
\newcommand{\bmu}{\bar\mu}
\newcommand{\bx}{{\bar x}}
\newcommand{\tm}{{\tilde m}}
\newcommand{\FA}{F_A}

\newcommand{\hFA}{{\hat F_A}}
\newcommand{\hFB}{{\hat F_B}}

\newcommand{\hpFA}{{\hat F_{\pA}}}

\newcommand{\hSigma}{\hat\Sigma}
\newcommand{\hsig}{\hat\Sigma}
\newcommand{\fkI}{f^k_I}

\newcommand{\Jmuf}{J f}
\newcommand{\Jmux}{J_x}
\newcommand{\Jmufk}{J f^k}
\newcommand{\Jmufn}{J f^n}
\newcommand{\JmufkI}{J f^k_I}

\newcommand{\xn}{{x^n}}
\newcommand{\yn}{{y^n}}

\newcommand{\xj}{{x^j}}
\newcommand{\yj}{{y^j}}

\newcommand{\psig}{\Sigma^u}
\newcommand{\nsig}{\Sigma^s}
\newcommand{\px}{{x^u}}
\newcommand{\nx}{{x^s}}
\newcommand{\py}{{y^u}}
\newcommand{\ny}{{y^s}}
\newcommand{\pz}{{z^u}}
\newcommand{\nz}{{z^s}}
\newcommand{\pP}{P^u}
\newcommand{\nP}{P^s}
\newcommand{\pI}{[I]^u}

\newcommand{\pJ}{[J]^u}

\newcommand{\pmu}{\mu^u}
\newcommand{\nmu}{\mu^s}

\newcommand{\pom}{m^u}

\newcommand{\pA}{A^u}

\newcommand{\hpA}{\hA^u}
\newcommand{\hnA}{\hA^s}

\newcommand{\pf}{f^u}
\newcommand{\nf}{f^s}
\newcommand{\pphi}{\phi^u}

\newcommand{\pfk}{f^{u,k}}
\newcommand{\nfk}{f^{s,k}}
\newcommand{\pfl}{f^{u,l}}

\newcommand{\pfkI}{\pfk_I}
\newcommand{\nfkI}{\nfk_I}
\newcommand{\pflJ}{\pfl_J}

\newcommand{\pfJI}{f^{u,k+l}_{JI}}

\newcommand{\pJmufk}{J \pfk}

\newcommand{\pJmufkI}{J \pfkI}

\newcommand{\pJmuflJ}{J \pflJ}

\newcommand{\pJmufJI}{J f^{u,k+l}_{JI}}

\newcommand{\xplus}{{w}}
\newcommand{\hxplus}{{\hw}}
\newcommand{\xminus}{{x_-}}
\newcommand{\Wsloc}{W^s_{loc}}
\newcommand{\Wuloc}{W^u_{loc}}
\newcommand{\Hxy}{H_{\hx,\hy}}
\newcommand{\Hyx}{H_{\hy,\hx}}
\newcommand{\Hpz}{H_{\hp,\hz}}
\newcommand{\Hpw}{H_{\hp,\hw}}
\newcommand{\Hsxy}{H^s_{\hx,\hy}}
\newcommand{\Huxy}{H^u_{\hx,\hy}}

\newcommand{\Huyx}{H^u_{\hy,\hx}}
\newcommand{\Hsxx}{H^s_{\hx,\hx}}

\newcommand{\Hupz}{H^u_{\hp,\hz}}
\newcommand{\Hspzf}{H^s_{\hf^l(\hz),\hp}}
\newcommand{\psipz}{\psi_{p,z}}
\newcommand{\Hsxz}{H^s_{\hx,\hz}}

\newcommand{\Hszy}{H^s_{\hz,\hy}}

\newcommand{\Hsxyf}{H^s_{\hf(\hx),\hf(\hy)}}

\newcommand{\Huxyfn}{H^u_{\hf^n(\hx),\hf^n(\hx)}}
\newcommand{\Hxplus}{H_{\hxplus,\hx}}
\newcommand{\Hplusx}{H_{\hx,\hxplus}}

\newcommand{\weak}{weak$^*$\ }
\newcommand{\supp}{{\operatorname{supp\,}}}

\newcommand{\id}{{\operatorname{id}}}

\newcommand{\GL}{\operatorname{GL}(d,\CC)}

\newcommand{\Ec}{\operatorname{Ecc}}

\newcommand{\extl}[1]{\Lambda^\ell(#1)}
\newcommand{\extdl}[1]{\Lambda^{d-\ell}(#1)}
\newcommand{\extlv}[1]{\Lambda^\ell_v(#1)}
\newcommand{\extdlv}[1]{\Lambda^{d-\ell}_v(#1)}

\newcommand{\pextlv}[1]{\mathbb{P}\Lambda^\ell_v(#1)}

\newcommand{\diag}{\operatorname{diag}}
\newcommand{\osc}{\operatorname{osc}}
\newcommand{\UU}{\operatorname{U}}

\newcommand{\quand}{\quad\text{and}\quad}

\newcommand{\vep}{\varepsilon}

\title[Simplicity of Lyapunov spectra]
      {Simplicity of Lyapunov spectra: \\ A sufficient criterion}
\author{Artur Avila and Marcelo Viana}
\date{\today}

\thanks{Work supported by the Brazil-France Agreement in
Mathematics.  M.V. is also supported by Pronex and Faperj.  This research
was partially conducted during the period A.A. served as a Clay Research
Fellow.}

\address{
CNRS UMR 7599, Laboratoire de Probabilit\'es et Mod\`eles
Al\'eatoires, Universit\'e Pierre et Marie Curie, Bo\^\i te
Postale 188, 75252 Paris Cedex 05, France}
\email{artur@ccr.jussieu.fr}

\address{
IMPA -- Estrada D. Castorina 110, Jardim Bot\^anico, 22460-320 Rio
de Janeiro -- Brazil. }
\email{viana@impa.br}

\begin{document}

\begin{abstract}
We exhibit an explicit sufficient condition for the Lyapunov
exponents of a linear cocycle over a Markov map to have
multiplicity~$1$. This builds on work of Guivarc'h-Raugi and
Gol'dsheid-Margulis, who considered products of random matrices,
and of Bonatti-Viana, who dealt with the case when the base
dynamics is a subshift of finite type. Here the Markov structure
may have infinitely many symbols and the ambient space needs not
be compact. As an application, in another paper we prove the
Zorich-Kontsevich conjecture on the Lyapunov spectrum of the
Teichm\"uller flow in the space of translation surfaces.
\end{abstract}

\maketitle

\setcounter{tocdepth}{1} \tableofcontents

\section{Introduction and statements}

Let $\hf:\hSigma\to\hSigma$ be an invertible measurable map and
$\hA:\hSigma\to\GL$ be a measurable function with values in the
group of invertible $d\times d$ complex matrices. These data
define a linear cocycle $\hFA$ over the map $\hf$, through
$$
\hFA:\hSigma\times\CC^d\to\hSigma\times\CC^d, \quad
\hFA(\hx,v)=\big(\hf(\hx),\hA(\hx)v\big).
$$
Note that $\hFA^n(x,v)=(\hf^n(\hx),\hA^n(\hx))$, where
$\hA^n(\hx)=\hA(\hf^{n-1}(\hx))\cdots\hA(\hf(\hx))\hA(\hx)$ and
$\hA^n(\hx)$ is the inverse of $\hA^{-n}(\hf^n(\hx))$ if $n < 0$.

Let $\hmu$ be an $\hf$-invariant probability measure on $\hSigma$
relative to which the logarithms of the norms of $\hA$ and its
inverse are integrable. By the theorem of Oseledets~\cite{Os68},
at $\mu$-almost every $\hx\in\hSigma$ there exist numbers
$\lambda_1(\hx)>\lambda_2(\hx)>\cdots>\lambda_k(\hx)$ and a
decomposition $\CC^d=E^1_\hx\oplus E_\hx^2\oplus\cdots\oplus
E^k_\hx$ into vector subspaces such that
$$
\hA(\hx)E_\hx^i=E_{\hf(\hx)}^i \quand \lambda_i(\hx)=
\lim_{|n|\to\infty} \frac 1n \log \|\hA^n(\hx)v\|
$$
for every non-zero $v\in E^i_\hx$ and $1\le i\le k$. We call $\dim
E^i_\hx$ the \emph{multiplicity} of $\lambda_i(\hx)$.

We assume that $\hmu$ is ergodic. Then the \emph{Lyapunov
exponents} $\lambda_i(\hx)$ are constant on a full measure subset
of $\hSigma$ and so are the dimensions of the Oseledets subspaces
$E_\hx^i$. The \emph{Lyapunov spectrum} of $\hA$ is the set of all
Lyapunov exponents. We say that the Lyapunov spectrum is
\emph{simple} if it contains exactly $d$ distinct values ($k=d$)
or, equivalently, if every Lyapunov exponent $\lambda_i$ has
multiplicity $1$. The main result in this paper, to be stated
below, provides an explicit sufficient condition for the Lyapunov
spectrum to be simple. We begin by describing the class of
cocycles to which it applies. In Appendix~\ref{s.AppA} we discuss
some extensions and applications.

\subsection{Symbolic dynamics}\label{sh.dynamics}

We take $\hSigma=\NN^{\,\ZZ}$, the full shift space with countably
many symbols, and $\hf:\hsig\to\hsig$ to be the shift map:
$$
\hf\big((x_n)_{n\in\ZZ}\big) = (x_{n+1})_{n\in\ZZ}.
$$
Let us call \emph{cylinder} of $\hsig$ any subset of the form
$$
[\iota_m, \ldots, \iota_{-1} ; \iota_0 ; \iota_1, \ldots, \iota_n] = \{\hx : x_j = \iota_j
\text{ for } j=m, \ldots, n\}.
$$
Cylinders of $\psig=\NN^{\{n\ge 0\}}$ and $\nsig=\NN^{\{n < 0\}}$
are defined similarly, corresponding to the cases $m=0$ and $n=-1$,
respectively, and they are represented as $[\iota_0, \iota_1,
\ldots, \iota_n]^u$ and $[\iota_m, \ldots, \iota_{-1}]^s$,
respectively. We endow $\hsig$, $\psig$, $\nsig$ with the topologies
generated by the corresponding cylinders. Let $\pP:\hsig\to\psig$
and $\nP:\hsig\to\nsig$ be the natural projections. We also consider
the one-sided shift maps $\pf:\psig\to\psig$ and $\nf:\nsig\to\nsig$
defined by
$$
\pf \circ \pP = \pP \circ \hf \quand \nf \circ \nP = \nP \circ
\hf^{-1}.
$$

For each  $\hx=(x_n)_{n\in\ZZ}$ in $\hsig$, we denote $\px=\pP(\hx)$
and $\nx=\nP(\hx)$. Then $\hx \mapsto (\nx,\px)$ is a homeomorphism
from $\hsig$ to the product $\nsig\times\psig$. In what follows we
often identify the two sets through this homeomorphism. When there
is no risk of ambiguity, we also identify the \emph{local stable set}
$$
 \Wsloc(\px) = \Wsloc(\hx) = \{(y_n)_{n\in\ZZ} : x_n=y_n \text{ for all } n\ge 0\}
 \quad\text{with}\quad\nsig
$$
and the \emph{local unstable set}
$$
 \Wuloc(\nx) = \Wuloc(\hx) = \{(y_n)_{n\in\ZZ} : x_n=y_n \text{ for all } n < 0\}
 \quad\text{with}\quad\psig,
$$
via the projections $\nP$ and $\pP$.

In Section~\ref{ss.inducing} we shall discuss how more general
situations may often be reduced to this one.

\subsection{Product structure}\label{sh.measure}

Let $\pmu=\pP_*\hmu$ and $\nmu=\nP_*\hmu$ be the images of the
ergodic $\hf$-invariant probability measure $\hmu$ under the
natural projections. It is easy to see that these are ergodic
invariant probabilities for $\pf$ and $\nf$, respectively. We take
$\nmu$ and $\pmu$ to be positive on cylinders. Moreover, we assume
$\hmu$ to be equivalent to their product, meaning there exists a
measurable function $\rho:\hsig\to(0,\infty)$ such that
$$
\hmu = \rho(\hx) \, (\nmu \times \pmu), \quad\hx\in\hsig.
$$

We assume that $\rho$ is bounded from zero and infinity. For
convenience of notation, we state this condition as follows: there
exists some constant $K>0$ such that
\begin{equation}\label{eq.boundedweak}
\frac 1K \le \frac{\rho(\nz,\px)}{\rho(\nz,\py)} \le K \quand
\frac 1K \le \frac{\rho(\nx,\pz)}{\rho(\ny,\pz)} \le K
\end{equation}
for all $\nx, \ny, \nz \in \nsig$ and $\px, \py, \pz \in \psig$.
Notice that $\{\hmu_\px=\rho(\cdot,\px)\nmu: \px\in\psig\}$
is a \emph{disintegration} of $\hmu$ into conditional probabilities
along local stable sets. By this we mean (see Rokhlin~\cite{Ro62}
or \cite[Appendix~C]{Beyond}) that $\hmu_\px(\Wsloc(\px))=1$ for
$\pmu$-almost every $\px$ and
$$
\hmu(D) = \int \hmu_x \big(D \cap \Wsloc(\px)\big) \, d\pmu(\px)
$$
for any measurable set $D\subset\hsig$. Analogously,
$\{\hmu_\nx=\rho(\nx,\cdot)\pmu: \nx\in\nsig\}$ is a
disintegration of $\hmu$ along local unstable sets. Since the
density $\rho$ is positive, the measures $\hmu_\px$, $\px\in\psig$
are all equivalent, and so are all $\hmu_\nx$, $\nx\in\nsig$.
Condition~\eqref{eq.boundedweak} just means that the Radon-Nikodym
derivatives
$$
 \frac{d\hmu_\px}{d\hmu_\py} \quad\text{with $\px, \py \in \psig$}
 \quand
 \frac{d\hmu_\nx}{d\hmu_\ny} \quad\text{with $\nx, \ny \in \nsig$}
$$
are uniformly bounded from zero and infinity. This will be used to
obtain the bounded distortion properties \eqref{eq.boundedjacnew}
and \eqref{eq.convX} below.

We also assume that the conditional probabilities $\hmu_\px$ and
$\hmu_\nx$ vary continuously with the base point, in the sense
that the functions
\begin{equation}\label{eq.continuousweak}
\psig\ni\px\mapsto\int\phi\,d\hmu_\px \quand
\nsig\ni\nx\mapsto\int\psi\,d\hmu_\nx
\end{equation}
are continuous, for any bounded measurable functions
$\phi : \nsig \to \RR$ and $\psi : \psig \to \RR$. Equivalently,
$$
 \px\mapsto
 \hmu_\px([\iota_m, \ldots, \iota_{-1}]^s)
 \quand
 \nx\mapsto\hmu_\px([\iota_0, \iota_1, \ldots, \iota_n]^u)
$$
are continuous for every choice of the $\iota_j$'s. This will be
used to obtain \eqref{eq.continuousjacnew} and
Lemma~\ref{l.continuity}.

In Section~\ref{ss.summable} we show that these hypotheses hold,
in particular, whenever the system satisfies a distortion
summability condition. Indeed, in that case the density $\rho$ may
be taken continuous and bounded from zero and infinity. In
general, the hypothesis~\eqref{eq.continuousweak} can probably be
avoided: that is the case at least when the cocycle is locally
constant; see the appendix of \cite{AV2} and also
Remark~\ref{r.avoidcontinuity} below.

\subsection{Invariant holonomies}\label{sh.cocycle}

Concerning the function $\hA:\hsig\to\GL$, we assume that it is
continuous and admits stable and unstable holonomies:

\begin{definition}\label{d.suholonomies}
We say $\hA$ admits \emph{stable holonomies}  if the limit
$$
\Hsxy = \lim_{n\to+\infty} \hA^n(\hy)^{-1}\hA^n(\hx)
$$
exists for any pair of points $\hx$ and $\hy$ in the same local
stable set, and depends continuously on $(\hx,\hy)$. \emph{Unstable
holonomies} $\Huxy$ are defined in a similar way, with $n\to-\infty$
and $\hx$ and $\hy$ in the same local unstable set.
\end{definition}
Notice that stable holonomies $\Hsxy:\CC^d\to\CC^d$ are linear maps
and they satisfy
\begin{itemize}
 \item[(a)] $\Hsxz = \Hszy \cdot \Hsxz$ and $\Hsxx = \id$,
 \item[(b)] $\hA(\hy) \cdot \Hsxy = \Hsxyf \cdot \hA(\hx)$,
\end{itemize}
over all points for which the relations make sense. Similar remarks
apply for the unstable holonomies.

For example, if $\hA$ is \emph{locally constant}, meaning that it is
constant on each cylinder $[\iota]$, $\iota\in\NN$, then
$\Hsxy\equiv\id$ and $\Huxy\equiv\id$.
In Section~\ref{ss.dominated} we discuss other situations where these
structures occur.

\subsection{Statement of main result}\label{ss.statement}
Let $\hp\in\hsig$ be a periodic point of $\hf$ and $q\ge 1$ be its
period. We call $\hz\in\hsig$ a \emph{homoclinic point} of $\hp$ if
$\hz\in \Wuloc(\hp)$ and there exists some multiple $l\ge 1$ of $q$
such that $\hf^l(\hz)\in \Wsloc(\hp)$. Then we define the
\emph{transition map}
$$
\psipz:\CC^d\to\CC^d, \quad \psipz = \Hspzf \cdot \hA^l(\hz)
\cdot \Hupz.
$$
The following notion is our main criterion for simplicity of the
Lyapunov spectrum. We refer to (p) as the \emph{pinching property}
and to (t) as the \emph{twisting property}.

\begin{definition}\label{d.simple}
We say that $\hA:\hsig\to\GL$ is \emph{simple} for $\hf$ if there
exists some periodic point $\hp\in\hsig$ of $\hf$ and some
homoclinic point $\hz\in\hsig$ of $\hp$ such that
\begin{itemize}
\item[(p)] All the eigenvalues of $\hA^q(\hp)$ have distinct
absolute values.

\item[(t)] For any invariant subspaces (sums of eigenspaces) $E$
and $F$ of $\hA^q(\hp)$ with $\dim E + \dim F = d$, we have
$
\psipz (E) \cap F = \{0\}.
$
\end{itemize}
\end{definition}

\begin{remark}
Let $\theta_j$, $j=1, \ldots, d$ represent the eigenspaces of
$\hA^q(\hp)$. For $d=2$ the twisting condition means that
$\psipz(\theta_i) \neq \theta_j$ for all $1\le i , j \le 2$. For
$d=3$ it means that $\psipz(\theta_i)$ is outside the plane
$\theta_j\oplus\theta_k$ and $\theta_i$ is outside the plane
$\psipz(\theta_j\oplus\theta_k)$, for all choices of $1\le i,  j,
k\le 3$. In general, this condition is equivalent to saying that
the matrix of the transition map in a basis of eigenvectors of
$\hA^q(\hp)$ has all its algebraic minors different from zero.
Indeed, it may be restated as saying that the determinant of the
square matrix
$$
\left(\begin{array}{cccccc} B_{1,i_1} & \cdots & B_{1,i_r} &
\delta_{1,j_1} & \cdots & \delta_{1,j_s} \\
\cdots & \cdots & \cdots &  \cdots & \cdots & \cdots \\
B_{d,i_1} & \cdots & B_{d,i_r} & \delta_{d,j_1} & \cdots &
\delta_{d,j_s}
\end{array}\right)
$$
is non-zero for any $I=\{i_1, \ldots, i_s\}$ and $J=\{j_1, \ldots,
j_r\}$ with $r+s=d$, where the $\delta_{i,j}$ are Dirac symbols
and the $B_{i,j}$ are the entries of the matrix of $\psipz$ in
the basis of eigenvectors. Up to sign, this determinant is the
algebraic minor $B[J^c\times I]$ corresponding to the lines
$j\notin J$ and columns $i\in I$.
\end{remark}

\begin{maintheorem}\label{main}
If $\hA:\hsig\to\GL$ is simple for $\hf$ then all the Lyapunov
exponents of the cocycle $\hFA$ for the measure $\hmu$ have
multiplicity~$1$.
\end{maintheorem}

Simplicity of the Lyapunov spectrum for independent random
matrices was investigated in the eighties by Guivarc'h,
Raugi~\cite{GR86}, and Gol'dsheid, Margulis~\cite{GM89}.
Theorem~\ref{main} also extends the main conclusions of Bonatti,
Viana~\cite{BoV04}, who treated the case when the base dynamics
$f$ is a subshift of finite type.

The present extension has been carried out to include in the
theory such examples as the Zorich cocycles, whose base dynamics
are not of finite type. It has been conjectured by Zorich and
Kontsevich~~\cite{KZ97,Zo96,Zo99} that the corresponding Lyapunov
exponents have multiplicity $1$. As an application of these ideas,
in \cite{AV2} we prove this conjecture. See also the comments in
Appendix~\ref{s.AppA} to the present paper.

Let us point out that we improve \cite{BoV04} not only in that
here we allow for infinite Markov structures and non-compact
ambient spaces, but also because our criterion is sharper: whereas
we only ask the cocycle to be simple, \cite{BoV04} needed a
similar hypothesis on all exterior powers as well.

\subsection{Outline of the proof}\label{ss.outline}

The starting point is the following observation. Let $\ell\in\{1,
\ldots, d-1\}$ be fixed and assume the cocycle has $\ell$ Lyapunov
exponents that are strictly larger than the remaining ones. Let
$E(\hx)$ be the sum of the Oseledets subspaces associated to those
largest exponents at a generic point $\hx\in\hsig$. Then
$\hx\mapsto E(\hx)$ defines a measurable invariant section of the
Grassmannian space of $\ell$-dimensional subspaces of $\CC^d$.
This section is invariant along local unstable sets, meaning that
$$
E(\hy) =
 \Huxy\cdot E(\hx) \quad\text{for all } \hy \in \Wuloc(\hx),
$$
because the hypotheses in Section~\ref{sh.cocycle} imply that
$$
\hA^n(\hy) = \Huxyfn \cdot \hA^n(\hx) \cdot \Huyx
 \quad\text{for all } n < 0,
$$
and the norms of the unstable holonomies are bounded. Let $\hm$ be
the probability measure on $\hsig\times\Grass$ which projects down
to $\hmu$ and has the Dirac measures $\delta_{E(\hx)}$ as
conditional probabilities along the Grassmannian fibers. Then
$\hm$ is an invariant measure for the action of $\hA$ on the
Grassmannian bundle $\hsig\times\Grass$ and, typically, it is the
unique one whose conditional probabilities are invariant under
unstable holonomies.

To try and prove the theorem, we consider the space of all
probability measures $\hm$ on $\hsig\times\Grass$ that project down
to $\hmu$, are invariant under the action of the cocycle, and whose
conditional probabilities $\hm_\hx$ along the Grassmannian fibers
are invariant under unstable holonomies.
Proposition~\ref{p.existence} ensures that such \emph{invariant
$u$-states} do exist. In Proposition~\ref{p.continuity} we prove
that the projection $\pom$ of any $u$-state $\hm$ to
$\psig\times\Grass$ admits conditional probabilities $\pom_\px$
along the Grassmannian fibers that depend continuously on the base
point $\px$. This is very important for our arguments: continuity
allows us to show that the kind of behavior the cocycle exhibits on
the periodic point $\hp$ in Definition~\ref{d.simple} propagates to
almost all orbits on the whole $\hsig$. Let us explain this.

Firstly, in Proposition~\ref{p.convergencia1}, we use a simple
martingale argument to show that the measure $\hm$ may be
recovered from $\pom$ through
\begin{equation}\label{eq.lim1}
\hm_\hx = \lim_{n\to\infty} \hA^n(\hf^{-n}(\hx))_* \,
\pom_{\pP(\hf^{-n}(\hx))} \qquad\text{$\hmu$-almost everywhere.}
\end{equation}
The assumption that $\hA^q(\hp)$ has $\ell$ largest eigenvalues
implies that $\hA^{qn}(\hp)_*\eta$ converges to the Dirac measure on
the sum of the eigenspaces associated to the largest eigenvalues,
for any probability measure $\eta$ on $\Grass$ that gives zero
weight to the hyperplane section defined by the other invariant
subspaces. A crucial step, carried out in Section~\ref{s.Dirac}, is
to prove that \emph{the limit on the right hand side of
\eqref{eq.lim1} is a Dirac measure for almost every $\hx$.} The
proof has two main parts. In Proposition~\ref{p.noatoms} we use the
assumption that the cocycle is simple to show that the conditional
probabilities of $m$ give zero weight to hyperplane sections of the
Grassmannian. Then, in Proposition~\ref{p.toDirac}, we use the
continuity property in the previous paragraph, and the assumption
that the cocycle is simple, to show that the behavior on the
periodic point we just described does propagate to almost every
orbit.

This proves that $\hm_\hx = \delta_{\xi(\hx)}$ almost everywhere,
where $\xi(\hx)$ is some $\ell$-subspace. In view of what we wrote
before, $\xi(\hx)$ should correspond to the subspace $E(\hx)$
associated to the largest Lyapunov exponents. To prove that this
is indeed so, we must also find the complementary invariant
subspace.
This is done by applying the previous theory to the adjoint
(relative to some Hermitian form) cocycle $\hB=\hA^*$ over the
inverse map $\hf^{-1}$. Since our hypotheses are symmetric under
time reversion, the same arguments as before yield an
$\ell$-dimensional section $\hx\mapsto\xi^*(\hx)$ which is
invariant under the action of $\hB$ and under stable holonomies.

Let $\eta(\hx)$ be the orthogonal complement of $\xi^*(\hx)$. Then
$\xi$ and $\eta$ are $\hA$-invariant sections with complementary
dimensions. Using the simplicity assumption once more, we check that
$\xi(\hx)$ and $\eta(\hx)$ are transverse to each other at almost
every point. The final step is to deduce from \eqref{eq.lim1} that
the Lyapunov exponents of $\hA$ along $\xi$ are strictly larger than
those along $\eta$.

\medskip\noindent
{\bf Acknowledgments.} We are grateful to A. Arbieto, C. Matheus,
and J.-C. Yoccoz, for several useful conversations, and to E.
Esteves for explanations on the structure of the Grassmannian
manifolds.

\section{Preliminary observations}

Here we recall a few basic notions and prove a number of technical
facts that will be useful in the sequel. The reader may be well
advised to skip this section in a first reading, and then come
back to it when a specific result or concept is needed.

\subsection{Exterior powers and Grassmannians}\label{ss.exterior}

Fix any $\ell\in\{1, \ldots, d-1\}$. The $\ell$th exterior power of
$\CC^d$, denoted by $\extl{\CC^d}$, is the vector space of alternate
$\ell$-forms $\omega:(\CC^d)^*\times \cdots \times (\CC^d)^*\to\CC$
on the dual space $(\CC^d)^*$. It has
$$
\dim \extl{\CC^d} = \left(\begin{array}{c} d \\ \ell
\end{array}\right).
$$
Every element of $\extl{\CC^d}$ may be written as a sum of elements
of the form $\omega_1\wedge\cdots\wedge\omega_\ell$ with
$\omega_i\in (\CC^d)^{**}$. We represent by $\extlv{\CC^d}$ the
subset of elements of this latter form, that we call
\emph{$\ell$-vectors}. Any $\ell$-vector may be written as
$c\,w_1\wedge\cdots\wedge w_\ell$, where $c\in\CC$ and the $w_i$ are
orthogonal unit vectors (relative to any fixed Hermitian form).
Hence, $\extlv{\CC^d}$ is a closed subset of $\extl{\CC^d}$.

Since the bi-dual space is canonically isomorphic to $\CC^d$, we may
think of the $\omega_i$ as vectors in $\CC^d$. Thus, there is a
natural projection $\pi_v$ from $\extlv{\CC^d}\setminus\{0\}$ to the
Grassmannian $\Grass$ of $\ell$-dimensional subspaces of $\CC^d$,
associating to each non-zero $\ell$-vector
$\omega_1\wedge\cdots\wedge\omega_\ell$ the subspace generated by
$\{\omega_1, \dots, \omega_\ell\}$. Two $\ell$-vectors have the same
image under $\pi_v$ if and only if one is a multiple of the other.
In other words, $\pi_v$ induces a bijection between $\Grass$ and the
projective space $\pextlv{\CC^d}$ of the space of $\ell$-vectors.

The $\ell$th exterior power $\extl{B}:\extl{\CC^d}\to \extl{\CC^d}$
of an operator $B:\CC^d\to\CC^d$ is defined by
$$
\extl{B}(\omega)(\phi_1,\ldots,\phi_\ell) = \omega\big(\phi_1\circ
B,\ldots,\phi_\ell\circ B\big).
$$
Notice that $\extl{B}(\omega_1\wedge\cdots\wedge\omega_\ell) =
B(\omega_1)\wedge\cdots\wedge B(\omega_\ell)$, and so $\extl{B}$
preserves the set $\extlv{\CC^d}$ of $\ell$-vectors. Moreover,
assuming $B$ is invertible,
\begin{equation}\label{eq.downtograss}
\pi_v\circ \extl{B} = B_\# \circ\pi_v
 \quad\text{on $\extlv{\CC^d}$},
\end{equation}
where $B_\#$ denotes the action of $B$ on the Grassmannian.

Let $H$ be a hyperplane, that is, a codimension $1$ linear subspace
of the vector space $\extl{\CC^d}$. Then $H$ may be written as
$$
H=\{\omega\in \extl{\CC^d} : \omega\wedge\upsilon=0\}
$$
for some non-zero $\upsilon\in \extdl{\CC^d}$. We call the
hyperplane \emph{geometric} if $\upsilon$ may be chosen a
$(d-\ell)$-vector, that is,
$\upsilon=\upsilon_{\ell+1}\wedge\cdots\wedge\upsilon_d$ for some
choice of vectors $\upsilon_i$ in $\CC^d=(\CC^d)^{**}$. By
definition, a \emph{hyperplane section} of $\Grass$ is the image
under the projection $\pi_v$ of the intersection of $\extlv{\CC^d}$
with some geometric hyperplane $H$ of $\extl{\CC^d}$. Note that,
given any $\ell$-vector
$\omega=\omega_1\wedge\cdots\wedge\omega_\ell$,
$$
\omega\in H \Leftrightarrow \omega \wedge \upsilon =0
\Leftrightarrow \pi_v(\omega) \cap \pi_v(\upsilon) \neq \{0\}.
$$
Hence, the hyperplane section of $\Grass$ associated to $H$
contains precisely the $\ell$-dimensional subspaces that have
non-trivial intersection with the $(d-\ell)$-dimensional subspace
generated by $\upsilon$. The \emph{orthogonal hyperplane section}
to $V \in \Grass$ is the hyperplane section associated to its
orthogonal complement $V^\perp$.

To any Hermitian form on $\CC^d$ there is a canonically associated
one on $\extl{\CC^d}$ such that the set of $\ell$-vectors
$e_{i_1}\wedge\cdots\wedge e_{i_\ell}$, $1\le i_1 < \cdots < i_\ell
\le d$ obtained from an arbitrary orthonormal basis $e_1, \ldots,
e_d$ of the space $E$ is an orthonormal basis of its exterior power.
If $B$ is a unitary operator then so is $\extl{B}$. Let $e_1,
\ldots, e_d$ be an orthonormal basis of $\CC^d$. We use the
\emph{polar decomposition} $B=K' D K$ of a linear isomorphism
$B:\CC^d\to\CC^d$, where $K$ and $K'$ are unitary operators, and $D$
is a diagonal operator (with respect to the chosen basis) with
positive eigenvalues $a_1, \ldots, a_d$. The $a_i$ are called
\emph{singular values} of $B$; we always take them to be numbered in
non-increasing order.

\subsection{Eccentricity of linear maps}\label{ss.eccentricity}

Let $L:\CC^d \to \CC^d$ be a linear isomorphism and $1\le \ell \le
d$. The \emph{$\ell$-dimensional eccentricity} of $L$ is defined
by
$$
\cE(\ell,L) = \sup\big\{\frac{m(L\mid\xi)}{\|L\mid\xi^\perp\|} :
 \xi\in\Grass\big\}, \qquad m(L\mid\xi) = \|(L\mid\xi)^{-1}\|^{-1}.
$$
We call \emph{most expanded $\ell$-subspace} any $\xi\in\Grass$
that realizes the supremum. These always exist, since the
Grassmannian is compact and the expression depends continuously on
$\xi$. These notions may be expressed in terms of the polar
decomposition of $L$ with respect to any orthonormal basis:
denoting by $a_1, \ldots, a_d$ the eigenvalues of the diagonal
operator $D$, in non-increasing order, then
$\cE(\ell,L)=a_{\ell}/a_{\ell+1}$. The supremum is realized by any
subspace $\xi$ whose image under $K$ is a sum of $\ell$
eigenspaces of $D$ such that the product of the eigenvalues is
$a_1\cdots a_\ell$. It follows that $\cE(\ell,L)\ge 1$, and the
most expanded $\ell$-subspace is unique if and only if the
eccentricity is larger than $1$.

Let $e_1, \ldots, e_d$ be a basis of eigenvectors of $D$
corresponding to the eigenvalues $a_1, \ldots, a_d$. For any
$I\subset\{1, \ldots, d\}$ we represent $E_I=\oplus_{i\in I} e_i$.
Given any $\eta\in\Grass$ one may find a subset $I=\{i_1, \ldots,
i_\ell\}$ of $\{1, \ldots, d\}$ such that $\eta$ is the graph of a
linear map
$$
E_I \to E_J, \quad e_i \mapsto \sum_{j\in J} \eta(i,j)\, e_j,
$$
where $J$ is the complement of $I$. We say that $\eta'\in\Grass$
is in the $\vep$-neighbor\-hood $B_\vep(\eta)$ of $\eta$ if (for
some choice of $I$) it may also be written as the graph of a
linear map from $E_I$ to $E_J$ such that all corresponding
coefficients $\eta(i,j)$ and $\eta'(i,j)$ differ by less than
$\vep$. Given a hyperplane section $H$ of $\Grass$, defined by
some $(d-\ell)$-vector $\upsilon$, and given $\delta>0$, we
represent by $H_\delta$ the union of the hyperplane sections
defined by all the $(d-\ell)$-vectors in the $B_\delta(\eta)$.

\begin{lemma}\label{l.vizhyper}
Given $C\ge 1$ and $\delta>0$ there exists $\vep>0$ such that, for
any $\eta\in\Grass$ and any diagonal operator $D$ with
eccentricity $\cE(\ell,D)\le C$, one may find a hyperplane section
$H$ of $\Grass$ such that $D^{-1}(B_\vep(\eta)) \subset H_\delta$.
\end{lemma}

\begin{proof}
Choose $I=\{i_1, \ldots, i_\ell\}$ such that $\eta$ is a graph
over the subspace generated by $e_{i_1}, \cdots, e_{i_\ell}$. In
other words, $\eta$ admits a basis of the form
$$
\{e_i + \sum_{j\in J} \eta(i,j) e_j : i \in I\},
$$
where $J=\{j_1, \ldots, j_{\ell-d}\}$ is the complement of $I$
inside $\{1, \ldots, d\}$. Let $a_1, \ldots, a_d$ be the
eigenvalues of $D$, in non-increasing order. Then
$$
\{f_i = e_i + \sum_{j\in J} \frac{a_i}{a_j} \, \eta(i,j) \, e_j :
i \in I\},
$$
is a basis of $D^{-1}(\eta)$. We claim that there exist $\alpha\in
I$ and $\beta\in J$ such that $a_\alpha/a_\beta\le K$: if $I=\{1,
\ldots, \ell\}$ it suffices to take $\alpha=\ell$ and
$\beta=\ell+1$; otherwise, we may always choose $\beta\in\{1,
\ldots, \ell\}\setminus I$ and $\alpha \in I\setminus\{1, \ldots,
\ell\}$, and then we even have $a_\alpha/a_\beta\le 1$. This
proves the claim. Now let
$$
\upsilon = e_{j_1} \wedge \cdots \wedge e_{\alpha,\beta} \wedge
 \cdots \wedge e_{j_\ell},
  \quad e_{\alpha,\beta}=e_\alpha\pm \frac{a_\alpha}{a_\beta}\,\eta(\alpha, \beta) \, e_\beta
$$
be the $(d-\ell)$-vector given by the wedge products of all $e_j$,
$j\in J$ except that $e_\beta$ is replaced by $e_{\alpha,\beta}$.
Notice that
\begin{equation*}\begin{aligned}
D^{-1}(\eta)\wedge\upsilon
 = & f_{i_1} \wedge \cdots \wedge f_{i_\ell} \wedge e_{j_1} \wedge
   \cdots \wedge e_{\alpha,\beta} \wedge \cdots \wedge e_{j_\ell} \\
 = & \big[ e_{i_1} \wedge \cdots \wedge e_{i_\ell}\big]\,
     \wedge \big[e_{j_1} \wedge \cdots \wedge \pm (a_\alpha/a_\beta) \, \eta(\alpha, \beta)
\, e_\beta \wedge \cdots \wedge e_{j_\ell} \big] \\
+ \big[e_{i_1} \wedge \cdots & \wedge (a_\alpha/a_\beta) \,
\eta(\alpha,\beta) \, e_\beta \wedge \cdots \wedge e_{i_\ell}\big]
\wedge \big[e_{j_1} \wedge \cdots \wedge e_\alpha \wedge\cdots
\wedge e_{j_\ell}\big]\,.
\end{aligned}\end{equation*}
Choosing the sign $\pm$ appropriately, the two terms cancel out
and so $D^{-1}(\eta)\wedge\upsilon=0$. This means that
$D^{-1}(\eta)$ belongs to the hyperplane section $H$ defined by
$\upsilon$. In just the same way, given any $\eta'$ in the
$\vep$-neighborhood of $\eta$ we may find a $(d-\ell)$-vector
$$
\upsilon' = e_{j_1} \wedge \cdots \wedge e'_{\alpha,\beta} \wedge
 \cdots \wedge e_{j_\ell},
  \quad e'_{\alpha,\beta} = e_\alpha\pm \frac{a_\alpha}{a_\beta}\,\eta'(\alpha, \beta) \, e_\beta
$$
such that $D^{-1}(\eta')$ belongs to the hyperplane section
defined by $\upsilon'$. Since $a_\alpha/a_\beta\le K$ and
$|\eta(\alpha,\beta)-\eta'(\alpha,\beta)|<\vep$, we have that
$\upsilon'\in B_\delta(\upsilon)$ as long as $\vep$ is small
enough. Then $D^{-1}(\eta')\in H_\delta$ for all $\eta'$ in the
$\vep$-neighborhood of $\eta$, as claimed.
\end{proof}

\begin{proposition}\label{p.eccentricity}
Let $\cN$ be a \weak compact family of probabilities on $\Grass$
such that all $\nu\in\cN$ give zero weight to all hyperplane
sections. Let $L_n:\CC^d\to\CC^d$ be linear isomorphisms such that
$(L_n)_* \nu_n $ converges to a Dirac measure $\delta_\xi$ as
$n\to\infty$, for some sequence $\nu_n$ in $\cN$. Then the
eccentricity $\cE(\ell,L_n)$ goes to infinity and the image
$L_n(\zeta^a_n)$ of the most expanding $\ell$-subspace of $L_n$
converges to $\xi$.
\end{proposition}

\begin{proof}
Let $L_n:\CC^d\to\CC^d$, $\nu_n\in\cN$, and $\xi\in\Grass$ be as
in the statement. Consider the polar decomposition $L_n=K_n' D_n
K_n$, where $D_n$ has eigenvalues $a_1, \ldots, a_d$, in
non-increasing order.

We begin by reducing to the case $K_n=K_n'=\id$. Let
$\cM=\UU(\ell,d)_*\cN$, where $\UU(\ell,d)$ is the group of
transformations induced on $\Grass$ by the unitary group. It is
clear that all $\mu\in\cM$ give zero weight to every hyperplane
section of $\Grass$. Notice also that $\cM$ is \weak compact: given
any sequence $\mu_j=(U_j)_*\nu_j$ with $\nu_j\in\cN$ and
$U_j\in\UU(\ell,d)$, up to considering subsequences one may assume
that $\nu_j$ converges to some $\nu\in\cN$ in the \weak topology and
$U_j$ converges to some $U\in\UU(\ell,d)$ uniformly on $\Grass$, and
then $(U_j)_*\nu_j$ converges to $U_*\nu\in\cM$ in the \weak
topology.  Let $\mu_n=(K_n)_*\nu_n\in\cM$. Then $(K_n'D_n)_*\mu_n$
converges to $\delta_\xi$.  In addition, up to considering a
subsequence, we may assume that $K_n'$ converges to some
$K'\in\UU(\ell,d)$ uniformly on $\Grass$. Note that
$((K')^{-1}K_n'D_n)_*\mu_n$ converges to $\delta_\eta$, where
$\eta=(K')^{-1}(\xi)$. Since $(K')^{-1}K_n'$ converges uniformly to
the identity, this implies that $(D_n)_*\mu_n$ also converges to the
Dirac measure at $\eta$.

Now, since $\cM$ and the space of hyperplane sections of $\Grass$
are compact, we may find $\delta>0$ such that $\nu(H_\delta) <
1/2$ for every $\mu\in\cN$ and every hyperplane section $H$ of
$\Grass$. On the other hand, given any $\vep>0$ we have
$$\mu_n(D_n^{-1}(B_\vep(\eta)))= (D_n)_*\mu_n(B_\vep(\eta))> 1/2$$
for every large $n$. Then $D_n^{-1}(B_\vep(\eta))$ can not
contained in $H_\delta$, for any hyperplane section $H$. In view
of Lemma~\ref{l.vizhyper}, this implies that
$\cE(\ell,L_n)=\cE(\ell,D_n)$ goes to infinity as $n\to\infty$, as
claimed in the first part of the lemma.

The second part is a consequence, through similar arguments. Given
any $\vep>0$, fix $\delta>0$ small enough so that
$\nu(H_\delta)<\vep$ for any $\nu\in\cN$ and any hyperplane section
$H$ of $\Grass$. Let $H^n\subset\Grass$ be the hyperplane section
orthogonal to the most expanding direction $\zeta^a_n$ of $L_n$. By
definition, the complement $\Grass\setminus H^n_\delta$ of the
$\delta$-neighborhood of $H^n$ consists of the elements of $\Grass$
that avoid any $(d-\ell)$-subspace $\delta$-close to
$(\zeta^a_n)^\perp$. Since the eccentricity of $L_n$ goes to
infinity,
$$
L_n\big(\Grass\setminus H^n_\delta\big)
 \subset B_\vep(L_n(\zeta^a_n))
$$
for every large $n$. Then, the $(L_n)_*\nu_n$-measure of
$B_\vep(L_n(\zeta^a_n))$ is larger than $1-\vep$. Since
$(L_n)_*\nu_n$ converges to the Dirac measure at $\xi$, it follows
that $\xi\in B_\vep(L_n(\zeta^a_n))$ for every large $n$. As
$\vep>0$ is arbitrary, this proves the second claim in the
proposition.
\end{proof}

\subsection{Quasi-projective maps}\label{ss.quasiprojective}

Let $v\mapsto [v]$ be the canonical projection from $\CC^d$ minus
the origin to the projective space $\PC$. We call $P_\#:\PC\to\PC$
a \emph{projective map} if there is some $P\in \GL$ that induces
$P_\#$ through $P_\#([v]) = [P(v)]$. It was pointed out by
Furstenberg~\cite{Fu73} that the space of projective maps has a
natural compactification, the space of \emph{quasi-projective
maps}, defined as follows. The quasi-projective map $Q_\#$ induced
by a non-zero, possibly non-invertible, linear map $Q:\CC^d \to
\CC^d$ is given by $Q_\#([v]) = [Q(v_1)]$ where $v_1$ is any
vector such that $v-v_1$ is in $\ker Q$. Observe that $Q_\#$ is
defined and continuous on the complement of the projective
subspace $\ker Q_\# = \{[v]:v\in\ker Q\}$. The space of
quasi-projective maps inherits a topology from the space of
non-zero linear maps, through the natural projection $Q\mapsto
Q_\#$. Clearly, every quasi-projective map $Q_\#$ is induced by
some linear map $Q$ such that $\|Q\| = 1$. It follows that the
space of quasi-projective maps in $\PC$ is compact for this
topology.

This notion has been extended to transformations on Grassmannian
manifolds, by Gol'dsheid, Margulis~\cite{GM89}. Namely, one calls
$P_\#:\Grass\to\Grass$ a \emph{projective map} if there is $P\in
\GL$ that induces $P_\#$ through $P_\#(\xi) = P(\xi)$. Note that
$P$ may always be taken such that the map $\extl{P}$ it induces on
$\extl{\CC^d}$ has norm $1$. Let $\cQ$ be the closure of the set
of all transformations $\extl{P}$ with $P$ invertible. Since every
$\extl{P}$ preserves the closed subset $\extlv{\CC^d}$, so does
every $Q\in\cQ$. The \emph{quasi-projective} map $Q_\#$ induced on
$\Grass$ by a map $Q\in\cQ$ is given by
$Q_\#(\pi_v(\omega))=\pi_v(Q(\omega))$ for any $\ell$-vector
$\omega$ in the complement of $\ker Q$. The space of all
quasi-projective maps on $\Grass$ inherits a topology from $\cQ$,
through the natural projection $Q\mapsto Q_\#$, and it is compact
for this topology, since we may always take $Q$ with norm equal to
$1$.

\begin{lemma}\label{l.kernel}
The kernel $\ker Q_\# = \pi_v(\ker Q)$ of any quasi-projective map
is contained in some hyperplane section of $\Grass$.
\end{lemma}

\begin{proof}
We only have to check that $\ker Q$ is contained in a geometric
hyperplane of $\extl{\CC^d}$. Let $P_n$ be any sequence of linear
invertible maps such that every $\extl{P_n}$ has norm $1$ and they
converge to $Q$. Consider the polar decomposition $P_n=K_n' D_n
K_n$ where $D_n=\diag[a_1^n, \ldots, a_d^n]$ relative to some
orthonormal basis $e_1, \ldots, e_d$. Then $\extl{P_n}=
\extl{K_n'}\extl{D_n}\extl{K_n}$ is the polar decomposition of
$\extl{P_n}$, where $\extl{D_n}$ is diagonal relative to the basis
$e_{i_1}\wedge\cdots\wedge e_{i_\ell}$, $i_1 < \cdots < i_\ell$ of
$\extl{\CC^d}$. Denote $e=e_1\wedge\cdots\wedge e_\ell$. Since the
eigenvalues $a_i^n$, $i=1, \ldots, d$ are in non-increasing order,
$$
a_1^n \cdots a_\ell^n = \|\extl{D_n}(e)\| = \|\extl{D_n}\| =
\|\extl{P_n}\| = 1.
$$
Taking the limit over a convenient subsequence, we get that
$Q=\extl{K'}\cD\extl{K}$ for some unitary operators $K$, $K'$ and
some norm $1$ operator $\cD$ diagonal with respect to the basis
$e_{i_1}\wedge\cdots\wedge e_{i_\ell}$. Moreover, $\|\cD(e)\|=1$
and the kernel of $\cD$ is contained in the hyperplane section
$H(e)$ orthogonal to $e$. Let $\omega=\extl{K}^{-1}(e)$ and
$H=\extl{K}^{-1}(H(e))$ be the hyperplane section orthogonal to
$\omega$. Then
$$
\eta\in\ker Q \Leftrightarrow \extl{K}\eta\in\ker\cD \Rightarrow
\extl{K}\eta\in H(e) \Leftrightarrow \eta\in H,
$$
and this proves the statement.
\end{proof}

The \emph{\weak topology} in the space of probability measures on $\Grass$
is characterized by the property that a sequence $(\nu_n)_n$ converges to a
probability $\nu$ if and only if, given any continuous function
$g:\Grass\to\RR$, the integrals $\int g\,d\nu_n$ converge to $\int g\,d\nu$.
It is well-known that this topology is metrizable and compact, because the
space of continuous functions on the Grassmannian contains countable dense
subsets.

\begin{lemma}\label{l.topology2}
If $(P_{n})_n$ is a sequence of projective maps converging to some
quasi-projective map $Q$ of $\Grass$, and $(\nu_n)_n$ is a
sequence of probability measures in $\Grass$ converging weakly to
some probability $\nu$ with $\nu(\ker Q)=0$, then
$(P_{n})_*\nu_n$ converges weakly to $Q_*\nu$.
\end{lemma}

\begin{proof}
Let $(K_m)_m$ be a basis of neighborhoods of $\ker Q$ such that
$\nu(\partial K_m)=0$ for all $m$. Given any continuous
$g:\Grass\to\RR$, and given $\vep>0$, fix $m\ge 1$ large enough
so that $\nu(K_m)\le\vep$. Then fix $n_0\ge m$ so that
$\nu_n(K_m)\le \nu(K_m)+\vep \le 2\vep$,
$$
\Big|\int_{K_m^{\,c}} (g\circ Q)\,d\nu_n - \int_{K_m^{\,c}}
(g\circ Q)\,d\nu\Big| \le \vep \quad\text{and}\quad
\sup_{K_m^{\,c}} \Big|g\circ P_n - g\circ Q\Big| \le
\vep
$$
for all $n\ge n_0$. Then, splitting into integrals over $K_m$ and
over $K_m^c$,
$$
\big|\int (g\circ P_n)\,d\nu_n - \int (g\circ
Q)\,d\nu \big|\le 2\vep + 3 \vep \sup |g|
$$
for all $n\ge n_0$. This proves the lemma.
\end{proof}

For notational simplicity, in what follows we drop the subscript
$\#$ and use the same symbol to represent a linear transformation
and its action on any of the spaces $\Grass$, $0 < \ell < d$. In
particular, we also denote by $\hFA$ the Grassmannian cocycles
$\hsig\times\Grass\to\hsig\times\Grass$ defined by $\hA$ over $\hf$.

\subsection{Bounded distortion}\label{ss.jacobians}

Let $k\ge 1$ be fixed. For each $I=(\iota_0, \ldots, \iota_{k-1})$
denote by $\pfkI:\psig \to \pI$ the inverse branch of
$\pfk=(\pf)^k$ with values in the cylinder $\pI=[\iota_0, \ldots,
\iota_{k-1}]^u$. Moreover, define
\begin{equation}\label{eq.defjac1}
\pJmufkI(\px) = \hmu_\px([I]^s) \quad\text{for each $\px\in\psig$,}
\end{equation}
where $[I]^s=[\iota_0, \ldots \iota_{k-1}]^s$.
The boundedness condition \eqref{eq.boundedweak} gives
\begin{equation}\label{eq.boundedjacnew}
\frac{1}{K}
 \le \frac{\pJmufkI(\px)}
         {\pJmufkI(\py)} \le K
\end{equation}
for every $I$ and any pair of points $\px$ and $\py$ in $\psig$.
This will be used in the proof of Lemma~\ref{l.todos1} and
Corollary~\ref{c.invcont1}.
Moreover, the continuity condition \eqref{eq.continuousweak} implies
that the function
\begin{equation}\label{eq.continuousjacnew}
\px \mapsto \pJmufkI(\px)
\end{equation}
is continuous on $\psig$, for every choice of $I$. In both cases,
we also have dual objects and statements for inverse branches
$\nfkI$ of the iterates of $\nf$. From \eqref{eq.continuousweak}
we also get the following fact, which will be used in the proof of
Proposition~\ref{p.continuity}.

\begin{lemma}\label{l.continuity}
Let $\Phi:\hsig\to\RR$ be a bounded measurable function such that,
for every fixed $\nx\in\nsig$, the function $\px\mapsto\Phi(\nx,\px)$
is continuous at some $\pz\in\psig$. Then
$$
\px \mapsto \int \Phi(\nx,\px) \,d\hmu_\px(\nx)
 \quad\text{is continuous at $\pz$}.
$$
There is also a dual statement obtained by interchanging the roles
of $\nx$ and $\px$.
\end{lemma}

\begin{proof}
Let $\pz\in\psig$ and $\vep>0$ be fixed. Define $\phi(\nx)=\Phi(\nx,\pz)$
for every $\nx\in\nsig$. The continuity condition~\eqref{eq.continuousweak}
gives that
\begin{equation}\label{eq.lcont1}
|\int \phi(\nx)\,d\hmu_{\px}(\nx) - \int
\phi(\nx)\,d\hmu_{\pz}(\nx)| < \vep
\end{equation}
for any $\px$ in some neighborhood $Z_0$ of the point $\pz$. Let
$Z_n$, $n\ge 0$ be a decreasing basis of neighborhoods of $\pz$. The
assumption that $\Phi$ is continuous on the second variable means
that for every $\nx$ there exists some $n\ge 1$ such that
$$
|\Phi(\nx,\px)-\phi(\nx)| < \vep \quad \text{for all $\px\in Z_n$.}
$$
Let $V(k,\vep)\subset\nsig$ be the set of points $\nx\in\nsig$ for which
we may take $n\le k$. Consider $k$ large enough so that the
$\hmu_\pz$-measure of $V(k,\vep)^c$ is less than $\vep$. Then, using
condition \eqref{eq.boundedweak},
$$
 \hmu_\px(V(k,\vep)^c) < K \vep
 \quad\text{for every $\px\in\psig$.}
$$
The difference $|\int \Phi(\nx,\px)\,d\hmu_{\px}(\nx) - \int
\phi(\nx)\,d\hmu_{\px}(\nx)|$ is bounded above by
$$
\int_{V(k,\vep)} \big| \Phi(\nx,\px) - \phi(\nx) \big |
\,d\hmu_{\px}(\nx)
 + 2 \sup|\Phi| \, \hmu_{\px}(V(k,\vep)^c)
 $$
and so, for any $\px \in Z_k$,
\begin{equation} \label{eq.lcont2}
 |\int \Phi(\nx,\px)\,d\hmu_{\px}(\nx) - \int \phi(\nx)\,d\hmu_{\px}(\nx)|
 < \vep + 2 K \vep \sup|\Phi| \,.
\end{equation}
Putting \eqref{eq.lcont1} and \eqref{eq.lcont2} together, we
conclude that
$$
|\int \Phi(\nx,\px)\,d\hmu_{\px}(\nx) - \int
\Phi(\nx,\pz)\,d\hmu_{\pz}(\nx)| < 2\vep + 2 K \vep \sup|\Phi|
$$
for every $\px$ in the neighborhood $Z_k$ of $\pz$. This proves the
lemma.
\end{proof}

Given any measurable set $F\subset\psig$ and any $I=(\iota_0,
\ldots, \iota_{k-1})$, we have
$$
\hf^{-k}([I]^s \times F) = \nsig \times \pfkI(F) = (\pP)^{-1}(\pfkI(F)).
$$
Consequently, since $\hmu$ is invariant under $\hf$ and
$\pmu=\pP_*\hmu$,
$$
\int_F \JmufkI(\px)\, d\pmu(\px)
 = \int_F \hmu_\px ([I]^s) \, d\pmu(\px)
 = \hmu([I]^s \times F)
 = \pmu(\pfkI(F)).
$$
Thus, $\pJmufkI$ is a Radon-Nikodym derivative of the measure
$F\mapsto\pmu\big(\pfkI(F)\big)$ with respect to $\pmu$. An
equivalent formulation is
$$
\int (\psi \cdot \pJmufkI) \, d\pmu = \int_{\pI} (\psi \circ \pfk)
\, d\pmu.
$$
for any bounded measurable function $\psi:\psig\to\RR$, the previous
equality corresponding to the case $\psi=\Chi_F$.
Considering $F=\{\px\in\psig: \pJmufkI(\px)=0\}$, we get that
$\pJmufkI(\pfk(\pz))>0$ for $\pmu$-almost every $\pz \in \pI$.
Therefore,
\begin{equation}\label{eq.defjac2}
\pJmufk:\psig \to (0,+\infty), \quad \pJmufk(\pz)
 = \frac{1}{\pJmufkI(\pfk(\pz))} \text{\ \ when $\pz \in \pI$}
\end{equation}
is well defined $\pmu$-almost everywhere. Moreover, given any
bounded measurable function $\xi:\pI\to\RR$ and denoting
$\psi=(\xi \cdot \pJmufk)\circ \pfkI$, we have that
$$
 \int(\xi\circ \pfkI)\,d\pmu
 = \int(\psi \cdot \pJmufkI)\,d\pmu
 = \int_{\pI}(\psi \circ \pfk)\,d\pmu
 = \int(\xi \cdot \pJmufk)\,d\pmu.
 $$
In particular, taking $\xi=\Chi_B$,
$$
\mu(\pfk(B)) = \int_B \pJmufk \, d\pmu
 \quad\text{for every measurable $B\subset \pI$.}
$$
In other words, $\pJmufk$ is a \emph{Jacobian} of $\pmu$ for the
$k$th iterate of $\pf$.

\begin{lemma}\label{l.todos1}
Given any $I=(\iota_0, \ldots, \iota_{k-1})$ and any
$\pz\in [I]^u$,
$$
 \hf_*^k \hmu_\pz = \pJmufk(\pz)\big( \hmu_{\pfk(\pz)} \mid [I]^s).
$$
Moreover, a dual statement is true for $\hf_*^{-k} \hmu_\nz$.
\end{lemma}

\begin{proof}
Let $\px=\pfk(\pz)$. Clearly, $\pz=\pfkI(\px)$ and $\hf^k$ maps
$\Wsloc(\pz)$ bijectively to $[I]^s \times \{\px\} \subset
\Wsloc(\px)$. Consider any $J=(\iota_l, \ldots, \iota_{-1})$,
where $l < 0$, and denote $JI=(\iota_l, \ldots, \iota_0, \ldots,
\iota_{k-1})$. By the definition \eqref{eq.defjac1},
$$
\hmu_{\px}([JI]^s) = \pJmufJI(\px)
\quad\text{and}\quad
(\hf_*^k \hmu_\pz)([JI]^s) = \hmu_\pz([J]^s) = \pJmuflJ(\pz).
$$
Since $\pfJI = \pflJ \circ \pfkI$, we have that
\begin{equation}\label{eq.jaciscont}
\pJmufJI(\px) = \pJmufkI(\px) \pJmuflJ(\pz)
 \quad\text{at $\pmu$-almost every point.}
\end{equation}
Using the continuity property \eqref{eq.continuousjacnew}, one
concludes that the equality in \eqref{eq.jaciscont} holds
everywhere on $\supp\pmu=\psig$. Replacing the previous pair of
relations, we find that
$$
 \hmu_{\px}([JI]^s)
 = \pJmufkI(\px) (\hf_*^k \hmu_\pz)([JI]^s)
$$
for every $\pz\in\psig$ and any choice of $J=(\iota_l, \ldots, \iota_{-1})$.
This means that
$$
(\hmu_{\px} \mid [I]^s) = \pJmufkI(\px)(\hf_*^k \hmu_\pz),
$$
which, in view of the definition \eqref{eq.defjac2}, is just another
way of writing the claim in the lemma. The dual statement is proved
in just the same way.
\end{proof}

\subsection{Backward averages}\label{ss.backwardaverages}

For each $\px\in\psig$ and $k\ge 1$ let the \emph{backward
average} measure $\pmu_{k,\px}$ of the map $\pf$ be defined on
$\psig$ by
$$
\pmu_{k,\px}
 = \sum_{\pfk(\pz) = \px} \frac{1}{\pJmufk(\pz)} \, \delta_{\pz}
 = \sum_{I} \pJmufkI(\px) \, \delta_{\pfkI(\px)},
$$
where the last sum is over all $I=(\iota_0, \ldots, \iota_{k-1})$.
From \eqref{eq.defjac1} we get that
\begin{equation}\label{eq.todos2}
\sum_{\pfk(\pz) = \px} \frac{1}{\pJmufk(\pz)}
 = \sum_{I} \pJmufkI(\px)
 = \sum_{I} \hmu_\px([I]^s)
 = 1
\end{equation}
for every $\px\in\psig$. In other words, every $\pmu_{k,\px}$ is a
probability measure. The definition also implies that
$$
\int\pmu_{k,\px}(F)\,d\pmu(\px) = \sum_{I} \int_{\pfk(F\cap \pI)}
\pJmufkI \, d\pmu = \sum_{I} \pmu(F \cap \pI) = \pmu(F)
$$
for every measurable subset $F$ of $\psig$. Thus,
\begin{equation}\label{eq.invariante2}
\int \int \psi(\pz)\, d\pmu_{k,\px}(\pz)\,d\pmu(\px) =
\int\psi(\px)\,d\pmu(\px)
\end{equation}
for any bounded measurable function $\psi$ on $\psig$. It is important
to notice that the next result is stated for every (not just almost every)
point $\px$:

\begin{lemma}\label{l.backwardconvergence}
For every $\px \in \psig$ and every cylinder $\pJ\subset\psig$,
$$
K \pmu(\pJ) \ge \limsup_n \frac 1n \sum_{k=0}^{n-1}
\pmu_{k,\px}(\pJ) \ge \liminf_n \frac 1n \sum_{k=0}^{n-1}
\pmu_{k,\px}(\pJ) \ge \frac{1}{K}\pmu(\pJ)
$$
\end{lemma}

\begin{proof}
Given any positive $\pmu$-measure set $X\subset\psig$, define
$$
\pmu_{k,X}= \frac{1}{\pmu(X)} \int_X \pmu_{k,\pz}\,d\pmu(\pz).
$$
From the definition of the Jacobian one gets that
$$
\pmu_{k,X}(F)=\frac{1}{\pmu(X)} \pmu(F\cap (\pf)^{-k}(X))
$$
for every measurable set $F$ and every $k\ge 1$. Since $\pmu$ is
ergodic, it follows that
\begin{equation}\label{eq.convX}
\frac 1 n \sum_{k=0}^{n-1} \pmu_{k,X}(F) \to \pmu(F).
\end{equation}
Take $F=\pJ$ and $X=\psig$. Assuming $k$ is larger than the length
of $J$, we have that $\pfkI(X)=[I]^u$ intersects $\pJ$ if and only
if it is contained in it. Then, $\pfkI(\py)\in \pJ$ if and only if
$\pfkI(\px)\in \pJ$, for any $\py\in X$. Together with
\eqref{eq.boundedjacnew}, this implies that
$$
\frac{1}{K} \le \frac{\pmu_{k,\py}(\pJ)}{\pmu_{k,\px}(\pJ)} \le K
$$
for all $\py \in X$, and so
$$
\frac{1}{K} \le \frac{\pmu_{k,X}(\pJ)}{\pmu_{k,\px}(\pJ)} \le K.
$$
Combined with \eqref{eq.convX}, this implies the statement of the lemma.
\end{proof}

As a direct consequence, for every cylinder $\pJ\subset\psig$ and
every $\px\in\psig$,
\begin{equation}\label{eq.weak2}
\limsup_k \pmu_{k,\px}(\pJ) \ge K^{-1}\pmu(\pJ).
\end{equation}
This fact will be used in the proof of Lemma~\ref{l.tipic0}.


\subsection{Holonomy reduction}\label{ss.reduction}

Fix an arbitrary point $\xminus\in\nsig$ and then, for each $\hx\in\hsig$,
denote by $\pphi(\hx)$ the unique point in $\Wuloc(\xminus)\cap\Wsloc(\hx)$.
Using the stable holonomies in Definition~\ref{d.suholonomies}, define
$\hpA:\hsig\to\GL$ by
\begin{equation}\label{eq.holonomyreduction}
\hpA(\hx)
 = H^s_{\hf(\hx),\pphi(\hf(\hx))} \cdot \hA(\hx) \cdot H^s_{\pphi(\hx),\hx}.
 = H^s_{\hf(\pphi(\hx)),\pphi(\hf(\hx))} \cdot \hA(\pphi(\hx))
\end{equation}
Equivalently, the cocycle $\hpFA$ defined by $\hpA$ over $f$ is
conjugate to the cocycle $\hFA$ defined by $\hA$ through the
conjugacy
$$
\Phi : \hsig \times \CC^d \to \hsig \times \CC^d,
 \quad \Phi(\hx,v) = (\hx, H^s_{\hx,\pphi(\hx)}).
$$
Consequently, the two cocycles have the same Lyapunov exponents,
and either one is simple if and only if the other one is. So, for
the purpose of proving Theorem~\ref{main} one may replace $\hA$ by
either $\hpA$. On the other hand, the second equality in
\eqref{eq.holonomyreduction} implies that $\hpA$ is constant on
every local stable set, and so
$$\hpA(\hx)=\pA(\px)
 \quad \text{for some $\pA:\psig\to\GL$.}
$$
There is a dual construction, using unstable holonomies, where one
finds a map $\hnA:\hsig\to\GL$ that is constant on every local
unstable set and such that the cocycle it defines  over $f$ is
also conjugate to $\hFA$.

From now on, and until the end of Section~\ref{s.Dirac}, we consider
$\hpA$ instead of $\hA$. Notice that the corresponding stable holonomies
are trivial
$$
\Hsxy=\id \quad\text{for all $\hx$ and $\hy$,}
$$
because $\hpA$ is constant on local stable sets. For simplicity,
we omit the superscripts $u$ in the notations for
$\hpA$, $\pA$, $\hpFA$, $\psig$, $\pP$, $\pf$,
$\px$, $\pmu$, $\pom$, $H^u_{\hx,\hy}$, $\pfkI$, etc, that is, we
just represent these objects as $\hA$, $A$, $\hFA$, $\Sigma$, $P$,
$f$, $x$, $\mu$, $m$, $H_{\hx,\hy}$, $\fkI$, etc.

\section{Convergence of conditional probabilities}\label{s.convergence}

Let $\hpi:\hsig\times\Grass\to\hsig$ and
$\pi:\psig\times\Grass\to\psig$ be the natural projections. The
value of $\ell\in\{1, \ldots, d-1\}$ will be fixed till very near
the end. Note that if $\hm$ is an $\hFA$-invariant probability on
$\hsig\times\Grass$ then $m=(P\times\id)_*\hm$ is an
$\FA$-invariant probability on $\Sigma\times\Grass$. Moreover, if
$\hpi_*\hm=\hmu$ then $\pi_* m=\mu$.
Given $\hx\in\hsig$ we denote $\xn=P(\hf^{-n}(\hx))$ for $n\ge 0$.

\begin{proposition}\label{p.convergencia1}
Let $\hm$ be any $\hFA$-invariant probability on
$\hsig\times\Grass$ such that $\hpi_*\hm=\hmu$. Let $\{m_x :
x\in\Sigma\}$ be a disintegration of the measure
$m=(P\times\id)_*\hm$ along the Grassmannian fibers. Then the
sequence of probability measures
$$
A^n(\xn)_* m_{\xn}
$$
on $\Grass$ converges in the \weak topology as $n\to\infty$,
for $\hmu$-almost every $\hx\in\hsig$.
\end{proposition}

Starting the proof, let $\cB$ be the Borel $\sigma$-algebra of
$\Sigma$. Consider the sequence $(\cB_n)_n$ of $\sigma$-algebras of
$\hsig$ defined by $\cB_0=P^{-1}(\cB)$ and $\cB_n =\hf(\cB_{n-1})$
for $n\ge 1$. In other words, $\cB_n$ is the $\sigma$-algebra
generated by all cylinders $[\iota_{-n}, \ldots ; \iota_0; \ldots,
\iota_m]$ with $m\ge 0$ and $\iota_j\in\NN$. Fix any continuous
function $g:\Grass\to\RR$. For $\hx\in\hsig$ and $n\ge 0$, define
$$
\hI_n(\hx) = \hI_n(g,\hx) = \int g \, d\left(A^n(\xn)_*
m_{\xn}\right) = \int \left(g \circ A^n(\xn)\right) dm_{\xn}.
$$
Notice that $\hI_n$ is $\cB_n$-measurable: it can be written as
$\hI_n=I_n \circ P\circ \hf^{-n}$, where $I_n$ is the
$\cB$-measurable function
$$
I_n(x) = I_n(g,x) = \int \left(g \circ A^n(x)\right) dm_{x}.
$$

\begin{lemma}\label{l.convergence1}
For $\mu$-almost every $x\in\Sigma$ and any $n \ge 0$ and $k \ge
1$,
$$
I_n(x)
 = \sum_{z \in f^{-k}(x)} \frac{1}{\Jmufk(z)} I_{n+k}(z)
 = \int I_{n+k}(z)\,d\mu_{k,x}(z).
$$
\end{lemma}

\begin{proof}
Since the measure $m$ is invariant under $\FA^k$, its disintegration
must satisfy
\begin{equation}\label{eq.invar2}
m_x = \sum_{z \in f^{-k}(x)} \frac{1}{\Jmufk(z)} \, A^k(z)_* m_z
 = \int \Big(A^k(z)_*m_z\Big) \,d\mu_{k,x}(z)
\end{equation}
for $\mu$-almost every $x\in\Sigma$. Then,
\begin{equation*}
\begin{aligned}
I_n (x) & = \int \Big(g\circ A^n(x)\Big) dm_x = \int \Big(g\circ
A^n(x)\Big)
       d\Big(\sum_{z\in f^{-k}(x)} \frac{1}{\Jmufk(z)} A^k(z)_* m_z\Big)
\\
& = \sum_{z\in f^{-k}(x)} \frac{1}{\Jmufk(z)}
       \int \Big(g\circ A^{n+k}(z)\Big) dm_z
  = \sum_{z\in f^{-k}(x)} \frac{1}{\Jmufk(z)} I_{n+k} (z),
\end{aligned}
\end{equation*}
for $\mu$-almost every $x\in\Sigma$, as claimed.
\end{proof}

The next lemma means that each $\hI_n$ is the \emph{conditional
expectation} of $\hI_{n+k}$ with respect to the $\sigma$-algebra
$\cB_n$ for all $k\ge 1$, and so the sequence $(\hI_n,\cB_n)_n$ is
a martingale.

\begin{lemma}\label{l.convergence2}
For any $n\ge 0$ and $k\ge 1$ and any $\cB_n$-measurable function
$\psi:\hsig\to\RR$,
$$
\int \hI_{n+k}(\hx)\psi(\hx)\,d\hmu(\hx) = \int
\hI_{n}(\hx)\psi(\hx)\,d\hmu(\hx).
$$
\end{lemma}

\begin{proof}
Let us write $\psi=\psi_n\circ P\circ \hf^{-n}$, for some
$\cB$-measurable function $\psi_n$. Since $\hmu$ is $\hf$-invariant and
$\mu=P_*\hmu$,
\begin{equation}\label{eq.eg1}
\int \hI_n(\hx)\psi(\hx)\,d\hmu(\hx) = \int
I_n(x)\psi_n(x)\,d\mu(x).
\end{equation}
Analogously, using the relation $\psi=(\psi_n\circ f^k) \circ P\circ
\hf^{-(n+k)}$,
\begin{equation}\label{eq.eg2}
\int \hI_{n+k}(\hx)\psi(\hx)\,d\hmu(\hx) = \int
I_{n+k}(x)\psi_n(f^k(x))\,d\mu(x).
\end{equation}
By Lemma~\ref{l.convergence1}, the expression on the right hand side
of \eqref{eq.eg1} is equal to
$$
\int\int I_{n+k} (z) d\mu_{k,x}(z) \psi_n(x) \,d\mu(x)
 = \int \int I_{n+k} (z) \psi_n(f^k(z)) d\mu_{k,x}(z) \,d\mu(x).
$$
According to the relation \eqref{eq.invariante2}, this last
expression is the equal to the right hand side of \eqref{eq.eg2}.
This proves the claim of the lemma.
\end{proof}

\begin{proof}[Proof of Proposition~\ref{p.convergencia1}]
By Lemma~\ref{l.convergence2} and the martingale convergence
theorem (see Durret~\cite{Dur96}), the sequence
$\hI_n=\hI_n(g,\cdot)$ converges $\hmu$-almost everywhere to some
measurable function $\cI(g,\cdot)$. Notice that $|\hI_n(g,\hx)|
\le \sup|g|$ for every $n\ge 1$, and so $|\cI(g,\hx)|$ is also
bounded above by $\sup|g|$, for $\hmu$-almost every $\hx\in\hsig$.
Considering a countable dense subset of the space of continuous
functions, we find a full $\hmu$-measure set of points $\hx$ such
that
$$
\hI_n(g,\hx)=\int g \, d\left(A^n(\xn)_* m_{\xn}\right) \to
\cI(g,\hx)
$$
for every continuous function $g:\Grass\to\RR$. Let $\tm_\hx$ be
the probability measure on $\Grass$ defined by
$$
\int g \, d\tm_\hx = \cI(g,\hx) \quad\text{for every continuous
$g:\Grass\to\RR$.}
$$
Then the previous relation means that $A^n(\xn)_* m_{\xn}$
converges weakly to $\tm_\hx$.
\end{proof}

\begin{corollary}\label{c.gauche1}
For $\hmu$-almost every $\hx\in\hsig$, the limit of $A^n(\xn)_*
m_{\xn}$ coincides with the conditional probability $\hm_\hx$ of
the measure $\hm$.
\end{corollary}

\begin{proof}
Taking the limit $k\to\infty$ in Lemma~\ref{l.convergence2}, and
using the dominated convergence theorem, we get that
$$
\int \cI(g,\hx)\psi(\hx)\,d\hmu(\hx) = \int
\hI_n(g,\hx)\psi(\hx)\,d\hmu(\hx)
$$
for every $\cB_n$-measurable integrable function $\psi$. This may
be rewritten as
$$
\int \psi(\hx) \int g(\xi) \,d\tm_\hx(\xi) \,d\hmu(\hx)
 = \int \psi(\hx) \int g(A^n(\xn)\xi) \, d m_{\xn}(\xi) \,d\hmu(\hx).
$$
Let $\psi=\Chi_{[I]}$ be the characteristic function of a generic
cylinder $[I]$ in $\cB_n$. Changing variables $\hx=\hf^n(\hz)$, and
using the fact that $\hmu$ is $\hf$-invariant, we get that the right
hand side of the previous equality is equal to
$$
 \int \Chi_{[I]}(\hf^n(\hz)) \int g(A^n(z)\xi) \, d m_{z}(\xi) \,d\hmu(\hz)
$$
where $z=P(\hz)$. Moreover, since the inner integrand $z\mapsto
g(A^n(z)\xi)$ is constant on local stable leaves, this may be
rewritten as
$$
\int \Chi_{[I]}(\hf^n(\hz)) \int g(A^n(\hz)\xi) \, d
\hm_{\hz}(\xi) \, d\hmu(\hz) = \int \Chi_{[I]}(\hx) \int g(\eta)
\, d\hm_{\hx}(\eta) \, d\hmu(\hx).
$$
In the last step we changed variables $(\hx,\eta)=\hFA^n(\hz,\xi)$
and used the fact that $\hm$ is invariant under $\hFA$.
Summarizing, at this point we have shown that
$$
\int\int \Chi_{[I]}(\hx) g(\xi) \,d\tm_\hx(\xi) \,d\hmu(\hx)
 =
\int\int \Chi_{[I]} (\hx) g(\eta) \, d \hm_\hx(\eta) \,
d\hmu(\hx).
$$
This relation extends immediately to linear combinations of
functions $\Chi_{[I]} \times g$. Since these linear combinations
form a dense subset of all bounded measurable functions on
$\hsig\times\Grass$, this implies that $\tm_\hx=\hm_\hx$ for
$\hmu$-almost every $\hx$, as claimed.
\end{proof}

\section{Properties of $u$-states}\label{s.ustates}

Let $\hm$ be a probability measure on $\hsig\times\Grass$ that
projects down to $\hmu$ on $\hsig$, in the sense that
$\hpi_*\hm=\hmu$. We call $\hm$ a \emph{$u$-state} if it admits some
disintegration $\{\hm_\hx : \hx\in\hsig\}$ into conditional
probabilities along the fibers $\{\hx\}\times\Grass$ that is
invariant under unstable holonomies:
$$
\hm_\hy = (\Hxy)_* \hm_\hx \quad\text{whenever } y \in
\Wuloc(\hx).
$$
We call the $u$-state \emph{invariant} if, in addition, it is
invariant under $\hFA$. We also call (invariant) $u$-states the
projections $m=(P\times\id)_* \hm$ down to $\Sigma\times\Grass$ of
the (invariant) $u$-states $\hm$ on $\hsig\times\Grass$. Notice
that $\pi_* m = \mu$, and $m$ is invariant under $\FA$ if $\hm$ is
invariant under $\hFA$.

Here we prove that invariant $u$-states $m$ do exist. Moreover,
every $u$-state admits some disintegration $\{m_x : x\in\Sigma\}$
into conditional probabilities along the fibers
$\{x\}\times\Grass$ varying continuously with the base point $x$,
relative to the \weak topology. The formal statements are in
Propositions~\ref{p.existence} and~\ref{p.continuity}. The proofs
use the assumption that $\hmu$ has product structure (recall
Section~\ref{sh.measure}).

\subsection{Existence of invariant $u$-states} \label{ss.existence}

Let $\cM$ be the space of probability measures on
$\hsig\times\Grass$ that project down to $\hmu$ on $\hsig$. The
\emph{\weak topology} on $\cM$ is the smallest topology such that
the map $\eta \mapsto \int \psi d\eta$ is continuous, for every
bounded continuous function $\psi:\hsig \times \Grass \to \RR$.
Notice that $\cM$ is a compact separable space for this topology.
This is easy to see from the following alternative description of
the topology. Let $K_n \subset \hsig$, $n\ge 1$ be pairwise disjoint
compact sets such that $\hmu(K_n)>0$ and $\sum \hat \mu(K_n)=1$. Let
$\cM_n$ be the space of measures on $K_n \times \Grass$ that project
down to $(\hmu \mid K_n)$. The usual \weak topology makes $\cM_n$ a
compact separable space. Given $\eta \in \cM$, let $\eta_n \in
\cM_n$ be obtained by restriction of $\eta$. The correspondence
$\eta\mapsto (\eta_n)_n$ identifies $\cM$ with $\prod \cM_n$ and the
product topology on $\prod \cM_n$ corresponds to the \weak topology
on $\cM$ under this identification. Thus, the latter is a compact
separable space, as claimed.

\begin{remark}\label{r.boundedfactor}
If $\eta^j$ converges to $\eta$ in the \weak topology then
\begin{equation}\label{eq.convfactor}
\int \psi(\hx,\xi)J(\hx)\,d\eta^j(\hx,\xi) \to \int
\psi(\hx,\xi)J(\hx)\,d\eta(\hx,\xi)
\end{equation}
for any continuous function $\psi:\hsig \times \Grass \to \RR$ and
any measurable bounded (or even $\hmu$-integrable) function $J:\hsig
\to \RR$. To prove this it suffices to consider the case when
$J=\Chi_B$ for some measurable set $B$, because every bounded
measurable function is the uniform limit of linear combinations of
characteristic functions. Now, using that $\hmu$ is a regular
measure (see Theorem~6.1 in \cite{Wa82}), we may find continuous
functions $J_n:\hsig\to [0,1]$ such that $\hmu(\{\hx\in\hsig:
J_n(\hx)\neq J(\hx)\}$ is arbitrarily small. By the definition of
the topology,
$$
 \int \psi(\hx,\xi)J_n(\hx)\,d\eta^j(\hx,\xi)
 \to \int \psi(\hx,\xi)J_n(\hx)\,d\eta(\hx,\xi)
 \quad\text{as $j\to\infty$.}
$$
This implies the convergence in \eqref{eq.convfactor}, because
corresponding terms in these two relations differ by not more than
$\sup|\psi| \, \hmu(\{\hx\in\hsig : J_n(\hx)\neq J(\hx)\}$, which
can be made arbitrarily small.
\end{remark}

Remark also, for future use, that in these arguments $\hmu$ may be
replaced by any other probability in $\hsig$.

\begin{proposition}\label{p.existence}
There exists some invariant $u$-state $\hm$ on
$\hsig\times\Grass$.
\end{proposition}

Here is an outline of the proof. The space $\cU$ of all $u$-states
is non-empty and forward invariant under the cocycle. Every Cesaro
\weak limit of the forward iterates of an element of $\cU$ is an
invariant $u$-state. The proposition follows by noting that \weak
limits do exist, because $\cU$ is compact relative to the \weak
topology. The last step demands some caution, because conditional
probabilities do not behave well under \weak limits, in general. We
fix an arbitrary point $\xplus\in\Sigma$ and observe that,
restricted to the cylinder, the space $\cU$ may be identified with
the space $\cN$ of probabilities on $\Wsloc(\xplus)\times\Grass$
that project down to $\hmu_{\xplus}$. Then it suffices to use that
the latter space is \weak compact.

Let us fill the details. Let $\{\hmu_x : x \in \Sigma\}$ be the
disintegration of $\hmu$ along local stable sets in Section~\ref{sh.measure}.
Denote by $\Jmux$ the Radon-Nikodym derivatives of the conditional
measure $\hmu_x$ with respect to $\hmu_\xplus$, for each $x\in\Sigma$.
According to \eqref{eq.boundedjacnew}, these $J_x$ are uniformly bounded
from zero and infinity. We use $\hx$ and $\hxplus$ to denote generic
points in $\Wsloc(x)$ and $\Wsloc(\xplus)$, respectively, with the
convention that whenever they appear in the same expression they are
related by
$$
\hxplus \in \Wsloc(\xplus) \cap \Wuloc(\hx).
$$
Let $\cN$ be the space of all probability measures $\lambda$ on
$\Wsloc(\xplus)\times\Grass$ that project down to $\hmu_{\xplus}$ on
$\Wsloc(\xplus)$. Recall, from the observation at the beginning of
this section, that $\cN$ is \weak compact. We denote by $\cU$ the
space of all $u$-states, that is, all probability measures $\eta$ on
$\hsig\times\Grass$ that project down to $\hmu$ and admit some
disintegration $\{\eta_\hx : \hx\in\hsig\}$ along the Grassmannian
fibers that is invariant under unstable holonomy:
\begin{equation}\label{eq.pt1}
\eta_\hx = (\Hxplus)_* \eta_\hxplus \quad\text{for all }\hx\in\hsig.
\end{equation}

\begin{lemma}\label{l.compact2}
$\cU$ is homeomorphic to $\cN$.
\end{lemma}

\begin{proof}
Every $\lambda\in\cN$ may be lifted to some $\eta\in\cU$ in the
following natural fashion: choose a disintegration $\{\lambda_\hw
: \hw \in \Wsloc(\xplus)\}$ of $\lambda$ and then let $\eta$ be
the measure on $\hsig\times\Grass$ whose projection coincides with
$\hmu$ and which admits
\begin{equation}\label{eq.pt2}
\eta_\hx=(\Hxplus)_* \lambda_{\hxplus}
\end{equation}
as conditional probabilities along the fibers
$\{\hx\}\times\Grass$. This definition does not depend on the
choice of the disintegration of $\lambda$. Indeed, let
$\{\tilde{\lambda}_\hw: \hw\in\Wsloc(\xplus)\}$ be any other
disintegration. By essential uniqueness, we have
$$
\tilde{\lambda}_\hw = \lambda_\hw \quad
 \text{for $\hmu_{\xplus}$-almost every $\hw\in W_{loc}^s(\xplus)$.}
$$
Since the measures $\hmu_x$, $x\in\Sigma$ are all equivalent, it
follows that $\tilde\eta_\hx = \eta_\hx$ for $\hmu_x$-almost every
$\hx\in W_{loc}^s(x)$ and every $x\in\Sigma$. So, the lifts
constructed from the two disintegrations do coincide. It is clear
from the construction that $\eta\in\cU$.

Let $\Psi:\cN\to\cU$, $\Psi(\lambda)=\eta$ be the map
defined in this way. We are going to prove that $\Psi$ is a
homeomorphism. To prove injectivity, suppose $\Psi(\lambda) = \hm =
\Psi(\theta)$. By \eqref{eq.pt2}, this means that
$$
(\Hxplus)_* \lambda_{\hxplus} = \hm_\hx = (\Hxplus)_* \theta_{\hxplus}
$$
for $\hmu$-almost every $\hx\in\hsig$. Since the conditional
probabilities $\hmu_x$ are all equivalent, this is the same as
$\lambda_\hw=\theta_\hw$ for $\hmu_\xi$-almost every $\hw\in
W_{loc}^s(\xplus)$. In other words, $\lambda=\theta$. To prove
surjectivity, consider any measure $\eta\in\cU$. By definition,
$\eta$ admits some disintegration $\{\eta_x : x\in\Sigma\}$
satisfying \eqref{eq.pt1}. Define
$\lambda_\hxplus=(\Hplusx)_*\eta_{\hx}$ for any
$\hx\in\Wuloc(\hxplus)$, and then let $\lambda$ be the measure on
$W_{loc}^s(\xplus)\times\Grass$ that projects down to
$\hmu_\xplus$ and has these $\lambda_\hxplus$ as conditional
probabilities along the fibers. Then $\lambda\in\cN$ and
$\eta=\Psi(\lambda)$.

We are left to check that $\Psi$ is continuous. Let
$\psi:\hsig\times\Grass\to\RR$ be any bounded continuous function
and let $\lambda^j$ be any sequence of measures converging to some
$\lambda$ in $\cN$. Using Remark~\ref{r.boundedfactor},
$$
\int \psi(x,\hx,\xi) \, d\lambda^j_\hx(\xi) \,d\hmu_x(\hx)
 = \int \psi(x,\hx,\xi) J_x(\hxplus)\,d\lambda^j_\hw(\xi)\,d\hmu_{\xplus}(\hxplus)
$$
converges to
$$
\int \psi(x,\hx,\xi) \, d\lambda_\hx(\xi) \,d\hmu_x(\hx) =
 \int \psi(x,\hx,\xi) J_x(\hxplus) \, d\lambda_\hw(\xi) \,d\hmu_{\xplus}(\hxplus)
$$
as $j\to \infty$, for every $x\in\Sigma$. Integrating with respect
to $\mu$, and using the bounded convergence theorem, we get that
$$
 \int\int \psi(x,\hx,\xi) \, d\lambda^j_\hx(\xi) \, d\hmu_x(\hx) \, d\mu(x)
 \to
 \int\int \psi(x,\hx,\xi) \, d\lambda_\hx(\xi) \, d\hmu_x(\hx) \, d\mu(x)
$$
as $j\to \infty$. This means that $\Psi(\lambda^j)$ converges to
$\Psi(\lambda)$ as $j\to\infty$.
\end{proof}

\begin{proof}[Proof of Proposition~\ref{p.existence}]
In view of the previous lemma, $\cU$ is non-empty and compact
relative to the \weak topology. Moreover, $\cU$ is invariant under
iteration by $\hFA$: this follows from the invariance property (b)
in Section~\ref{sh.cocycle} for unstable holonomies, together with
the fact that local unstable sets are mapped inside local unstable
sets by the inverse of $\hf$. Consider any probability measure $\bar
m\in\cU$. The sequence
$$
\hm_n = \frac 1n \sum_{j=0}^{n-1}(\hFA^j)_* \bar m
$$
has accumulation points $\hm$ in $\cU$. Since $\hFA$ is a
continuous map, the push-forward operator $(\hFA)_*$ is continuous
relative to the \weak topology. It follows that any such
accumulation point is $\hFA$-invariant and, consequently, an
invariant $u$-state.
\end{proof}

\subsection{Continuity of conditional probabilities}
\label{ss.continuity}

Now we prove that conditional probabilities of $u$-states along
the Grassmannian fibers depend continuously on the base point:

\begin{proposition}\label{p.continuity}
Any $u$-state $m$ in $\Sigma\times\Grass$ admits some
disintegration $\{m_x : x\in\Sigma\}$ into conditional
probabilities along the Grassmannian fibers varying continuously
with $x\in\Sigma$ in the \weak topology.
\end{proposition}

This continuous disintegration is necessarily unique, because
disintegrations are essentially unique and $\mu$ is supported on the
whole $\Sigma$. For the proof of the proposition we need the
following simple observation:

\begin{lemma}\label{l.conditional2}
Let $\{\hm_\hx : \hx\in\hsig\}$ be a disintegration along
$\{\{\hx\}\times\Grass:\hx\in\hsig\}$ of some probability measure
$\hm$ on $\hsig\times\Grass$ such that $\hpi_*\hm=\hmu$. Then
$$
m_x=\int \hm_\hx\,d\hmu_x(\hx)
$$
is a disintegration of $m=(P\times\id)_* \hm$ along
$\{\{x\}\times\Grass:x\in\Sigma\}$.
\end{lemma}

\begin{proof}
For any $\varphi:\Sigma\times\Grass\to\RR$ and
$\hat\varphi=\varphi\circ(P\times\id)$,
\begin{equation*}
\begin{aligned}
\int \int \varphi \, dm_x \,d\mu(x) & = \int \int\Big( \int
\varphi (x,v) \, d\hm_\hx(v)  d\hmu_x(\hx)\Big) d\mu(x)
\\
& = \int \int \Big(\int \hat\varphi (x,v) \, d\hm_\hx (v) \Big)
d\hmu_x(\hx)\,d\mu(x)
\\
& = \int \left(\int \hat\varphi (x,v) \, d\hm_\hx (v) \right)
d\hmu(\hx)
 = \int \hat\varphi d\hm = \int \varphi dm
\end{aligned}
\end{equation*}
and this proves that $\{m_x : x\in\Sigma\}$ is a disintegration of
$m$.
\end{proof}

\begin{remark}\label{r.avoidcontinuity}
For $u$-states this gives that, for any measurable set
$E\subset\Grass$,
$$
m_x(E)
 = \int \hm_\hx (E) \,d\hmu_x(\hx)
 = \int \hm_\hy (\Hsxy(E)) \frac{d\hmu_x}{d\hmu_y} (\hy) \,d\hmu_y(\hy)
$$
for any pair of points $x$ and $y$ in the same cylinder. When the
cocycle is \emph{locally constant} the stable holonomies
$\Hsxy=\id$. In this case it immediately follows that the
conditional probabilities $m_x$ and $m_y$ are all equivalent.
Moreover, their Radon-Nikodym derivatives are uniformly bounded, as
a consequence of the boundedness condition \eqref{eq.boundedweak}.
Starting from this observation, in the appendix of \cite{AV2} we
give a version of Theorem~\ref{main} for locally constant cocycles
that does not require the continuity hypothesis
\eqref{eq.continuousweak}.
\end{remark}

\begin{proof}[Proof of Proposition~\ref{p.continuity}]
Let $\{\hm_\hx : \hx\in\hsig\}$ be a disintegration of $\hm$ into
conditional probabilities that are invariant under unstable
holonomies: $\hm_\hy=(\Hxy)_*\hm_\hx$ for every $\hx, \hy$ in the
same local unstable set. Let $\{m_x : x\in\Sigma\}$ be the
disintegration of $m$ given by Lemma~\ref{l.conditional2}. For any
continuous $g:\Grass\to\RR$ and any points $x$ and $y$ in the same
cylinder of $\Sigma$, we have
$$
\int g(\xi)\,dm_x(\xi)
 = \int \int g(\xi) \, d\hm_\hx(\xi) \, d\hmu_x(\hx)
 = \int \int g(\Hyx(\eta)) \, d\hm_\hy(\eta) \, d\hmu_x(\hx)
 $$
where $\hy$ denotes the unique point in $\Wsloc(y) \cap
\Wuloc(\hx)$. Fix $y$ and consider the function
$$
\Phi(\nx,\px) = \int g(\Hyx(\eta)) \, d\hm_\hy(\eta),
 \quad\text{where } \hx = (\nx,\px).
$$
It is clear that $\Phi$ is measurable and bounded by the
$\sup|g|$. Moreover, it is continuous on $\px$ for each fixed
$\nx$. To check this it suffices to note that $\hm_\hy$ does not
depend on $\px$, while the function $g$ and the holonomies depend
continuously on $\hx$ (recall Definition~\ref{d.suholonomies}). It
follows from Lemma~\ref{l.continuity} that
$$
x\mapsto \int g(\xi)\,dm_x(\xi) = \int \Phi(\nx,x) \, d\hmu_x(\nx)
$$
is continuous. This proves the claim of the proposition.
\end{proof}

\begin{corollary}\label{c.invcont1}
If $m$ is an invariant $u$-state and $\{m_x : x\in\Sigma\}$ is the
continuous disintegration of $m$, then
$$
m_x
 = \sum_{z \in f^{-k}(x)} \frac{1}{\Jmufk(z)} \, A^k(z)_* m_z
 = \int A^k(z)_*m_z\,d\mu_{k,x}(z)
$$
for every $x\in\Sigma$ and every $k\ge 1$.
\end{corollary}

\begin{proof}
The second equality is just the definition of the backward averages,
see Section~\ref{ss.backwardaverages}. As for the first equality, it
must hold for every $k\ge 1$ and $\mu$-almost every $x$, because $m$
is invariant under $f$. Moreover, all the expressions involved vary
continuously with $x\in\Sigma$: this follows from
Proposition~\ref{p.continuity}, property
\eqref{eq.continuousjacnew}, and our assumption that the cocycle is
continuous. Hence, the first equality must hold at every point of
$\supp\mu=\Sigma$.
\end{proof}

\begin{corollary}\label{c.invcont2}
If $\{\hm_\hx : \hx\in\hsig\}$ is a disintegration of an invariant $u$-state $\hm$
into conditional probabilities invariant under unstable holonomies then
$$
\hm_{\hf^n(\hx)} = A^n(x)_* \hm_{\hx}
$$
for every $n\ge 1$, every $x\in\Sigma$, and $\hmu_x$-almost every $\hx\in \Wsloc(x)$.
\end{corollary}

\begin{proof}
Since $\hm$ is $\hFA$-invariant, the equality is true for all $n\ge
1$ and $\hmu$-almost all $\hz\in\hsig$ or, equivalently, for
$\hmu_z$-almost every $\hz\in\Wsloc(z)$ and $\mu$-almost every
$z\in\Sigma$. Consider an arbitrary point $x\in\Sigma$. Since $\mu$
is positive on open sets, $x$ may be approximated by points $z$ such
that
$$
\hm_{\hf^n(\hz)} = A^n(z)_* \hm_{\hz}
$$
for every $n\ge 1$ and $\hmu_z$-almost every $\hz\in\Wsloc(z)$. Since the conditional
probabilities of $\hm$ are invariant under unstable holonomies, it follows that
$$
\hm_{\hf^n(\hx)}
 = (H_{\hf^n(z),\hf^n(x)})_* A^n(z)_* \hm_{\hz}
 = A^n(x)_* (H_{\hz,\hx})_* \hm_{\hz}
 = A^n(x)_* \hm_{\hx}
$$
for $\hmu_z$-almost every $\hz\in\Wsloc(z)$, where $\hx$ is the unique point in
$\Wsloc(x)\cap\Wuloc(\hz)$. Since the measures $\hmu_x$ and $\hmu_z$ are equivalent,
this is the same as saying that the last equality holds for $\hmu_x$-almost every
$\hx\in \Wsloc(x)$, as claimed.
\end{proof}

\section{Invariant measures of simple cocycles}\label{s.simple}

In this section we prove that invariant $u$-states of simple cocycles are fairly
smooth along the Grassmannian fibers: they give zero weight to every hyperplane section.

\begin{proposition}\label{p.noatoms}
Suppose that $\hA$ is simple. Let $m$ be any invariant $u$-state in
$\Sigma\times\Grass$ and $\{m_x : x\in\Sigma\}$ be the continuous
disintegration of $m$. Then $m_x(V)=0$ for every $x\in\Sigma$ and
any hyperplane section $V$ of $\Grass$.
\end{proposition}

In Section~\ref{ss.smoothness} we argue by contradiction to reduce
the proof of Proposition~\ref{p.noatoms} to
Proposition~\ref{p.artur}, a combinatorial result about
intersections of hyperplane sections. The latter is proved in
Section~\ref{ss.disjointness}. See also Appendix~\ref{s.AppB}.

\subsection{Smoothness of conditional probabilities} \label{ss.smoothness}

Suppose there is some point of $\Sigma$ and some hyperplane section
of the corresponding Grassmannian fiber which has positive
conditional measure. Let $\gamma_0>0$ be the supremum of the values
of $m_x(V)\ge\gamma$ over all $x\in\Sigma$ and all hyperplane
sections $V$. The supremum is attained at every point:

\begin{lemma}\label{l.tipic0}
For every $x\in\Sigma$ there exists some hyperplane section $V$ of
$\Grass$ such that $m_x(V)=\gamma_0$.
\end{lemma}

\begin{proof}
Fix any cylinder $[J]\subset\Sigma$ and any positive constant
$c<\mu([J])/K$, where $K$ is the constant in \eqref{eq.weak2}. Let
$z\in\Sigma$ and $V$ be a hyperplane section with $m_z(V)>0$. For
each $y\in f^{-k}(z)$, let $V_y=A^k(y)^{-1}(V)$. By
Corollary~\ref{c.invcont1},
$$
m_z(V)
 = \int m_y(V_y) d\mu_{k,z}(y)
\le \mu_{k,z}([J]) \sup\{m_y(V_y):y\in [J]\}
 + (1-\mu_{k,z}([J])) \gamma_0.
$$
By \eqref{eq.weak2}, there exist arbitrarily large values of $k$
such that $\mu_{k,z}([J])\ge c$. Then
$$
m_z(V) \le c \sup\{m_y(V_y) : y\in [J]\} + (1-c) \gamma_0\,.
$$
Varying the point $z\in\Sigma$ and the hyperplane section $V$, we
can make the left hand side arbitrarily close to $\gamma_0$. It
follows that
$$
\sup\{m_y(V_y) : y\in [J]\}\ge \gamma_0.
$$
This proves that the supremum over any cylinder $[J]$ coincides with
$\gamma_0$. Then, given any $x\in\Sigma$ we may find a sequence
$x_n\to x$ and hyperplane sections $V_n$ such that
$m_{x_n}(V_n)\to\gamma_0$. Moreover, we may assume that $V_n$
converges to some hyperplane section $V$ in the Hausdorff topology.
Given any neighborhood $U$ of $V$, we have $m_{x_n}(U)\ge
m_{x_n}(V_n)$ for all large $n$. By Proposition~\ref{p.continuity},
the conditional probabilities $m_{x_n}$ converge weakly to $m_x$.
Assuming $U$ is closed, it follows that
$$
m_{x}(U)\ge\limsup_n m_{x_n}(U) \ge \limsup_n m_{x_n}(V_n) \ge
\gamma_0.
$$
Making $U\to V$, we conclude that $m_{x}(V)\ge\gamma_0$. This proves
that the supremum $\gamma_0$ is realized at $x$, as claimed.
\end{proof}

\begin{lemma}\label{l.tipic1}
For any $x\in\Sigma$ and any hyperplane section $V$ of $\Grass$,
we have $m_x(V)=\gamma_0$ if and only if
$m_y(A(y)^{-1}V)=\gamma_0$ for every $y\in f^{-1}(x)$.
\end{lemma}

\begin{proof}
This is a direct consequence of Corollary~\ref{c.invcont1} and the
relation \eqref{eq.todos2}: for every $x\in\Sigma$,
$$
m_x(V) = \sum_{y\in f^{-1}(x)} \frac{1}{\Jmuf(y)} \,
m_y(A(y)^{-1}V) \quand
 \sum_{y\in f^{-1}(x)} \frac{1}{\Jmuf(y)}=1.
$$
Since $\gamma_0$ is the maximum value of the measure of any
hyperplane section, we get that $m_x(V)=\gamma_0$ if and only if
$m_y(A(y)^{-1}V)=\gamma_0$ for every $y\in f^{-1}(x)$, as stated.
\end{proof}

\begin{lemma}\label{l.tipic2}
For any $x\in\Sigma$ and any hyperplane section $V$ of $\Grass$,
we have $\hm_\hx(V) \le \gamma_0$ for $\hmu_x$-almost every
$\hx\in \Wsloc(x)$. Hence, $m_x(V)=\gamma_0$ if and only if
$\hm_\hx(V) = \gamma_0$ for $\hmu_x$-almost every $\hx\in
\Wsloc(x)$.
\end{lemma}

\begin{proof}
Suppose there is $y\in\Sigma$, a hyperplane section $V$, a constant
$\gamma_1>\gamma_0$, and a positive $\hmu$-measure subset $X$ of
$\Wsloc(y)$ such that $\hm_\hy(V) \ge \gamma_1$ for every $\hy\in
X$. For each $m<0$, consider the partition of
$\Wsloc(y)\approx\nsig$ determined by the cylinders $[I]^s=[\iota_m,
\ldots, \iota_{-1}]^s$, with $\iota_j\in\NN$. Since these partitions
generate the $\sigma$-algebra of the local stable set, given any
$\vep>0$ we may find $m$ and $I$ such that
$$
\hmu_y(X \cap [I]^s) \ge (1-\vep) \hmu_y([I]^s).
$$
Observe that $[I]^s \approx [I]^s \times \{y\}$ coincides with
$\hf^n(W_{loc}^s(x))$, where $x=f_I^{-n}(y)$. So, using also
Lemma~\ref{l.todos1},
$$
\hmu_x\big(\hf^{-n}(X) \cap \Wsloc(x)\big)
 = (\hf_*^n\hmu_x)(X \cap [I]^s)
 = J_\mu f^n(x) \, \hmu_y\big(X \cap [I]^s).
$$
By the previous inequality and Lemma~\ref{l.todos1}, this is
bounded below by
$$
 (1-\vep) \, J_\mu f^n(x) \, \hmu_y\big([I]^s)
  = (1-\vep) \, (\hf_*^n\hmu_x) \big([I]^s\big)
  = \hmu_x \big(\Wsloc(x)\big) = 1 -\vep.
$$
In this way we have shown that
$$
\hmu_y\big(\hf^{-n}(X) \cap \Wsloc(x)\big) \ge (1-\vep).
$$
Fix $\vep>0$ small enough so that $(1-\vep)\gamma_1>\gamma_0$. Using
Corollary~\ref{c.invcont2}, we find that
$$
\hm_\hx(A^n(x)^{-1}V)
 = \hm_{\hy}(V)\ge \gamma_1
$$
for $\hmu_x$-almost every $\hx\in \hf^{-n}(X)\cap \Wsloc(x)$. It
follows that
$$
m_x(A^n(x)^{-1}V) = \int \hm_\hx(A^n(x)^{-1}V)\,d\hmu_x(\hx) \ge
(1-\vep) \gamma_1 >\gamma_0,
$$
which contradicts the definition of $\gamma_0$. This contradiction
proves the first part of the lemma. The second one is a direct consequence,
using the fact that $m_x(V)$ is the $\hmu_x$-average of all $\hm_\hx(V)$.
\end{proof}

Before we proceed, let us introduce some useful terminology.
Recall that a hyperplane section $V$ of $\Grass$ is the image of
$\extlv{\CC^d}\cap H$ under the projection $\pi_v$, where $H$ is
the geometric hyperplane of $\extl{\CC^d}$ defined by some
non-zero element $\upsilon\in\extdlv{\CC^d}$. Notice that, for any
linear isomorphism $B$ of $\CC^d$,
$$
 \extl{B} H = \{\omega : \omega \wedge \extdl{B}(\upsilon) = 0\}
 \quand
 B(V)=\pi_v\big(\extlv{\CC^d}\cap \extl{B} H\big).
$$
Suppose $B$ is diagonalizable. Then we say $V$ is \emph{invariant}
for $B$ if the subspace $\pi_v(\upsilon)$ is a sum of eigenspaces
of $B$. Likewise, we say $V$ \emph{contains no eigenspace} of $B$
if $\pi_v(\upsilon)$ intersects any sum of $\ell$ eigenspaces of
$B$ at the origin only or, equivalently, if $H$ contains no
$\ell$-vector $\omega$ such that $\pi_v(\omega)$ is a sum of
eigenspaces of $B$. A subset $J$ of $\{0, 1, \ldots, N-1\}$ is
called \emph{$\vep$-dense} if $\# J \ge N\vep$.


\begin{proposition}\label{p.artur}
For any $\vep>0$ there exists $N\ge 1$ such that
$$\bigcap_{j \in J} B^j(V)=\emptyset$$ for every $\vep$-dense set
$J\subset\{0, 1, \ldots, N-1\}$, every linear isomorphism
$B:\CC^d \to \CC^d$ whose eigenvalues all have distinct absolute
values, and every hyperplane section $V$ of $\Grass$ containing
no eigenspace of $B$.
\end{proposition}

This proposition will be proved in Section~\ref{ss.disjointness}.
Right now let us explain how it can be used to finish the proof of
Proposition~\ref{p.noatoms}.

Fix a periodic point $\hp\in\hsig$ of $\hf$ and a homoclinic point
$\hz\in\hsig$ as in Definition~\ref{d.simple}. Let $p=P(\hp)$ be
the corresponding periodic point of $f$ and let $z=P(\hz)$. By
Lemma~\ref{l.tipic0}, we may find a hyperplane section $V$ of
$\Grass$ with $m_p(V)=\gamma_0$. Write $V=\pi_v(\extlv{\CC^d}\cap
H)$, where $H$ is the geometric hyperplane defined by some
non-zero $(d-\ell)$-vector $\upsilon$. Let $V^n=A^{-nq}(p) V$ and
$H^n$ be the geometric hyperplane defined by $A^{-nq}(p)\upsilon$.
Then, $V^n = \pi_v(\extlv{\CC^d} \cap H^n)$ for each $n\ge 0$.
Since all the eigenvalues of $A^q(p)$ have distinct absolute
values, $A^{-nq}(p)\upsilon$ converges to some $(d-\ell)$-vector
$\upsilon_1$ such that $\pi_v(\upsilon_1)$ is a sum of eigenspaces
of $A^q(p)$. This means that $V^n$ converges to
$V_1=\pi_v(\extlv{\CC^d} \cap H_1)$, where $H_1 = \{\omega:
\omega\wedge \upsilon_1=0\}$ is the geometric hyperplane section
defined by $\upsilon_1$. On the other hand, using
Lemma~\ref{l.tipic1} we find that $m_p(V^n)=\gamma_0$ for all
$n\ge 0$. By lower semi-continuity of the measure, it follows that
$m_p(V_1)=\gamma_0$. Note that $V_1$ is invariant for $A^q(p)$.
This shows that we may suppose, right from the start, that $V$ is
invariant for $A^q(p)$.

Now define $W=A^l(z)^{-1}V$. From Lemmas~\ref{l.tipic1} and~\ref{l.tipic2}
we get that $m_z(W)=\gamma_0$ and $\hm_\zeta(W)=\gamma_0$ for
$\hmu_z$-almost every $\zeta\in W^s_{loc}(z)$. For each $\eta\in W^s_{loc}(p)$,
define $W_\eta=H_{\zeta,\eta}(W)$, where $\zeta$ is the unique point in
$\Wuloc(\eta)\cap\Wsloc(\hz)$. Since $\hm$ is a $u$-state and the measures
$\hmu_z$ and $\hmu_p$ are equivalent, we have
$\hm_\eta(W_\eta)=\hm_\zeta(W)=\gamma_0$ for $\hmu_p$-almost every
$\eta$. For each $j\ge 0$, let
$$
W^j_\eta
 =A^{-jq}(p)W_{\hf^{jq}(\eta)}
\qquad\text{(in particular, $W^j_\hp=A^{-jq}(q)(W_\hp)$).}
$$
Using Corollary~\ref{c.invcont2}, we get that
$\hm_\eta(W^j_\eta)=\hm_{\hf^{jq}(\eta)}(W_{\hf^{jq}(\eta)})=\gamma_0$
for every $j\ge 0$ and $\hmu_p$-almost every $\eta$. It is clear
that every $W_\eta^j$ is an $\ell$-dimensional projective subspace.
Moreover, it depends continuously on $\eta$, for each fixed $j$,
because unstable holonomies vary continuously with the base points
(Definition~\ref{d.suholonomies}). Notice that
$$
W_\hp = H_{\hz,\hp} A^l(z)^{-1} V = \psipz^{-1} V
 \quad\text{(recall $H_{\hz,\hp}=H^u_{\hz,\hp}$ and
 $H^s_{\hp,\hf^l(\hz)}=\id$).}
$$
Thus, the second condition in Definition~\ref{d.simple} implies
that $W_\hp$ contains no eigenspace of $A^q(p)$.

Taking $\vep=\gamma_0$, $V=W_\hp$, $B=A^q(p)$ in
Proposition~\ref{p.artur} we find $N\ge 1$ such that
$$
\bigcap_{j\in J} W^j_\hp = \emptyset
\quad\text{for every $\gamma_0$-dense subset $J$ of $\{0, 1, \ldots, N-1\}$.}
$$
Since the family of sets $J$ is finite, we may use continuity to
conclude that
\begin{equation}\label{eq.disjointW}
\bigcap_{j\in J} W^j_\eta = \emptyset
\quad\text{for every $\gamma_0$-dense subset $J$ of $\{0, 1, \ldots, N-1\}$}
\end{equation}
and any $\eta$ in a neighborhood of $\hp$ inside the local stable
set. On the other hand, the fact that
$\hm_\eta(W_\eta^j)=\gamma_0$ for all $j\ge 0$ implies (use a
Fubini argument) that there exists some $\omega\in\Grass$ such
that the set
$$
J=\{0\le j\le N-1: \omega\in W_\eta^j\}
$$
is $\gamma_0$-dense in $\{0, 1, \ldots, N-1\}$. This contradicts
\eqref{eq.disjointW}. This contradiction shows that we have reduced
the proof of Proposition~\ref{p.noatoms} to proving
Proposition~\ref{p.artur}.

\subsection{Intersections of hyperplane sections}
\label{ss.disjointness}

Now we prove Proposition~\ref{p.artur}. We say that $I \subset
\NN$ is a \emph{$k$-cube} of sides $c_1, \ldots,c_k \in \NN$ based on
$c \in \NN\cup\{0\}$ if $I$ is the set of all $x\in\NN$ that can
be written as $x=c+\sum_i a_i c_i$ with $a_i=0$ or $a_i=1$. We
shall need the following couple of lemmas on $k$-cubes.

\begin{lemma}\label{l.first}
Let $H \subset \CC^d$ be a codimension $1$ subspace, $B:\CC^d \to
\CC^d$ be a linear isomorphism, and $I$ be a $k$-cube for some $1
\leq k \leq d$.  If $H(I)=\cap_{i \in I} B^i(H)$ has codimension
at most $k$ then there exists a subcube $I' \subset I$ and an
integer $l\ge 1$ such that $B^l(H(I'))=H(I')$.
\end{lemma}

\begin{proof}
The proof is by induction on $k$.  The case $k=1$ is easy. Indeed,
the $1$-cube $I=\{c,c+c_1\}$ and so $H(I)=B^c(H)\cap
B^{c+c_1}(H)$. Since $H$ has codimension $1$, if $H(I)$ has
codimension at most $1$ then all the subspaces involved must
coincide:
$$
H(I)=B^c(H)=B^{c+c_1}(H),
$$
and this gives the claim with $l=c_1$ and $I'=\{c\}$. Now assume the
statement holds for $k-1$. Let $I$ be a $k$-cube of sides
$c_1, \ldots, c_k$ based on $c$. Let $I_1$ and $I_2$ be the $(k-1)$-cubes
of sides $c_1, \ldots, c_{k-1}$ based on $c$ and on $c+c_k$,
respectively. Then $I=I_1\cup I_2$. If either $H(I_1)$ or $H(I_2)$
has codimension at most $k-1$, then the conclusion follows from the
induction hypothesis. Otherwise, both $H(I^1)$ and $H(I^2)$ have
codimension at least $k$. Since their intersection $H(I)$ has
codimension at most $k$, they must all coincide:
$$
H(I)=H(I_1)=H(I_2)=B^{c_k}(H(I_1))
$$
and the conclusion follows, with $l=c_k$ and $I'=I_1$.
\end{proof}


\begin{lemma}\label{l.second}
For every $\vep>0$ and $k\ge 1$ there exists $\delta>0$ such that
for all sufficiently large $N\ge 1$ the following holds: for every
$\vep$-dense subset $J$ of $\{0, 1, \ldots, N-1\}$ there exist
$c_1, \ldots, c_k \in \NN$ and  a $\delta$-dense subset $J^k$ of $\{0,
1, \ldots, N-1\}$ such that for every $c\in J^k$ the set $J$
contains the $k$-cube with sides $c_1, \ldots, c_k$ based on $c$.
\end{lemma}

\begin{proof}
The proof is by induction on $k$. Let us start with the case $k=1$.
Let $a_j$, $j= 1, \ldots, \#J$ be the elements of $J$, in increasing
order. By assumption, $\#J\ge\vep N$. Then, clearly,
$$
\frac 1{\#J-1} \sum_{i=1}^{\#J-1} a_{i+1}-a_i \le
\frac{N-1}{\#J-1} \le \frac{2N}{\#J} \le \frac{2}{\vep}
$$
(assume $N$ is large enough so that $\#J \ge\vep N\ge 2$). Then at
least half of these differences are less than twice the average:
there exists $I'\subset\{1, \ldots, \#J-1\}$ with
$\#I'\ge(\#J-1)/2\ge\#J/4$ such that $a_{i+1}-a_i\le 4/\vep$ for
all $i\in I'$. Then there must be some $c_1\ge 1$ and a subset
$I''$ of $I'$ such that
$$
a_{i+1}-a_i = c_1 \quad\text{for all $i\in I''$} \quand \#I''\ge
\frac{\vep\#I'}{4}\ge\frac{\vep\#J}{16}\ge\frac{\vep^2}{16}N.
$$
It follows that $\delta=\vep^2/16$ and $J^1=\{a_i:i\in I''\}$
satisfy the conclusion of the lemma for $k=1$.

Now assume the conclusion holds for $k-1$. Then there exists
$\delta_{k-1}=\delta(\vep)>0$, positive integers
$c_1, \ldots, c_{k-1}$, and a $\delta_{k-1}$-dense subset $J^{k-1}$ of
$\{0, 1, \ldots, N-1\}$ such that $J$ contains a $(k-1)$-cube of
sides $c_1, \ldots, c_{k-1}$ based on every $c\in J^{k-1}$. Applying
case $k=1$ of the lemma with $\delta_{k-1}$ in the place of
$\vep$, we find $\delta=\delta(\vep)>0$, a positive integer $c_k$,
and a $\delta$-dense subset $J^k$ of $\{0, 1, \ldots, N-1\}$ such
that $\{c, c+c_k\} \in J^{k-1}$ for every $c \in J^k$. Then
 $c_1, \ldots, c_k$ and $J^k$ satisfy the conclusion of the lemma.
\end{proof}

We now conclude the proof of the proposition. Fix
$k=\dim\extl{\CC^d}-1$. Assume $N$ is large enough so that
Lemma~\ref{l.second} applies. It follows from the lemma that $J$
contains some $k$-cube $I$. Let $H$ be the geometric hyperplane
corresponding to $V$. If $\cap_{i \in I} B^j(V)\subset\Grass$ is
not empty then $H(I)=\cap_{i\in I} B^i(H)$ has positive dimension,
that is, its codimension in $\extl{\CC^d}$ is at most $k$. So,
Lemma~\ref{l.first} implies that there exists a subcube $I'
\subset I$ and an integer $l\ge 1$ such that $H(I')$ is invariant
under $B^l$. Thus, $\cap_{i \in I'} B^i(V)\subset\Grass$ is
non-empty and invariant under $B^l$. Since all the eigenvalues of
$B$ have different absolute values, for every $\ell$-subspace $W
\subset \CC^d$ we have that $B^j (W)$ converges to a sum of
eigenspaces of $B$ as $j\to\infty$. Since $\cap_{i \in I'} B^j(V)$
is non-empty, invariant, and closed, we conclude that it contains
some sum of eigenspaces of $B$. In particular, $V$ contains a sum
of eigenspaces of $B$, which contradicts the hypothesis. This
contradiction proves Proposition~\ref{p.artur}.

\section{Convergence to a Dirac measure}\label{s.Dirac}

In this section we prove that, for simple cocycles, the limit of
the iterates of any invariant $u$-state $m$ is a Dirac measure on
almost every Grassmannian fiber. Recall that, given any
$\hx\in\hsig$, we denote $\xn=P(\hf^{-n}(\hx))$ for $n\ge 0$.

\begin{proposition} \label{p.toDirac}
If $\hA$ is simple then, for every invariant $u$-state $m$ and
$\hmu$-almost every $\hx\in\hsig$, the sequence
$A^n(\xn)_*m_{\xn}$ converges to a Dirac measure
$\delta_{\xi(\hx)}$ in the fiber $\{x\}\times\Grass$ when
$n\to\infty$.
\end{proposition}

\begin{proof}
In view of Proposition~\ref{p.convergencia1}, we only have to show
that for $\hmu$-almost every $\hx\in\hsig$ there exists some
subsequence $(n_j)_j$ and a point $\xi(\hx)\in\Grass$ such that
\begin{equation}
\label{e.dirac0} A^{n_j}(x^{n_j})_*
m_{x^{n_j}}\to\delta_{\xi(\hx)} \quad\text{when\ } j\to \infty.
\end{equation}
Let $\hp\in\hsig$ be a periodic point, with period $q\ge 1$, and
$\hz\in\hsig$ be a homoclinic point as in
Definition~\ref{d.simple}. Denote $p=P(\hp)$ and $z=P(\hz)$. Let
$[I]=[\iota_0, \ldots, \iota_{q-1}]$ be the cylinder of $\Sigma$
that contains $p$. It is no restriction to assume that $z\in [I]$:
this may always be achieved replacing $\hz$ by some
$\hf^{-qi}(\hz)$ which, clearly, does not affect the conditions in
Definition~\ref{d.simple}.

\begin{figure}[phtb]
\begin{center}
\psfrag{pp}{$\hp$} \psfrag{zz}{$\hz$} \psfrag{xx}{$\hx$}
\psfrag{p}{$p$} \psfrag{z}{$z$} \psfrag{zk}{$z^{qk}$}
\psfrag{x0}{$x_0$} \psfrag{x1}{$x^{n_j+qk}$}
\psfrag{x2}{$x^{n_j}$} \psfrag{f1}{$\hf^{qk}$}
\psfrag{f2}{$\hf^{n_j}$}
\includegraphics[height=1.8in]{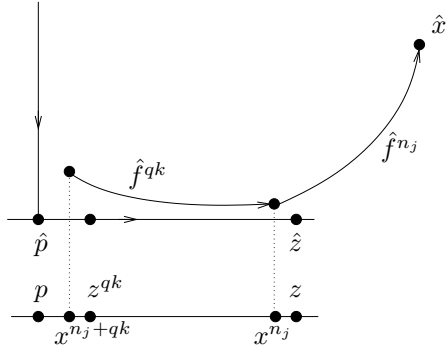}
\caption{\label{f.case1} Proof of Proposition~\ref{p.toDirac}:
case $\xi(\hz)$ not in $\ker Q$}
\end{center}
\end{figure}

For $\hmu$-almost every $\hx\in\hsig$ there exists a sequence
$(n_j)_j$ such that $\hf^{-n_j}(\hx)$ converges to $\hz$. That is
because $\hmu$ is ergodic and positive on open sets. Let $k\ge 1$ be
fixed. 
From Proposition~\ref{p.convergencia1} 
we conclude that
$$
\begin{aligned}
\lim_{j\to\infty} A^{n_j}(x^{n_j})_* m_{x^{n_j}}
 & = \lim_{j\to\infty} A^{n_j+qk}(x^{n_j+qk})_* m_{x^{n_j+qk}}
 \\
 & = \lim_{j\to\infty} A^{n_j}(x^{n_j})_* A^{qk} (x^{n_j+qk})_*m_{x^{n_j+qk}}.
\end{aligned}
$$
Note that $x^{n_j+qk}$ converges to $z^{qk}$ when $j\to\infty$. 
See Figure~\ref{f.case1}. Then, by Proposition~\ref{p.continuity},
the probability $m_{x^{n_j+qk}}$ converges to $m_{z^{qk}}$ when
$j\to\infty$. So, since $A$ is continuous,
$$
A^{qk}(x^{n_j+qk})_* m_{x^{n_j+qk}} \to A^{qk}(z^{qk})_* m_{z^{qk}}
\quad\text{when\ } j\to\infty.
$$
Since the space of quasi-projective maps is compact, up to replacing
$(n_j)_j$ by a subsequence we may suppose that $A^{n_j}(x^{n_j})$
converges to some quasi-projective map $Q$ on $\Grass$. By
Lemma~\ref{l.kernel}, the kernel of $Q$ is contained in some
hyperplane of $\Grass$. Hence, by Proposition~\ref{p.noatoms}, the
subspace $\ker Q$ has zero measure relative to $A^{qk}(z^{qk})_*
m_{z^{qk}}$. So, we may apply Lemma~\ref{l.topology2} to conclude
that
\begin{equation}\label{e.dirac2}
\lim_{j\to\infty} A^{n_j}(x^{n_j})_* m_{x^{n_j}}
= Q_* A^{qk}(z^{qk})_* m_{z^{qk}},
\end{equation}
for any $k\ge 1$ (in particular, the latter expression does not
depend on $k$).

Now, the pinching condition (p) in Definition~\ref{d.simple}
implies that $A^q(p)$ has $\ell$ eigenvalues that are strictly
larger, in norm, than all the other ones. Denote by
$\xi(\hp)\in\Grass$ the sum of the eigenspaces corresponding to
those largest eigenvalues, and define
$\xi(\hz)=\Hpz\cdot\xi(\hp)$.

\begin{lemma}\label{l.whatmeasure}
The sequence $A^{qk}(z^{qk})_* m_{z^{qk}}$ converges to
$\delta_{\xi(\hz)}$ when $k\to\infty$.
\end{lemma}

\begin{proof}
Using the relations $A^{qk}(p)^{-1} = \hA^{-qk}(\hp)$ and
                    $A^{qk}(z^{qk})^{-1}=\hA^{-qk}(\hz)$, we find that
$$
A^{qk}(z^{qk})_* m_{z^{qk}}
= \left(\hA^{-qk}(\hz)^{-1} \cdot \hA^{-qk}(\hp)\right)_*
A^{qk}(p)_* m_{z^{qk}}.
$$
By the Definition~\ref{d.suholonomies} of unstable holonomies,
$\hA^{-qk}(\hz)^{-1} \hA^{-qk}(\hp)$ converges to $H_{\hp,\hz}$
when $k\to\infty$. Observe also that $A^{qk}(p)_* m_{z^{qk}}$
converges to the Dirac measure at $\xi(\hp)\in\Grass$ when
$k\to\infty$. That is because $m_{z^{qk}}$ converges to $m_p$, by
Proposition~\ref{p.continuity}, and $m_p$ gives zero weight to the
hyperplane section defined by the sum of the eigenspaces of
$A^q(p)$ complementary to $\xi(\hp)$, by
Proposition~\ref{p.noatoms}. It follows that $A^{qk}(z^{qk})_*
m_{z^{qk}}$ converges to $(\Hpz)_* \, \delta_{\xi(\hp)} =
\delta_{\xi(\hz)}$ when $k\to\infty$, as stated in the lemma.
\end{proof}

Suppose, for the time being, that $\xi(\hz)$ is in the domain
$\Grass\setminus \ker Q$ of the quasi-projective map $Q$. From
Lemma~\ref{l.topology2} we get that
$$
Q_* A^{qk}(p)_* m_{z^{qk}} \to Q_* \delta_{\xi(\hz)}
 = \delta_{\xi(\hx)}
$$
when $k\to\infty$, where $\xi(\hx)=Q(\xi(\hz))$. Combined with the
relation \eqref{e.dirac2}, this gives that
$A^{n_j}(x^{n_j})_*m_{x^{n_j}}$ converges to the Dirac measure
$\delta_{\xi(\hx)}$ when $j\to\infty$. This proves
\eqref{e.dirac0} and Proposition~\ref{p.toDirac} in this case.

\begin{figure}[phtb]
\begin{center}
\psfrag{pp}{$\hp$} \psfrag{zz}{$\hz$} \psfrag{xx}{$\hx$}
\psfrag{yy}{$\hy$} \psfrag{p}{$p$} \psfrag{z}{$z$}
\psfrag{w}{$\hf^l(\hz)$} \psfrag{x2}{{\footnotesize $x^{n_j}$}}
\psfrag{x1}{{\footnotesize $x^{n_j+qk}$}}
\psfrag{x3}{{\footnotesize $y^{n_j+qk+l}$}}
\psfrag{x4}{{\footnotesize $y^{n_j+qk+l+qm}$}} \psfrag{f1}{{\tiny
$\hf^{qk}$}} \psfrag{f2}{{\tiny $\hf^{n_j}$}} \psfrag{f3}{{\tiny
$\hf^{l}$}} \psfrag{f4}{{\tiny $\hf^{qm}$}}
\includegraphics[height=2in]{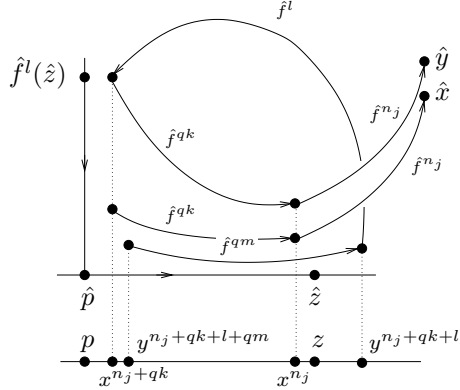}
\caption{\label{f.case2} Proof of Proposition~\ref{p.toDirac}:
avoiding $\ker Q$}
\end{center}
\end{figure}

Next, we show that one can always reduce the proof to the previous
case. Let $l\ge 1$ be as in Definition~\ref{d.simple}. For each
$j$ much larger than $k$, let $m_j=n_j+qk+l$ and $\hy=\hy(j,k)$ be
defined by
\begin{equation}\label{eq.yyy}
\hf^{-n_j-qk}(\hy) \in \Wsloc(\hf^{-n_j-qk}(\hx)) \cap
\Wuloc(\hf^l(\hz)).
\end{equation}
See Figure~\ref{f.case2}.
By construction, $y^{m_j+qm}$ is sent to $x^{n_j+qk}$ by the map $f^{l+qm}$.
Hence, using Proposition~\ref{p.convergencia1},
\begin{equation}\label{eq.lll}
\begin{aligned}
\lim_{j\to\infty} A^{n_j}(x^{n_j})_*m_{x^{n_j}} & =
\lim_{j\to\infty} A^{n_j+qk}(x^{n_j+qk})_*m_{x^{n_j+qk}}
\\
& = \lim_{j\to\infty} A^{m_j+qm}(y^{m_j+qm})_*m_{y^{m_j+qm}}.
\end{aligned}
\end{equation}
for any fixed $k$ and $m$. We are going to prove that the limit is
indeed a Dirac measure. For this, let $\hw=\hw(k)$ be defined by
\begin{equation}\label{eq.www}
\hf^l(\hw)\in \Wsloc(\hf^{-qk}(\hz))\cap \Wuloc(\hf^l(\hz)).
\end{equation}
Notice that $\hw$ is in $\Wuloc(\hz)=\Wuloc(\hp)$. Let $k$ and $m$
be fixed, for the time being. As $j\to\infty$, the sequence
$\hf^{-n_j-qk}(\hx)$ converges to $\hf^{-qk}(\hz)$ and so,
combining \eqref{eq.yyy} and \eqref{eq.www}, the sequence
$\hf^{-n_j-qk}(\hy)$ converges to $\hf^l(\hw)$. It follows that
$y^{m_j}$ converges to $w=P(\hw)$, and so
$$
A^{m_j}(y^{m_j})=A^{n_j}(x^{n_j}) A^{qk+l}(y^{m_j})
$$
converges to $\tilde Q = Q \circ A^{qk+l}(w)$ in the space of
quasi-projective maps, as $j\to\infty$. Define
$\xi(\hw)=\Hpw\cdot\xi(\hp)$. The key observation is

\begin{lemma}\label{l.keyobservation}
Assuming $k$ is large enough, $\xi(\hw)$ is not contained in $\ker
\tilde Q$.
\end{lemma}

\begin{proof}
From the definitions of $\tilde Q$ and $\hw$ we get that
$$
\ker \tilde Q
 = A^{kq+l}(w)^{-1} \cdot \ker Q
 = A^l(w)^{-1} \cdot A^{qk}(z^{qk})^{-1} \cdot \ker Q.
$$
By the invariance property of unstable holonomies, we have
$$
A^{qk}(z^{qk})^{-1}
 = \hA^{-qk}(\hz)
 = H_{\hp,\hf^{-qk}(\hz)} \cdot \hA^{-qk}(\hp) \cdot H_{\hz,\hp}.
$$
So, the previous equality may be rewritten as
$$
\ker \tilde Q
 = \hA^l(\hw)^{-1} \cdot H_{\hp,\hf^{-qk}(\hz)} \cdot
     \hA^{-qk}(\hp) \cdot H_{\hz,\hp} \cdot \ker Q.
$$
Notice that $\hf^{-qk}(\hz)$ converges to $\hp$ and so, by
\eqref{eq.www}, the point $\hw$ converges to $\hz$, as
$k\to\infty$. By the continuity of the cocycle and the holonomies,
it follows that $H_{\hp,\hf^{-qk}(\hz)}$ converges to the identity
and $\hA^l(\hw)$ converges to $\hA^l(\hz)$, as $k$ goes to
$\infty$. By Lemma~\ref{l.kernel}, the kernel of $Q$ is contained
in some hyperplane section of $\Grass$. Then the same is true for
$H_{\hz,\hp}\cdot\ker Q$: it is contained in the set of all
$\ell$-dimensional subspaces that intersect the
$(d-\ell)$-dimensional subspace $\pi_v(\upsilon)$ associated to
some $(d-\ell)$-vector $\upsilon$. Since all eigenvalues of
$\hA^q(\hp)$ have distinct absolute values, the backward iterates
of $\pi_v(\upsilon)$ under $\hA^{q}(\hp)$ converge to some
$(d-\ell)$-dimensional sum $\pi_v(\eta)$ of eigenspaces of
$\hA^{q}(\hp)$. It follows that, as $k\to\infty$, the sequence
$\hA^{-qk}(\hp)\cdot H_{\hz,\hp} \cdot\ker Q$ converges to some
subset $V_0$ of the hyperplane section $V$ defined by $\eta$.
Combining these two observations we get that, as $k\to\infty$,
\begin{equation}\label{eq.convker}
\ker \tilde Q \to \hA^l(\hz)^{-1} (V_0)
    \subset \hA^l(\hz)^{-1}(V).
\end{equation}
It is easy to see that $\xi(\hz)$ does not belong to
$\hA^l(\hz)^{-1}(V)$: otherwise,
$$
\hA^l(\hz)\cdot\xi(\hz)
 = \hA^l(\hz) \cdot H_{\hp,\hz} \cdot \xi(\hp)
 = \psi_{p,z} \cdot \xi(\hp)
$$
would intersect $\pi_v(\eta)$ and, since $\xi(\hp)$ and
$\pi_v(\eta)$ correspond to sums of eigenspaces with complementary
dimensions, that would contradict the twisting condition in
Definition~\ref{d.simple}. Using \eqref{eq.convker} and the fact
that $\xi(\hw)$ converges to $\xi(\hz)$ when $k\to\infty$, we
deduce that $\xi(\hw)$ is not in $\ker\tilde Q$ if $k$ is large
enough, as claimed.
\end{proof}

We can now finish the proof of Proposition~\ref{p.toDirac}. The
arguments are the same as in the previous case, with $n_j$ and $z$
replaced by $m_j=n_j+qk+l$ and $w$, respectively, and $qm$ in the
role of $qk$. Indeed, from \eqref{eq.yyy} and \eqref{eq.www} we get
that $\hf^{-m_j}(\hy)$ converges to $\hw$ as $j\to\infty$.
Consequently, $A^{qm}(y^{m_j+qm})$ converges to $A^{qm}(w^{qm})$
and, using also Proposition~\ref{p.continuity}, $m_{y^{m_j+qm}}$
converges to $m_{w^{qm}}$ as $j\to\infty$. So, in view of
\eqref{eq.lll},
\begin{equation*}\begin{aligned}
\lim_{j\to\infty} A^{n_j}(x^{n_j})_*m_{x^{n_j}}
 & = \lim_{j\to\infty} A^{m_j+qm}(y^{m_j+qm})_* m_{y^{m_j+qm}}
 \\ & = \tilde Q_* A^{qm}(p)_*m_{w^{qm}}
\end{aligned}\end{equation*}
for any $m\ge 1$, which is an analogue of \eqref{e.dirac2}. By
Proposition~\ref{p.continuity}, the measure $m_{w^{qm}}$ converges
to $m_p$ as $m\to\infty$. By Proposition~\ref{p.noatoms}, the
measure $m_p$ gives zero weight to the hyperplane section defined
by the sum of the eigenspaces complementary to $\xi(\hp)$.
Therefore, just as in Lemma~\ref{l.whatmeasure}, we conclude that
$A^{qm}(\hw^{qm})_*m_{\hw^{qm}}$ converges to $\delta_{\xi(\hw)}$
when $m\to\infty$. Hence, fixing $k$ as in
Lemma~\ref{l.keyobservation} and using Lemma~\ref{l.topology2},
$$
\lim_{m\to\infty}\tilde Q_* A^{qm}(p)_*m_{w^{qm}} =
\delta_{\xi(\hx)},
$$
where $\xi(\hx)=\tilde Q_* \delta_{\xi(\hw)}$. This shows that
$\lim_{j\to\infty} A^{n_j}(x^{n_j})_*m_{x^{n_j}} =
\delta_{\xi(\hx)}$. Now the proof of Proposition~\ref{p.toDirac}
is complete.
\end{proof}

In the next proposition we summarize some consequences of the
previous results that are needed for the next section:

\begin{proposition}\label{p.summary}
Suppose that $\hA$ is simple. Then there exists a measurable section
$\xi:\hsig\to\Grass$ such that, on a full $\hmu$-measure subset of
$\hsig$,
\begin{enumerate}
\item $\xi$ is invariant under the cocycle and under unstable
holonomies: $\hA(\hx)\xi(\hx)=\xi(\hf(\hx))$ and $\xi(\hy) =
H^u_{\hx,\hy} \cdot \xi(\hx)$ for $\hx$ and $\hy$ in the same
local unstable set \item for any compact set $\Gamma\subset\hsig$,
the eccentricity $\cE(\ell,\hA^n(\hf^{-n}(\hx)))\to\infty$, and
the image under $\hA^n(\hf^{-n}(\hx))$ of the $\ell$-subspace most
expanded by $\hA^n(\hf^{-n}(\hx))$ converges to $\xi(\hx)$,
restricted to the subsequence of iterates
$\hf^{-n}(\hx)\in\Gamma$.
\end{enumerate}
\end{proposition}

\begin{proof}
From Corollary~\ref{c.gauche1} and Proposition~\ref{p.toDirac} we
get that the conditional probabilities of the original measure $\hm$
along the Grassmannian fibers coincide with the Dirac measures
$\delta_{\xi(\hx)}$ almost everywhere. Since $\hm$ is an invariant
$u$-state, it follows that $\xi$ is almost everywhere invariant
under the cocycle and under the unstable holonomies, as stated in
part 1 of the proposition.

Part 2 follows from Proposition~\ref{p.eccentricity}, with
$\cN=\{m_{P(\hx)} : \hx\in\Gamma\}$, $\nu_n=m_{\xn}$,
$L_n=A^n(\xn)=\hA^n(\hf^{-n}(\hx))$, and $\xi=\xi(\hx)$.
By Proposition~\ref{p.continuity}, the family is $\cN$ is
\weak compact if $\Gamma$ is compact. It follows that the
eccentricity $\cE(\ell,\hA^n(\hf^{-n}(\hx)))$ tends to infinity,
and the image under $\hA^n(\hf^{-n}(\hx))$ of the subspace most
expanded by $\hA^n(\hf^{-n}(\hx))$ converges to $\xi(\hx)$, as
claimed.
\end{proof}

\begin{remark}\label{r.doesnotmatter}
In Section~\ref{ss.reduction} we replaced the original cocycle
$\hA$ by another one conjugate to it,
$$
\hpA(\hx)
 = H^s_{\hf(\hx),\pphi(\hf(\hx))} \cdot \hA(\hx) \cdot H^s_{\pphi(\hx),\hx}
 = H^s_{\hf(\pphi(\hx)),\pphi(\hf(\hx))} \cdot \hA(\pphi(\hx)),
$$
which is constant on local unstable sets and, consequently, whose
stable holonomies are trivial. The statement of
Proposition~\ref{p.summary} is not affected by such substitution.
Indeed, if $\xi$ is an invariant section for $\hpA$ as in the
proposition, then
$$
H^s_{\pphi(\hx),\hx} \xi(\hx)
$$
is an invariant section for $\hA$, and it is invariant also under
the corresponding unstable holonomies. In addition,
$$
\hA^n(\hf^{-n}(\hx))
 = H^s_{\pphi(\hx),\hx}
   \cdot (\hpA)^n(\hf^{-n}(\hx))
   \cdot H^s_{\hf^{-n}(\hx),\pphi(\hf^{-n}(\hx))}.
$$
Considering only iterates in a compact set, the corresponding
conjugating isomorphisms $H^s$
belong to a bounded family. Hence, the claims in part 2 of
Proposition~\ref{p.summary} hold for $\hA$ if and only if they
hold for $\hpA$.
\end{remark}

\section{Proof of the main theorem}\label{s.proof}

We are going to show that $\hx \mapsto \xi(\hx)\in\Grass$
corresponds to the sum of the Oseledets subspaces of the cocycle
associated to the $\ell$ largest (strictly) Lyapunov exponents. In
particular, $\xi(\hx)$ is uniquely defined almost everywhere. This
will also prove that the invariant $u$-state is unique if the
cocycle is simple.

The first step is to exhibit the sum $\eta(\hx)$ of the subspaces
associated to the remaining Lyapunov exponents. This is done in
Section~\ref{ss.adjoint}, through applying the previous theory to
the adjoint cocycle. Then, in Section~\ref{ss.directionofmaximum}
we use the second part of Proposition~\ref{p.summary} to show that
vectors along $\xi(\hx)$ are more expanded than those along $\eta(\hx)$.

\subsection{Adjoint cocycle}\label{ss.adjoint}
Let $\cdot$ be a Hermitian form on $\CC^d$, that is, a complex
$2$-form $(u,v) \mapsto u\cdot v$ which is linear on the first
variable and satisfies $u\cdot v = \overline{v \cdot u}$ for every
$u$ and $v$.
The \emph{adjoint} of a linear operator $L:\CC^d\to\CC^d$ relative
to the Hermitian form is the linear operator $L^*:\CC^d\to\CC^d$
defined by
$$
L^*(u) \cdot v = u \cdot L(v)
 \quad\text{for every $u$ and $v$ in $\CC^d$.}
$$
The matrix of $L^*$ in any orthonormal basis for the Hermitian form
is the conjugate transpose of the matrix of $L$ in that basis:
$L^*_{i,j} = \overline{L}_{j,i}$. The eigenvalues of $L^*$ are the
conjugates of the eigenvalues of $L$, and the operator norms of the
two operators coincide: $\|L^*\|=\|L\|$.

Let $\hB(\hx):\CC^d\to\CC^d$ be defined by $\hB(\hx) =
\hA(\hf^{-1}(\hx))^*$ or, equivalently,
\begin{equation}\label{eq.adjoint2}
\hB(\hx) u \cdot v = u \cdot \hA(\hf^{-1}(\hx)) v \quad\text{for
every $u$ and $v \in\CC^d$.}
\end{equation}
Consider the linear cocycle defined over $\hf^{-1}$ by
$$
\hFB:\hsig\times\CC^d\to\hsig\times\CC^d,
 \quad (\hx,u)\mapsto (\hf^{-1}(\hx), \hB(\hx) u),
$$
as well as the induced Grassmannian cocycle. Notice that
\begin{equation*}\begin{aligned}
\hB^n(\hx)
 & = \hA(\hf^{-n}(\hx))^* \cdots \hA(\hf^{-2}(\hx))^* \hA(\hf^{-1}(\hx))^*
 = \hA^n(\hf^{-n}(\hx))^*.
\end{aligned}\end{equation*}
The choice of the Hermitian form is not important: different choices
yield cocycles that are conjugate. For convenience, we fix once and
for all such that eigenvectors of $\hA^q(\hp)$ form an orthonormal
basis.

The integrability condition in the Oseledets theorem holds for $\hB$
if and only if it holds for $\hA$, because
$\|\hB(\hx)\|=\|\hA(\hf^{-1}(\hx))\|$ and the measure $\hmu$ is
invariant under $\hf$.
%
It is easy to check that the previous results apply to the cocycle
defined by $\hB$. To begin with, our hypotheses on the dynamics
(Section~\ref{sh.dynamics}) and on the invariant measure
(Section~\ref{sh.measure}) are, evidently, symmetric under time
reversion. The hypotheses on the cocycle (Section~\ref{sh.cocycle})
are also clearly satisfied: a simple calculation shows that $\hB$
admits stable and unstable holonomies given by
\begin{equation}\label{eq.holonomyreverse}
H_{\hx,\hy}^{u,\hB} = \Big(H_{\hy,\hx}^{s,\hA}\Big)^*
 \quand
H_{\hx,\hy}^{s,\hB} = \Big(H_{\hy,\hx}^{u,\hA}\Big)^*.
\end{equation}

\begin{lemma}\label{l.iff}
$\hB$ is simple for $\hf^{-1}$ if and only if $\hA$ is simple for
$\hf$.
\end{lemma}

\begin{proof}
Let $\hp$ be a periodic point of $\hf$. For any orthonormal basis of
$\CC^d$, the matrix of $\hB^q(\hp)=\hA^q(\hp)^*$ is the conjugate
transpose of the matrix of $\hA^q(\hp)$, and the eigenvalues of the
former are the complex conjugates of the eigenvalues of the latter.
Hence, the pinching condition in Definition~\ref{d.simple} holds for
any of them if and only if it holds for the other. Next, notice that
$\hz$ is a homoclinic point for $\hf$ if and only if
$\hw=\hf^l(\hz)$ is a homoclinic point for the inverse: $\hz \in
\Wuloc(\hp,\hf)$ and $\hf^l(\hz) \in \Wsloc(\hp,\hf)$ if and only if
$\hf^{-l}(\hw) \in \Wsloc(\hp,\hf^{-1})$ and $\hw \in
\Wuloc(\hp,\hf^{-1})$. We have chosen the Hermitian form in such a
way that eigenvectors of $\hA^q(\hp)$ form an orthonormal basis.
Then the matrix of $\hB^l(\hw)=\hA^l(\hz)^*$ in this basis is the
conjugate transpose of the matrix of $\hA^l(\hz)$, and so the
algebraic minors of the former are the complex conjugates of the
algebraic minors of the latter. Thus, the twisting condition in
Definition~\ref{d.simple} holds for $\hB$ if and only if it holds
for $\hA$.
\end{proof}

This ensures that the previous results do apply to $\hB$. From
Proposition~\ref{p.summary} we obtain that
\begin{itemize}
 \item[(i)] there exists a section $\xi^*: \hsig \to \Grass$ which is
invariant under the cocycle $\hFB$ and under the unstable holonomies of $\hB$
 \item[(ii)] given any compact $\Gamma\subset\hsig$, restricted to the
subsequence of iterates $\hf^n(\hx)$ in $\Gamma$, the eccentricity
$\cE(\ell, \hB^n(\hf^n(\hx)))=\cE(\ell, \hA^n(\hx))$ goes to
infinity and the image $\hB^n(\hf^n(\hx))\zeta^a_n(\hf^n(\hx))$ of
the $\ell$-subspace $\zeta^a_n(\hf^n(\hx))$ most expanded by
$\hB^n(\hf^n(\hx))$ tends to $\xi^*(\hx)$ as $n\to\infty$.
\end{itemize}
Let us show that $\xi(\hx)$ is outside the hyperplane section
orthogonal to $\xi^*(\hx)$:

\begin{lemma}\label{l.xi*xi}
For $\hmu$-almost every $\hx$, the subspace $\xi(\hx)$ is transverse
to the orthogonal complement of $\xi^*(\hx)$.
\end{lemma}

\begin{proof}
Recall, from Section~\ref{ss.reduction} and
Remark~\ref{r.doesnotmatter}, that we may take the stable holonomies
of $\hA$ to be trivial. Then, by \eqref{eq.holonomyreverse}, the
unstable holonomies of $\hB$ are also trivial. So, the fact that
$\xi^*$ is invariant under unstable holonomies just means that it is
constant on local unstable sets of $\hf^{-1}$, that is, on local
stable sets of $\hf$. Then the same is true about the orthogonal
complement of $\xi^*(\hx)$. In other words, the hyperplane section
of $\Grass$ orthogonal to $\xi^*(\hx)$ depends only on $x=P(\hx)$.
Denote it as $H_x$. Using Proposition~\ref{p.noatoms} and then
Proposition~\ref{p.summary}, we obtain
$$
0 = m_x(H_x) = \int \delta_{\xi(\hx)} (H_x) \, d\hmu_x(\hx) =
\hmu_x\big(\{\hx\in W_{loc}^s(x) : \xi(\hx)\in H_x\}\big),
$$
for $\mu$-almost every $x$. Consequently, $\hmu\big(\{\hx\in\hsig :
\xi(\hx)\in H_x\}\big)=0$. This means that, for almost every point,
the subspace $\xi(\hx)$ intersects the orthogonal complement of
$\xi^*(\hx)$ at the origin only, which is precisely the claim in the
lemma.
\end{proof}

Let $\eta(\hx)\in\operatorname{Grass}(d-\ell,d)$ denote the
orthogonal complement of $\xi^*(\hx)$. Recall that $\xi$ and $\xi^*$
are invariant under the corresponding cocycles:
$$
\hA(\hx)\xi(\hx) = \xi(\hf(\hx))
 \quad\text{and}\quad
\hB(\hx)\xi^*(\hx) = \xi^*(\hf^{-1}(\hx))
$$
$\hmu$-almost everywhere. The latter implies that $\eta(\hx)$ is
also invariant under $\hA$. According to Lemma~\ref{l.xi*xi}, we
have $\CC^d = \xi(\hx) \oplus \eta(\hx)$ at almost every point. We
want to prove that the Lyapunov exponents of $\hA$ along $\xi$ are
strictly bigger than those along $\eta$. To this end, let
$$
\xi(\hx) = \xi^1(\hx)\oplus \cdots \oplus \xi^u(\hx) \quand
 \eta(\hx) = \eta^s(\hx) \oplus \cdots \oplus \eta^1(\hx)
$$
be the Oseledets decompositions of $\hA$ restricted to the two
invariant subbundles. Take the factors to be numbered in such a way
that $\xi^u$ corresponds to the smallest Lyapunov exponent among all
$\xi^i$, and $\eta^s$ corresponds to the largest Lyapunov exponent
among all $\eta^j$. Denote $d_u=\dim \xi^u$ and $d_s=\dim \eta^s$, and
let $\lambda_u$ and $\lambda_s$ be the Lyapunov exponents associated
to these two subbundles, respectively.

\subsection{Direction of maximum expansion}
\label{ss.directionofmaximum}

Given a linear map $L:\CC^d\to\CC^d$ and a subspace $V$ of $\CC^d$,
we denote by $\det(L,V)$ the determinant of $L$ along $V$, defined
as the quotient of the volumes of the parallelograms determined by
$\{Lv_1, \dots, Lv_s\}$ and $\{v_1, \dots, v_s\}$, respectively, for
any basis $v_1, \dots, v_s$ of $V$. Then we define, for each $n\ge
1$,
\begin{equation}\label{eq.Dn}
\Delta^n(\hx) =
          \frac{\det(\hA^n(\hx), \xi^u(\hx))^{1/d_u}}
               {\det(\hA^n(\hx), W(\hx))^{1/(d_u+d_s)}}
               \quad \text{where}\quad W(\hx)= \xi^u(\hx)\oplus\eta^s(\hx)\, .
\end{equation}
According to the theorem of Oseledets~\cite{Os68},
$$
\frac{1}{n} \log \det(\hA^n(\hx), \xi^u(\hx))
 \to d_u \, \lambda_u
 \quand
\frac{1}{n} \log \det(\hA^n(\hx), W(\hx))
 \to d_u \, \lambda_u + d_s \, \lambda_s.
$$
Consequently,
\begin{equation}\label{eq.diferenca}
\lim_{n\to\infty} \frac{1}{n} \log \Delta^n(\hx) = \frac{d_s}{d_u+d_s}
(\lambda_u-\lambda_s).
\end{equation}
So, to prove that $\lambda_u$ is strictly larger than $\lambda_s$
we must show that $\log\Delta^n$ goes linearly to infinity at
almost every point. The main step is

\begin{proposition}\label{p.Delta0}
For any compact set $\Gamma\subset\hsig$ and for $\hmu$-almost every
$\hx\in\hsig$,
$$
\lim_{n\to\infty} \Delta^n(\hx) = + \infty
$$
restricted to the subsequence of values of $n$ for which
$\hf^n(\hx)\in\Gamma$.
\end{proposition}

\begin{proof}
Let $\xi^a_n(\hx)=\hB^n(\hf^n(\hx))\zeta^a_n(\hx)$ be the image of
the $\ell$-dimensional subspace most expanded by
$\hB^n(\hf^n(\hx))=\hA^n(\hx)^*$. Equivalently, $\xi^a_n(\hx)$ is
the $\ell$-dimensional subspace most expanded by $\hA^n(\hx)$.
Throughout, we consider only the values of $n$ for which
$\hf^{n}(\hx)\in\Gamma$. Then we may use property (ii) in
Section~\ref{ss.adjoint}: the eccentricity
$$
E_n
 =\cE(\ell, \hB^n(\hf^n(\hx)))
 =\cE(\ell, \hA^n(\hx))
$$
tends to infinity, and $\xi^a_n(\hx)$ tends to $\xi^*(\hx)$, as
$n\to\infty$. In view of Lemma~\ref{l.xi*xi}, the latter fact
implies that the subspace $\xi(\hx)$ is transverse to the orthogonal
complement of $\xi^a_n(\hx)$, with angle uniformly bounded from zero
for all large $n$. Let us consider the orthogonal splitting
$$
\CC^d = \xi^a_n(\hx) \oplus \xi^a_n(\hx)^\perp.
$$
Let $\xi^u_n(\hx)\subset\xi^a_n(\hx)$ be the image of the subspace
$\xi^u(\hx)\subset\xi(\hx)$ under the orthogonal projection. We
claim that
\begin{equation}\label{eq.gap1}
\det\big(\hA^n(\hx) \mid \xi_n^u(\hx)\big) \le C_1
\det\big(\hA^n(\hx) \mid \xi^u(\hx)\big),
\end{equation}
for some constant $C_1$ independent of $n$. To see this, observe
that any basis $\alpha$ of $\xi^u_n(\hx)$ may be lifted to a basis
$\beta$ of $\xi^u(\hx)$. This operation increases the volume of the
corresponding parallelogram by, at most, some factor $C_1$ that
depends only on a bound for the angle between $\xi(\hx)$ and the
orthogonal complement of $\xi^a_n(\hx)$. Note also that, the
$\hA^n(\hx)$-images of $\xi^a_n(\hx)$ and $\xi^a_n(\hx)^\perp$ are
orthogonal to each other, because $\xi^a_n(\hx)$ is the
$\ell$-subspace most expanded by $\hA^n(\hx)$. Hence, the
$\hA^n(\hx)$-image of $\alpha$ may be obtained from the
$\hA^n(\hx)$-image of $\beta$ by orthogonal projection, an operation
that can only decrease the volume of the parallelogram. Combining
these observations, we get \eqref{eq.gap1}. Next, let $\eta_n^s(\hx)$
be the subspace of $\xi^a_n(\hx)^\perp$ characterized by
$$
W(\hx)
 = \xi^u(\hx) \oplus \eta^s(\hx)
 = \xi^u(\hx) \oplus \eta^s_n(\hx).
$$
Equivalently, $\eta_n^s(\hx)$ is the projection of $\eta^s(\hx)$
to the orthogonal complement of $\xi^a_n(\hx)$ along the direction
of $\xi(\hx)$. Since the angle between $\xi^u(\hx)$ and
$\eta^s_n(\hx)$ is bounded from zero,
\begin{equation}\label{eq.gap2}
\det\big(\hA^n(\hx),W(\hx)\big) \le C_2
\det\big(\hA^n(\hx),\xi^u(\hx)\big)
\det\big(\hA^n(\hx),\eta^s_n(\hx)\big)
\end{equation}
where the constant $C_2$ is independent of $n$. Furthermore,
\begin{equation*}\begin{aligned}
\det\big(\hA^n(\hx),\eta^s_n(\hx)\big)
 & \le \|\hA^n(\hx) \mid \eta_n^s(\hx)\|^{d_s}
 \le \|\hA^n(\hx) \mid \xi^a_n(\hx)^\perp\|^{d_s} \\
\det\big(\hA^n(\hx),\xi^u_n(\hx)\big)
 & \ge m\big(\hA^n(\hx) \mid \xi_n^u(\hx)\big)^{d_u}
\ge m\big(\hA^n(\hx) \mid \xi^a_n(\hx)\big)^{d_u}
\end{aligned}\end{equation*}
because $\eta_n^s(\hx)\subset\xi^a_n(\hx)^\perp$ and
$\xi_n^u(\hx)\subset\xi^a_n(\hx)$. Consequently,
\begin{equation}\label{eq.gap3}
E_n
 = \frac{m(\hA^n(\hx) \mid \xi^a_n(\hx))}
        {\|\hA^n(\hx) \mid \xi^a_n(\hx)^\perp\|}
 \le \frac{\det\big(\hA^n(\hx),\xi^u_n(\hx)\big)^{1/d_u}}
          {\det\big(\hA^n(\hx),\eta^s_n(\hx)\big)^{1/d_s}} \,.
\end{equation}
From \eqref{eq.gap1}--\eqref{eq.gap3} we obtain
$$
\det\big(\hA^n(\hx),W(\hx)\big) \le C E_n^{-d_s}
\det\big(\hA^n(\hx),\xi^u(\hx)\big)^{1+d_s/d_u}
$$
with $C=C_1^{s/u} C_2$. Consequently,
$$
\Delta^n(\hx) =
          \frac{\det(\hA^n(\hx), \xi^u(\hx))^{1/d_u}}
               {\det(\hA^n(\hx), W(\hx))^{1/(d_u+d_s)}}
              \ge \big(C^{-1} E_n^s\big)^{1/d_u+d_s}
$$
and this goes to infinity when $n\to\infty$. The proof of the
proposition is complete.
\end{proof}

Now we are ready for the proof of Theorem~\ref{main}. Fix any
compact set $\Gamma\subset\hsig$ such that $\hmu(\Gamma)>0$.
By Poincar\'e recurrence, the first return map
$$
g:\Gamma\to\Gamma, \quad g(\hx) = \hf^{r(\hx)}(\hx)
$$
is well defined on a full $\hmu$-measure subset of $\hsig$.
The normalized restriction $\hmu/\hmu(\Gamma)$ of the measure
$\hmu$ to $\Gamma$ is invariant and ergodic for $g$. Moreover,
$\hFA$ induces a linear cocycle
$$
G:\Gamma\times\CC^d\to\Gamma\times\CC^d,
 \quad G(\hx,v)=(g(\hx),\cG(\hx)v)
$$
where $\cG(\hx)=\hA^{r(\hx)}(\hx)$. Clearly, this cocycle
preserves the subbundles $\xi(\hx)$ and $\eta(\hx)$, as well
as their Oseledets decompositions
$$
\xi(\hx) = \xi^1(\hx)\oplus \cdots \oplus \xi^u(\hx) \quand
 \eta(\hx) = \eta^s(\hx) \oplus \cdots \oplus \eta^1(\hx).
$$
It is also clear (see Section~\ref{ss.inducing}) that the Lyapunov
exponents of $G$ with respect to $\hmu/\hmu(\Gamma)$ are the
products of the exponents of $\hFA$ by the average return time
$1/\hmu(\Gamma)$.

Thus, to show that $\lambda_u>\lambda_s$ it suffices to prove the
corresponding statement for $G$. Define
$$
\cD^k(\hx) =
          \frac{\det(\cG^k(\hx), \xi^u(\hx))^{1/d_u}}
               {\det(\cG^k(\hx), W(\hx))^{1/(d_u+d_s)}}
               \quad \text{where}\quad W(\hx)= \xi^u(\hx)\oplus\eta^s(\hx)\, .
$$
Notice that, since $\xi^u$ and $\eta^s$ are both $G$-invariant,
$$
\cD^k(\hx) = \cD(\hx) \, \cD(g(\hx)) \cdots \cD(g^{k-1}(\hx))
$$
for all $k\ge 1$, where we write $\cD=\cD^1$. Notice also that
$\cD^k(\hx)$ is a subsequence of the sequence $\Delta^n(\hx)$
defined in \eqref{eq.Dn}. Since $g$ is a return map to $\Gamma$,
this subsequence corresponds to values of $n$ for which
$\hf^n(\hx)\in\Gamma$. So, Proposition~\ref{p.Delta0} may be applied
to conclude that
\begin{equation}\label{eq.cDgrows}
\lim_{k\to\infty}\sum_{j=0}^{k-1} \log \cD(g^j(\hx))
= \lim_{k\to\infty} \cD^k(\hx) = \infty
\quad\text{for $\hmu$-almost every $\hx\in\Gamma$.}
\end{equation}
We use the following well-known fact (see \cite[Corollary~6.10]{Kr85})
to conclude that the growth is even linear:

\begin{lemma}\label{l.ergodic}
Let $T:X \to X$ be a measurable transformation preserving a
probability measure $\nu$ in $X$, and $\varphi:X\to\RR$ be a
$\nu$-integrable function such that $\lim_{n\to\infty}
\sum_{j=0}^{n-1} (\varphi\circ T^j) = + \infty$ at $\nu$-almost
every point. Then $\int \varphi\,d\nu > 0$.
\end{lemma}

Applying the lemma to $T=g$ and $\varphi=\log\cD$, we find that
\begin{equation}\label{eq.fim1}
\lim_{k\to\infty} \frac{1}{k} \log \cD^k(\hx)
 = \lim_{k\to\infty}\frac{1}{k} \sum_{j=0}^{k-1} \log\cD(g^j(\hx))
 = \int \log \cD \frac{d\hmu}{\hmu(\Gamma)} > 0
\end{equation}
at $\hmu$-almost every point. On the other hand, from
\eqref{eq.diferenca} and the relation between the Lyapunov spectra
of $\hFA$ and $G$,
\begin{equation}\label{eq.fim2}
\lim_{k\to\infty} \frac 1k \log \cD^k(\hx)
 = \frac{d_s}{d_u+d_s} (\lambda_u - \lambda_s) \frac{1}{\hmu(\Gamma)}
\end{equation}
These two relations imply that $\lambda_u>\lambda_s$. In this way,
we have shown that there is a definite gap between the first $\ell$
Lyapunov exponents and the remaining $d-\ell$ ones. Since this
applies for every $1\le \ell <d$, we conclude that the Lyapunov
spectrum is simple. The proof of Theorem~\ref{main} is complete.

\begin{remark}
\emph{A posteriori}, we get from \eqref{eq.diferenca},
\eqref{eq.fim1}, \eqref{eq.fim2} that $\Delta^n(x)$ goes linearly to
infinity when $n\to\infty$, that is, we do not need to restrict to
$\hf^n(\hx)\in\Gamma$.
\end{remark}

\appendix

\section{Extensions and applications}\label{s.AppA}

In this appendix we check that our methods apply to the Zorich
cocycles introduced in \cite{Zo96,Zo99}. We start with a few
simple comments on our hypotheses.

\subsection{Inducing}\label{ss.inducing}

Here we explain how cocycles over more general maps can often be
reduced to the case of the full countable shift. We begin by
treating the case of subshifts of countable type. In particular,
we recover the main results of \cite{BoV04}, in a stronger form.

Let $\cI$ be a finite or countable set and $T=\big(t(i,j)\big)_{i, j
\in \cI}$ be a transition matrix, meaning that every entry $t(i,j)$
is either $0$ or $1$. Define
$$
\hsig_T = \{(\iota_n)_{n\in\ZZ} \in \cI^\ZZ :
t(\iota_n,\iota_{n+1})=1 \text{ for all } n\in\ZZ\}
$$
and let $\hf_T:\hsig_T\to\hsig_T$ be the restriction to $\hsig_T$ of
the shift map on $\cI^\ZZ$. By definition, the cylinders $[\ \cdot\
]$ of $\hsig_T$ are its intersections with the cylinders of the full
space $\cI^\ZZ$. One-sided shift spaces $\psig_T\subset\cI^{\{n\ge
0\}}$ and $\nsig_T\subset\cI^{\{n < 0\}}$, and cylinders $[\ \cdot\
]^u\subset\psig_T$ and $[\ \cdot\ ]^s\subset\nsig_T$ are defined
analogously.

Let $\hnu_T$ be a probability measure on $\hsig_T$ invariant under
$\hf_T$ and whose support contains some cylinder
$[I]=[\iota_0;\iota_1, \ldots, \iota_{k-1}]$ of $\hsig_T$. By
Poincar\'e recurrence, the subset $X$ of points that return to
$[I]$ infinitely many times in forward and backward time has full
measure. Let $r(\hx)\ge 1$ be the first return time and
$$
\hg(\hx)= \hf^{r(\hx)}(\hx), \quad\text{for $\hx\in X$.}
$$
This first return map $\hg:X \to X$ may be seen as a shift on
$\hsig=\NN^\ZZ$. Indeed, let $\{J(\ell): \ell\in\NN\}$ be an
enumeration of the family of cylinders of the form
\begin{equation}\label{eq.cylinderJ}
[\iota_0; \iota_1, \ldots, \iota_{r-1}, \iota_r, \ldots,
\iota_{r+k-1}], \quad\text{with $\iota_{r+i}=\iota_i$ for $i=0, 1,
\ldots, k-1$}
\end{equation}
and $r\ge 1$ minimum with this property. Then
$$
\NN^\ZZ \to X, \quad (\ell_n)_{n\in\ZZ}
 \mapsto \bigcap_{n\in\ZZ} \hg^{-n}(J(\ell_n))
$$
conjugates $\hg$ to the shift map. Let $\hnu$ be the normalized
restriction of $\hnu_T$ to $X$. Then $\hnu$ is a $\hg$-invariant
probability measure and, assuming $\hnu_T$ is ergodic for $\hf_T$,
it is $\hg$-ergodic. The measure $\hnu$ is positive on cylinders,
since $[I]$ is contained in the support of $\hnu_T$. It has product
structure if $\hnu_T$ has. The latter makes sense because every
cylinder $[\iota]$ of $\hsig_T$ is homeomorphic to a product of
cylinders of $\psig_T$ and $\nsig_T$.

To each cocycle defined over $\hf$ by some $\hA_T:\hsig_T\to\GL$
we may associate a cocycle defined over $\hg$ by
$$
\hB(\hx) = \hA_T^{r(\hx)}(\hx).
$$
Notice that $\hB$ is continuous if $\hA_T$ is, since the return
time $r(\hx)$ is constant on each cylinder as in
\eqref{eq.cylinderJ}. Also, $\hB$ admits stable and unstable
holonomies if $\hA_T$ does: the holonomy maps for the two cocycles
coincide on the domain of $\hB$. Furthermore, the Lyapunov
exponents of $\hB$ are obtained by multiplying those of $\hA_T$ by
the average return time. Indeed, given any non-zero vector $v$,
$$
\lim_{n\to\infty} \frac 1n \log\|\hB^n(\hx)v\|
 = \lim_{n\to\infty} \frac 1n \log\|\hA^{S_n r(\hx)}(\hx)v\|\,,
 \qquad S_n r(\hx) = \sum_{j=0}^{n-1} r(g^j(\hx)),
$$
and, for $\hnu$-almost every $\hx$, this is equal to
$$
 \lim_{n\to\infty} \frac 1n S_n r(\hx)
 \lim_{m\to\infty} \frac 1m \log\|\hA^{m}(\hx)v\|
 = \frac{1}{\hnu_T([I])}
 \lim_{m\to\infty} \frac 1m \log\|\hA^{m}(\hx)v\|\,,
$$
since $n^{-1} S_n r(\hx)$ converges $\int r\,d\hnu = 1/\hnu([I])$.
In particular, the Lyapunov spectrum of either cocycle is simple if
and only if the other one is.

Finally, the cocycle $\hB$ is simple for $\hg$ if $\hA_T$ is simple
for $\hf_T$. More precisely, suppose $\hf_T$ admits points $\hp$ and
$\hz$ satisfying the conditions in Definition~\ref{d.simple} for the
cocycle defined by $\hA_T$ and such that $\hp$ is in the interior of
the support of $\hnu_T$. Let $q\ge 1$ be the minimum period of $\hp$
and $[I]=[\iota_0; \iota_1, \ldots, \iota_{qs-1}]$ be a cylinder
that contains $\hp$, with $s\ge 1$. Taking $s$ sufficiently large,
we may assume that $[I]$ is contained in the support of $\hnu_T$.
Replacing $\hz$ and $\hf^l(\hz)$ by appropriate backward and forward
iterates, respectively, we may also assume that they are both in
$[I]$. Then $\hp$ is also a periodic point for $\hg$ and $\hz$ is an
associated homoclinic point. Since the holonomies of the cocycles
defined by $\hA_T$ and $\hB$ coincide, it follows that the pinching
and twisting conditions in Definition~\ref{d.simple} hold also for
the cocycle defined by $\hB$.

In this way we have shown that \emph{our simplicity criterion
extends directly to cocycles over any subshift of countable type
$f_T:\hsig_T\to\hsig_T$.} There is also a non-invertible version of
this construction, where one starts with a one-sided subshift of
countable type $f_T:\Sigma_T\to\Sigma_T$ and an invariant
probability $\nu_T$ on $\Sigma_T$, and one constructs a first return
map $g(x)=f^{r(x)}(x)$ to some cylinder $[I]$ contained in the
support of $\nu_T$. Then $g$ is conjugate to the shift map on
$\NN^{\{n\ge 0\}}$ and the normalized restriction $\nu$ of the
measure $\nu_T$ to its domain is a $g$-invariant probability.
Moreover, the measure $\nu$ is ergodic for $g$ if $\nu_T$ is ergodic
for $f_T$. The natural extension $\hg$ of the return map may be
realized as the shift map on $\NN^\ZZ$. The lift $\hnu$ of the
probability $\nu$ is a $\hg$-invariant measure, and it is
$\hg$-ergodic if $\nu$ is ergodic for $g$. In Section~\ref{ss.summable}
we discuss conditions on $\nu$ under which the lift has product structure.
Given any $A_T:\Sigma_T\to\GL$, the map $B(x)=A_T^{r(x)}(x)$ defines a cocycle
over $g$. Moreover, $B$ lifts canonically to a cocycle $\hB$ over
$\hg$, constant on local stable sets, and having the same Lyapunov
exponents. Thus, the Lyapunov spectrum of $A_T$ is simple if and
only if the Lyapunov spectrum of $\hB$ (or $B$) is.

More generally, let $f:M\to M$ be a transformation preserving a
probability $\nu_f$ and assume there exists a return map $g$ to some
domain $D\subset\supp\nu_f$ which is a \emph{Markov map}. By this we
mean that there exists a finite or countable partition $\{J(\ell):
\ell\in\NN\}$ of $D$ such that (i) $g$ maps each $J(\ell)$
bijectively to the whole domain $D$ and (ii) for any sequence
$(\ell_n)_n$ in $\NN^{\{n\ge 0\}}$ the intersection of
$g^{-n}(J(\ell_n))$ over all $n\ge 0$ consists of exactly one point.
Then $g$ may be seen as the shift map on $\NN^{\{n\ge 0\}}$. The
normalized restriction $\nu$ of $\nu_f$ to the domain of $g$ is a
$g$-invariant probability, and it is $g$-ergodic if $\nu_f$ is
ergodic for $f$. As before, to any cocycle over $f$ we may
associated a cocycle over $g$, or its natural extension, such that
the Lyapunov spectrum of either is simple if and only if the other
one is. This type of construction will be used in
Section~\ref{ss.zorich}.

\subsection{Bounded oscillation}\label{ss.summable}

Let $f:\Sigma\to\Sigma$ be the shift map on $\Sigma=\NN^{\{n\ge
0\}}$. The \emph{lift} of an $f$-invariant probability measure $\mu$
is the unique $\hf$-invariant measure $\hmu$ on $\hsig=\NN^\ZZ$ such
that $P_*\hmu=\mu$. The \emph{$k$-oscillation} of a function
$\psi:\Sigma\to\RR$ is defined by
$$
\osc_k(\psi) = \sup_I \sup\{\psi(x)-\psi(y): x, y \in [I]\}
$$
where the first supremum is over all sequences $I=(\iota_0, \ldots,
\iota_k)$ in $\NN^k$. We say $\psi$ has \emph{bounded oscillation}
if $\sum_{k=1}^\infty \osc_k(\psi) < \infty$. This implies
$\osc_k(\psi)\to 0$ and so $\psi$ is continuous, in a uniform sense.
We are going to prove

\begin{proposition}\label{p.lift}
If the Jacobian of $\nu$ for $f$ has bounded oscillation then the
lift $\hmu$ has product structure.
\end{proposition}

\begin{lemma}\label{l.newconditional1}
Let $x$ and $y$ be in $\Sigma=\psig$. For each point $\hx\in
W^s_{loc}(x)$, define $\hy\in \Wuloc(\hx) \cap\Wsloc(y)$. Then the
limit
$$
J_{x,y}(\hx)= \lim_{n\to\infty} \frac{\Jmufn(\xn)}{\Jmufn(\yn)}\,,
 \quad\text{where $\xn=P(\hf^{-n}(\hx))$ and $\yn=P(\hf^{-n}(\hy))$,}
 $$
exists, uniformly on $x$, $y$, and $\hx$. Moreover, the function
$(x,y,\hx)\mapsto J_{x,y}(\hx)$ is continuous and uniformly
bounded from zero and infinity.
\end{lemma}

\begin{proof}
The arguments are quite standard. Begin by noting that
\begin{equation}\label{e.newprevious}
\log \frac{\Jmufn(\xn)}{\Jmufn(\yn)} = \sum_{j=1}^n
\log\Jmuf(\xj)-\log\Jmuf(\yj).
\end{equation}
Notice that $\xj$ and $\yj$ are in the same cylinder $[\iota_{-j},
\ldots, \iota_{-1}]^u$, for each $j\ge 1$. Hence, the $j$th term in
the sum is bounded in norm by the $j$-oscillation of $\log\Jmuf$. It
follows that the series in \eqref{e.newprevious} converges
absolutely and uniformly, and the sum is bounded by $\sum_{j}
\osc_j(\log\Jmuf)$. This implies all the claims in the lemma.
\end{proof}

\begin{lemma}\label{l.newconditional2}
Let $\{\hmu_x: x\in\Sigma\}$ be any disintegration of the lift
$\hmu$ of $\mu$. For a full $\mu$-measure subset of points
$x\in\Sigma$, we have
$$
\hmu_x(\xi_n) = \frac{1}{\Jmuf^n(\xn)}
$$
for every cylinder $\xi_n=[\iota_{-n}, \ldots, \iota_{-1}]^s$, $n\ge
1$ and every point $\hx\in\xi_n\times\{x\}$.
\end{lemma}

\begin{proof}
Let $F$ be any measurable subset of $\Sigma$. Then $\hf^{-n}(\xi_n
\times F)=P^{-1}(F_n)$, where $F_n$ is the subset of $[\iota_{-n},
\ldots, \iota_{-1}]^u$ that is sent bijectively to $F$ by the map
$f^{n}$. Consequently,
\begin{equation}\label{e.newjacobian1}
\frac{\hat\mu(\xi_n \times F)}{\mu(F)}
 = \frac{\hat\mu(P^{-1}(F_n))}{\mu(F)}
 = \frac{\mu(F_n)}{\int_{F_n}\Jmufn\,d\mu}.
\end{equation}
On the other hand, for $\mu$-almost any point $x\in\Sigma$ and any
cylinder $\xi_n\subset\nsig$,
$$
\hmu_x(\xi_n)
 = \lim_{F\to x}\frac{\hat\mu(\xi_n \times F)}{\mu(F)}
$$
where the limit is over a basis of neighborhoods $F$ of $x$. As $F
\to x$, the sets $F_n$ converge to the unique point in $[\iota_{-n},
\ldots, \iota_{-1}]^u$ that is mapped to $x$ by $f^n$. This point is
precisely $\xn=P(\hf^{-n}(\hx))$, for any choice of
$\hx\in\xi_n\times\{x\}$. In view of \eqref{e.newjacobian1},
this gives that
$$
\hmu_x(\xi_n) = \frac{1}{\Jmufn(\xn)}
$$
for every cylinder $\xi_n$ and any $x$ in some full $\mu$-measure
subset.
\end{proof}

\begin{lemma}\label{l.newconditional3}
There exists a disintegration $\{\hmu_x:x\in\Sigma\}$ of the lift
$\hmu$ such that $\hat\mu_y = J_{x,y} \hat\mu_x$ for every $x$ and
$y$ in $\Sigma$.
\end{lemma}

\begin{proof} Let $\{\bmu_x: x\in\Sigma\}$ be an arbitrary disintegration.
By the previous lemma, there exists a full measure subset $S$ of
$\Sigma$ such that
\begin{equation}\label{eq.newjacobian2}
\frac{\bmu_y(\xi_n)}{\bmu_x(\xi_n)}
 = \frac{\Jmufn(\xn)}{\Jmufn(\yn)}
 \quad\text{for any $\xi_n=[\iota_{-n}, \ldots, \iota_{-1}]^s$
  and any $x, y \in S$,}
\end{equation}
where $\xn=P(\hf^{-n}(\hx))$ and $\yn=P(\hf^{-n}(\hy))$ for any
$\hx\in\xi_n\times\{x\}$ and $\hy\in\xi_n\times\{y\}$. Define
$J_{n,x,y}$ to be the function on $\Wsloc(x)$ which is constant
equal to the right hand side of \eqref{eq.newjacobian2} on each
$\xi_n\times\{x\}$. Given any cylinder $\eta\subset\nsig$ and any
large $n\ge 1$, we may write
$$
\bmu_y(\eta)
 = \sum_{\xi_n\subset\eta} \bmu_y(\xi_n)
 = \sum_{\xi_n\subset\eta} J_{n,x,y}(\hx) \bmu_x(\xi_n)
 = \int_{\eta} J_{n,x,y}(\hx) \, d\hmu_x(\hx),
$$
where the sum is over all the cylinders $\xi_n$ that form $\eta$.
Passing to the limit as $n\to\infty$, we obtain from
Lemma~\ref{l.newconditional1} that
$$
\bmu_y(\eta) = \int_{\eta} J_{x,y} \, d\bmu_x \quad\text{for any
cylinder $\eta\subset\nsig$.}
$$
This shows that $\bmu_y = J_{x,y} \bmu_x$ for every $x$ and $y$ in
the full measure set $S$. Fix any $\bx\in S$ and define $\hmu_y =
J_{\bx,y} \bmu_\bx$ for every $y\in\Sigma$. Then $\hmu_y=\bmu_y$ for
every $y\in S$, and so $\{\hmu_x\}$ is a disintegration of $\hmu$.
Moreover,
$$
\hmu_y
 = J_{\bx,y} \bmu_\bx
 = J_{x,y} J_{\bx,x} \bmu_\bx
 = J_{x,y}\hmu_x
$$
for any $x, y \in \Sigma$, as claimed in the lemma.
\end{proof}

\begin{proof}[Proof of Proposition~\ref{p.lift}]
Fix an arbitrary point $\xplus$ in $\Sigma$ and then define
$$
r(\nx,\px)=J_{\xplus,\px}(\nx,\px) \quad\text{for every
$\hx=(\nx,\px)\in\hsig$.}
$$
By the previous lemma, $\hmu_{\px} = r(\nx,\px) \hmu_\xplus$ for
every $\px\in\Sigma$. The lift $\hmu$ projects to $\pmu=\mu$ on
$\Sigma$, by definition. The projection $\nmu$ to $\nsig$ is given
by
$$
\nmu = \hmu_\xplus \int_{\Sigma} r(\nx,\px) \,d\mu(\px).
$$
It follows that $\hmu = \rho(\nx,\px) \nmu \times \pmu$, with
$$
\rho(\nx,\px) = \frac{r(\nx,\px)}{\int_{\Sigma} r(\nx,\px)
\,d\mu(\px)}\,.
$$
Since the function $r(\nx,\px)$ is continuous and uniformly bounded
from zero and infinity, so is the density $\rho$. This implies that
$\hat\mu$ has product structure.
\end{proof}

\subsection{Fiber bunched cocycles}\label{ss.dominated}

As pointed out in Section~\ref{sh.cocycle}, existence of stable and
unstable holonomies is automatic when the cocycle is locally
constant. Another, more robust, construction of cocycles with stable
and unstable holonomies was given in \cite{BGV03}. Let us recall it
briefly here.

\begin{definition}\label{d.dominated}
We say that $\hA:\hsig\to\GL$ is \emph{$s$-fiber bunched} (or
\emph{$s$-dominated}) for $\hf:\hsig\to\hsig$ if there exist
constants $N\ge 1$, $C>0$, $\nu \in (0,1]$, $\tau \in (0,1)$, and
$\theta\in (0,1)$ , and a distance $d$ on $\hsig$, such that
\begin{itemize}
 \item[(a)] $d(\hf^N(\hx), \hf^N(\hy))\leq \theta d(\hx,\hy)$
 if $\hx, \hy$ are in the same local stable set
 \item[(b)] $\|\hA^N(\hx)^{\pm 1}\|\le C$ \ and \ $\|\hA^N(\hx)-\hA^N(\hy)\|\le C d(\hx,\hy)^\nu$
 \item[(c)] $\|\hA^N(\hx)\|\|\hA^N(\hx)^{-1}\| \theta^{\,\nu} < \tau$
\end{itemize}
for every $\hx, \hy \in \hsig$. We say that $\hA$ is \emph{$u$-fiber
bunched} (or \emph{$u$-dominated}) for $\hf$ if $\hA^{-1}$ is
$s$-fiber bunched for $\hf^{-1}$.
\end{definition}

\begin{proposition}\label{p.conjug}
If $\hA$ is $s$-fiber bunched (respectively, $u$-fiber bunched) then
it admits stable holonomies (respectively, unstable holonomies).
\end{proposition}

\begin{proof}
Replacing $\hf$ by $\hf^N$ in Definition~\ref{d.dominated}, we may
assume $N=1$. Denote $H_n(\hx,\hy)= \hA^n(\hy)^{-1} \hA^n(\hx)$ for
each $n\ge 1$ and $\hx$ and $\hy$ in the same local stable set. Then
$$
H_{n+1}(\hx,\hy) - H_{n}(\hx,\hy)
 = \hA^{n}(\hy)^{-1} \hA(\hf^n(\hy))^{-1} \big[\hA(\hf^n(\hx)) - \hA(\hf^n(\hy))\big] \hA^n(\hx)
$$
By condition (a), we have $d(\hf^n(\hx), \hf^n(\hy)) \le \theta^n
d(\hx,\hy)$. Using condition (b), it follows that
$$
\|H_{n+1}(\hx,\hy) - H_{n}(\hx,\hy)\|
 \le C^2 d(\hx,\hy)^\nu \prod_{j=0}^{n-1} \Big(\|\hA(\hf^j(\hy))^{-1}\|\,\|\hA(\hf^j(\hx))\|
\, \theta^\nu\Big).
$$
Fix $\hat\tau \in (\tau, 1)$. By conditions (a) and (b),
$\hA(\hf^j(\hy))$ is close to $\hA(\hf^j(\hx))$ when $j$ is large,
uniformly on $\hx$ and $\hy$. Combining this with condition (c), we
get that there exists $k\ge 1$, independent of $\hx$ and $\hy$, such
that
$$
\|\hA(\hf^j(\hy))^{-1}\| \, \|\hA(\hf^j(\hx))\| \theta^\nu <
\hat\tau
$$
for all $j\ge k$. Thus, the previous inequality implies that
$$
\big\|H_{n+1}(\hx,\hy) - H_n(\hx,\hy)\big\|
 \le C^2 d(\hx,\hy)^\nu \, C^{2k} \theta^{k \nu} \, \hat\tau^{n-k}
 \le \hat{C} \hat\tau^n d(\hx,\hy)^\nu,
$$
for some appropriate constant $\hat{C}>0$. This implies that $H_n$
is a Cauchy sequence, uniformly on $(x,y)$. Hence, it is uniformly
convergent, as claimed. This proves that $\hA$ admits stable
holonomies if $\hA$ is $s$-fiber bunched. The dual statement is
proved in just the same way.
\end{proof}

We say that $\hA:\hsig\to\GL$ is \emph{fiber bunched} if it is
simultaneously $s$-fiber bunched and $u-$fiber bunched. From
Proposition~\ref{p.conjug} we immediately get that if $\hA$ is fiber
bunched then it admits stable and unstable holonomies.

\begin{remark}
In some cases it is possible to reduce non-fiber bunched cocycles to
the fiber bunched case. For instance, let $F=(f,A)$ be a linear
cocycle $F=(f,A)$ over a shift map, say, which is not fiber bunched
but whose Lyapunov spectrum is narrow, meaning that the difference
between all Lyapunov exponents is sufficiently small. Then we may
use inducing to construct from $F$ a fiber bunched cocycle.
\end{remark}

%

\subsection{Zorich cocycles}\label{ss.zorich}

Finally, we are going to explain how the methods in this paper can
be applied to Zorich cocycles~\cite{Zo96,Zo99}. We begin by
recalling the definition of these cocycles. Motivations and proofs
for the results we quote can be found in Kontsevich,
Zorich~\cite{KZ97}, Marmi, Moussa, Yoccoz~\cite{MMY05},
Rauzy~\cite{Ra79}, Veech~\cite{Ve82,Ve86}, Zorich~\cite{Zo96,Zo99},
and references therein. See also \cite{AV2}, where we show that
Zorich cocycles are simple, thus proving the Zorich-Kontsevich
conjecture that the corresponding Lyapunov spectra are simple.

\subsubsection{The Rauzy algorithm}

Fix some integer $d\ge 2$. Let $\Pi=\Pi_d$ be the set of all
irreducible pairs $\pi=(\pi_0,\pi_1)$ of permutations
$\pi_\vep=(\alpha^\vep_1, \alpha^\vep_2, \ldots, \alpha^\vep_d)$ of
the alphabet $\{1, \ldots, d\}$. By \emph{irreducible} we mean that
$\pi_1\circ\pi_0^{-1}$ preserves no subset $\{1, \ldots, k\}$ with
$k<d$. We shall denote the rightmost symbol $\alpha^\vep_d$ simply
as $\alpha(\vep)$ for $\vep\in\{0, 1\}$. Let $\Delta=\Delta_d$ be
the standard open simplex of dimension $d-1$, that is, the set of
all vectors $\lambda=\big(\lambda_1, \lambda_2, \ldots,
\lambda_d\big)$ such that $\lambda_j>0$ for all $j$ and
$\sum_{j=1}^d\lambda_j=1$. We call $g:\Delta\to\Delta$ a projective
map if there exists a linear isomorphism $G:\RR^d\to\RR^d$ with
non-negative coefficients such that
\begin{equation}\label{eq.projectivization}
g(\lambda)
 = \frac{G(\lambda)}{\sum_{i=1}^d G(\lambda)_i}
  = \frac{G(\lambda)}{\sum_{i, j=1}^d G_{i,j}\lambda_j}\,.
\end{equation}
If the coefficients of $G$ are strictly positive then the image of
$g$ is relatively compact in $\Delta$. In this case $g$ is a
contraction for the \emph{projective metric} defined in $\Delta$ by
$$
d(\lambda,\lambda')
 = \log \max\{\frac{\lambda_i\lambda_j'}{\lambda_j\lambda'_i}: i, j = 1, \ldots, d\}.
$$
The contraction rate depends only on a lower bound for the
coefficients of $G$ or, equivalently, for the Euclidean distance
from $g(\Delta)$ to the boundary of $\Delta$.

Let $\cR:(\pi,\lambda)\mapsto (\pi',\lambda')$ be defined on an open
dense subset of $\Pi\times\Delta$, as follows. For each $\pi\in\Pi$
and $\vep\in\{0,1\}$, let
$$
\Delta^\vep(\pi)=\{\lambda\in\Delta:
\lambda_{\alpha(\vep)}>\lambda_{\alpha(1-\vep)}\}.
$$
We say that $(\pi,\lambda)$ has \emph{type} $\vep$ if
$\lambda\in\Delta^\vep(\pi)$. Then, by definition,
$\pi_\vep'=\pi_\vep$ and
$$
\pi_{1-\vep}'=\big(\alpha_1^{1-\vep}, \ldots,
\alpha_{k-1}^{1-\vep}, \alpha(1-\vep), \alpha_{k}^{1-\vep},
\ldots, \alpha_{d-1}^{1-\vep}\big)
$$
where $k\in\{1, \ldots, d-1\}$ is defined by
$\alpha_k^{1-\vep}=\alpha(\vep)$. In other words, $\pi_{1-\vep}'$
is obtained from $\pi_{1-\vep}$ by looking for the position $k$
the last symbol of $\pi_\vep$ occupies in $\pi_{1-\vep}$, leaving
all symbols to the left of $k$ unchanged, and rotating the symbols
to the right of $k$ one position to the right. Moreover,
$$
\lambda'_{j}=\frac{1}{a}\,\lambda_{j}
 \text{\quad for\quad} j\neq\alpha(\vep),
 \qquad
 \lambda'_{j}=\frac{1}{a}\big(\lambda_{\alpha(\vep)}-\lambda_{\alpha(1-\vep)}\big)
 \text{\quad for\quad} j=\alpha(\vep)
$$
where the normalizing factor $a=1-\lambda_{\alpha(1-\vep)}$. Notice
that $\lambda\mapsto\lambda'$ sends each $\Delta^\vep$ bijectively
onto $\Delta$. Moreover, this map is just the projectivization of
the linear isomorphism $R_{\pi,\lambda}:\RR^d\to\RR^d$
$$
\big(\lambda_1, \ldots, \lambda_{l-1},
 \lambda_{\alpha(\vep)},
 \lambda_{l+1} \ldots, \lambda_d\big)
 \mapsto \big(\lambda_1, \ldots, \lambda_{l-1},
 \lambda_{\alpha(\vep)}-\lambda_{\alpha(1-\vep)},
 \lambda_{l+1}, \ldots,
 \lambda_{d}\big),
$$
in the sense that $\lambda'=(1/a) R_\pi(\lambda)$ with
$a=\sum_{i=1}^d (R_\pi\lambda)_i$. It is interesting to write this
also as $\lambda = a R_{\pi,\lambda}^{-1}(\lambda')$, because the
inverse operator
$$
\big(\lambda_1, \ldots, \lambda_{l-1},
 \lambda_{\alpha(\vep)},
 \lambda_{l+1} \ldots, \lambda_d\big)
 \mapsto \big(\lambda_1, \ldots, \lambda_{l-1},
 \lambda_{\alpha(\vep)}+\lambda_{\alpha(1-\vep)},
 \lambda_{l+1}, \ldots,
 \lambda_{d}\big).
$$
has non-negative integer coefficients.

Let us call a Rauzy component of $\Pi\times\Delta$ any smallest set
of the form $\Pi_0\times\Delta$ which is invariant under $\cR$. From
now on we always consider the restriction of the algorithm to some
Rauzy component. The map $\cR$ admits an absolutely continuous
invariant measure $\nu$, that is, an invariant measure such that the
restriction to each $\{\pi\}\times\Delta$ is absolutely continuous
with respect to Lebesgue measure on the standard simplex. However,
$\nu$ is usually infinite. This can be overcome by considering the
following accelerated algorithm.

\subsubsection{The Zorich algorithm}

Define $\cZ(\pi,\lambda)=(\cR^{n})(\pi,\lambda)$, where the
acceleration time $n=n(\pi,\lambda)\ge 1$ is the largest number of
consecutive iterates by the Rauzy algorithm during which the type
remains unchanged. In precise terms, $n=n(\pi,\lambda)$ is
characterized by (assume $(\pi^{(i)},\
\lambda^{(i)})=\cR^i(\pi,\lambda)$ is defined for all $0\le i \le
n$)
$$
 (\pi^{(i)},\ \lambda^{(i)}) \text{ has type $\vep$ for $0\le i < n$}
 \quand (\pi^{(n)},\ \lambda^{(n)}) \text{ has type $1-\vep$.}
$$
Since each
$\cR:\{\pi^{(i)}\}\times\Delta^\vep(\pi^{(i)})\to\{\pi^{(i+1)}\}\times\Delta$
is a projective bijection, the map $\cR^n$ sends some sub-simplex
$\{\pi\}\times D(\pi,\lambda)\subset\{\pi\}\times\Delta^\vep(\pi)$
containing $(\pi,\lambda)$ bijectively onto
$\{\pi^{(n)}\}\times\Delta^{1-\vep}(\pi^{(n)})$. Moreover, its
inverse is the restriction of a projective map
$\{\pi^{(n)}\}\times\Delta\to\{\pi\}\times\Delta$. By definition,
$\cZ=\cR^n$ restricted to $D(\pi,\lambda)$. Let $\cD$ be the
(countable) family of all these sub-simplices $D(\pi,\lambda)$. The
union of its elements has full measure on $\Pi\times\Delta$.

The transformation $\cZ$ admits an absolutely continuous invariant
probability measure $\mu$ on each Rauzy component, and this measure
is ergodic. Moreover, the density of $\mu$ is a rational function of
the form
\begin{equation}\label{eq.rational}
\frac{d\mu}{dm}(\lambda)
 = \sum_\alpha \frac{1}{\cP_\alpha(\lambda)}
 \quad\text{on each domain } \{\pi\}\times\Delta
\end{equation}
where the sum is over some finite set of polynomials $\cP_\alpha$
with non-negative coefficients and degree $d$. In particular, the
density is smooth and bounded from zero on every
$\{\pi\}\times\Delta$. In general, the density is not bounded from
infinity, because the $\cP_\alpha$ may have zeros on the boundary of
$\Delta$.

\subsubsection{Linear cocycles}
The \emph{Rauzy cocycle} over $\cR$ is defined by
$$
F_R:\Pi\times\Delta\times\RR^d\to\Pi\times\Delta\times\RR^d, \quad
\big(\pi,\lambda,v\big) \mapsto \big(\cR(\pi,\lambda),
R^{-1*}_{\pi,\lambda}(v)\big).
$$
Notice that this cocycle is constant on each $\Delta^\vep(\pi)$,
because $R_{\pi,\lambda}$ depends only on $\pi$ and the type $\vep$
of $\lambda$. The \emph{Zorich cocycle} over $\cZ$ is defined by
$$
F_Z:\Pi\times\Delta\times\RR^d\to\Pi\times\Delta\times\RR^d,
 \quad F_Z\big(\pi,\lambda,v\big)=F_R^{n(\pi,\lambda)}\big(\pi,\lambda,v\big)
$$
Notice that $F_Z\big(\pi,\lambda,v\big)=\big(\cZ(\pi,\lambda),
Z_{\pi,\lambda}(v))$ where $Z_{\pi,\lambda}$ is constant on each
element of $\cD$ and its inverse has non-negative integer
coefficients. The Zorich cocycle is integrable with respect that the
$\cZ$-invariant measure $\mu$, meaning that
$\log\|Z_{\pi,\lambda}^{\pm 1}\|$ are integrable functions. Thus,
its Lyapunov exponents are well-defined at $\mu$-almost every point.
By ergodicity, the exponents are constant $\mu$-almost everywhere.

Consider the linear map $\Omega_\pi:\RR^d\to\RR^d$ defined by
$$
\Omega_\pi(\lambda)_i = \sum_{j\ :\ \pi_1(j) < \pi_1(i)} \lambda_j -
                        \sum_{j\ :\ \pi_0(j) < \pi_0(i)} \lambda_j \, .
$$
This map $\Omega_\pi$ is anti-symmetric (not necessarily an
isomorphism), and so
$$
\omega_\pi\big(\Omega_\pi(u), \Omega_\pi(v)\big)
 = u \cdot \Omega_\pi(v)
$$
defines a symplectic form on the range $H_\pi=\Omega_\pi(\RR^d)$. In
particular, the dimension of $H_\pi$ is even. The map $\Omega_\pi$
also satisfies
\begin{equation}\label{eq.symplectic}
\Omega_{\pi'} \cdot R_{\pi,\lambda} = R_{\pi, \lambda}^{-1*} \cdot
\Omega_\pi.
\end{equation}
This implies that the Rauzy cocycle leaves invariant the subbundle
$$
\cH_\Pi = \{(\pi,\lambda,v)\in\Pi\times\Delta\times\RR^d : v \in
H_\pi\}
$$
and even preserves the symplectic form $\omega_\pi$ on it. Then the
same is true for the Zorich cocycle.

It follows that the Lyapunov spectrum of the Zorich cocycle
restricted to the subbundle $\cH_\Pi$ has the form
\begin{equation}\label{eq.spectrum}
 \lambda_1 \ge \cdots \ge \lambda_g \ge 0
 \ge -\lambda_g \ge \cdots \ge -\lambda_1
 \quad\text{(where $2g=\dim H_\pi$).}
\end{equation}
The other Lyapunov exponents of $F_Z$, corresponding to directions
transverse to $\cH_\pi$, vanish identically and are not of interest
here. The \emph{Zorich-Kontsevich conjecture} states that all the
inequalities in \eqref{eq.spectrum} are strict or, in other words,
\emph{the Lyapunov spectrum of the restricted Zorich cocycle is
simple.} We are going to argue that, modulo the simple observations
in Sections~\ref{ss.inducing} and~\ref{ss.summable}, all the
hypotheses of Theorem~\ref{main} are satisfied in the context of
Zorich cocycles, and so our methods can be used to prove this
conjecture.

\subsubsection{Inducing on a compact simplex}

Let $\cD$ be the family of sub-simplices introduced in the
definition of the Zorich algorithm: $\cZ$ maps each element of $\cD$
bijectively to some $\{\pi'\}\times\Delta^{1-\vep}$, and the inverse
is the restriction of a projective map
$\{\pi'\}\times\Delta\to\{\pi\}\times\Delta$. Pulling $\cD$ back
under $\cZ$ we obtain, for each $n\ge 1$, a countable family $\cD^n$
of sub-simplices each of which is mapped bijectively to some
$\{\pi^{(n)}\}\times\Delta^{1-\vep}$ by the iterate $\cZ^n$, the
inverse being the restriction of a projective map
$\{\pi^{(n)}\}\times\Delta\to\{\pi\}\times\Delta$. For $\mu$-almost
every $(\pi,\lambda)$, there exists some $n\ge 1$ for which this
projective map has strictly positive coefficients, and so the image
$\{\pi\}\times \Gamma$ is relatively compact in
$\{\pi\}\times\Delta$. Let us fix such $n$, $\pi$, $\lambda$ once
and for all, and denote by $\{\pi\}\times D_*$ the corresponding
element of $\cD^n$. In particular, $D_*\subset \Gamma$ is relatively
compact in $\Delta$. It follows that $D_*$ has finite diameter for
the projective metric of $\Delta$, and also that the density
$d\mu/dm$ is smooth and bounded from zero and infinity on $D_*$. For
notational simplicity, we identify $\{\pi\}\times\Delta \approx
\Delta$ and $\{\pi\}\times D_* \approx D_*$ in what follows.

By Poincar\'e recurrence, there exists a first return map $\cG$ of
the map $\cZ^n$ to the domain $D_*$. More precisely, using the
Markov structure of $\cZ^n$, there exists a countable family
$\{D_\iota : \iota\in\NN\} \subset \cup_{k\ge 1} \cD_{kn}$ of
sub-simplices of $D_*$ such that their union has full measure in
$D_*$, each $D_\iota$ is mapped bijectively to the whole $D_*$ by
$\cG$, and the inverse of each $\cG: D_\iota \to D_*$ is the
restriction of a projective map $\Delta\to\Delta$. By construction,
the images of these inverse branches are all contained in $\Gamma$,
and so they all contract the projective metric, with uniform
contraction rates. Let $D\subset D_*$ be the (full measure) subset
of points that return infinitely many times to $D_*$. In particular,
the map
$$
\Phi: \NN^{\{n\ge 0\}} \to D \quad
 (\iota_n)_n \mapsto \cap_{n\ge 0} \cG^{-n}(D_{\iota_n})
$$
is well defined (the intersection consists of exactly one point),
and it conjugates $\cG:D \to D$ to the shift map on $\NN^{\{n\ge
0\}}$. Then the natural extension of $\cG$ is realized by the shift
map on $\NN^\ZZ$.

On the one hand, as observed before, the invariant density $d\mu/dm$
is smooth and bounded from zero and infinity on $D$. It follows that
its logarithm is bounded and Lipschitz continuous, for either
Euclidean or projective metric, with uniform constants. On the other
hand, the inverse branches of $\cG$ are all projective maps with
range contained in the same relatively compact domain $\Gamma$. This
implies that the logarithms of their derivatives are also bounded
and Lipschitz continuous, for either metric, with uniform constants.
Putting these two facts together we get that the logarithm of the
Jacobian of $\cG$ with respect to the measure $\mu$ is uniformly
bounded and Lipschitz continuous on each $D_\iota$. Combining this
with the previous observation that inverse branches of $\cG$
contract the projective metric uniformly, we easily obtain that
$\log J\cG$ has bounded oscillation in the sense of
Section~\ref{ss.summable}. Consequently, the lift of $\mu \mid D$ to
the natural extension of $\cG$ has product structure.

Recall that the Zorich cocycle $F_Z$ is constant on each element of
$\cD$. It is clear from the construction that points in each
$D_\iota$ visit exactly the same elements of $\cD$ all the way up to
their return to $D_*$. Thus, the linear cocycle $F_G$ induced by
$F_Z$ over the return map $\cG$ is also locally constant, meaning
that it is constant on each $D_\iota$. In particular, the cocycle
$F_G$ is continuous for the shift topology, and it admits stable and
unstable holonomies.

\subsubsection{Pinching and twisting conditions}\label{sss.pinching}

The only missing ingredient to establish the Zorich-Kontsevich
conjecture is to prove that the Zorich cocycles are simple, in the
sense of Definition~\ref{d.simple}. This is done in \cite{AV2}.
In fact, the pinching and twisting conditions appear in a slightly
different guise in that paper, in terms of the monoid generated
by the cocycle.

In this context, a \emph{monoid} is just a subset of $\GL$ closed
under multiplication and containing the identity. The
\emph{associated monoid} $\cB=\cB(F)$ is the smallest monoid that
contains the image of $F$. We call $\cB$ is \emph{simple} if it is
both pinching and twisting, where $\cB$ is
\begin{itemize}
\item \emph{pinching} if it contains elements with arbitrarily large
eccentricity $\Ec(B)$
\item \emph{twisting} if for any $F\in\Grass$ and any finite family
$G_1, \ldots, G_N $ of elements of $\Grass$ there exists $B\in\cB$
such that $B(F) \cap G_i =\{0\}$ for all $j=1, \ldots, N$.
\end{itemize}
The \emph{eccentricity} of a linear map $B\in\GL$ is defined by
$$
\Ec(B) = \min_{1\le \ell < d}\frac{\sigma_\ell}{\sigma_{\ell+1}}
$$
where $\sigma_1^2 \ge \cdots \ge \sigma_d^2$ are the eigenvalues of
the self-adjoint operator $B^*B$, in non-increasing order.
Geometrically, the positive square roots $\sigma_1 \ge \cdots \ge
\sigma_d$ correspond to the lengths of the semi-axes of the
ellipsoid $\{B(v): \|v\|=1\}$. It is evident from the definition
that any monoid that contains a pinching submonoid is also pinching,
and analogously for twisting.

It is not difficult to see that the two formulations of the
definition of simplicity are equivalent, for locally constant real
cocycles. Indeed, Lemma ~A.5 in~\cite{AV2} states that if the
associated monoid is simple then there exists some periodic point
and some homoclinic point as in Definition~\ref{d.simple}.
Conversely, the conditions in Definition~\ref{d.simple} imply that
the associated monoid is simple. Indeed, the first condition implies
that $\cB$ contains some element $B_1$ whose eigenvalues all have
distinct norms. Then the powers $B_1^n$ have arbitrarily large
eccentricity as $n\to\infty$, and so $\cB$ is pinching. Moreover,
the second condition implies that the monoid contains some element
$B_2$ satisfying $B_2(V)\cap W =\{0\}$ for any pair of subspaces $V$
and $W$ which are sums of eigenspaces of $B_1$ and have
complementary dimensions. Given any $F$, $G_1$, \ldots, $G_n$ as in
the definition, we have that $B_1^n(F)$ is close to some sum $V$ of
$\ell$ eigenspaces of $B_1$, and every $B_1^{-n}(G_i)$ is close to
some sum $W_i$ of $d-\ell$ eigenspaces of $B_1$, as long as $n$ is
large enough. It follows that $B_2(B_1^n(F)) \cap
B_1^{-n}(G_i)=\{0\}$, that is, $B_1^n B_2 B_1^n(F) \cap G_i=\{0\}$.
This proves $\cB$ is twisting.


\section{Intersections of hyperplane sections}\label{s.AppB}

Here we give an alternative proof of Proposition~\ref{p.noatoms}
under the assumption that the eigenvalues of the cocycle at the
fixed point $p$ are real. Observe that this is automatic for real
cocycles, since we also assume that the absolute values of the
eigenvalues are all distinct. Instead of Proposition~\ref{p.artur}
we use the following result, which has a stronger conclusion.

\begin{proposition}\label{p.disjointness}
There exists $N=N(\ell,d)$ such that
$$
B^{-m_1}(V) \cap \cdots \cap B^{-m_N}(V) = \emptyset
$$
for any $B:\CC^d\to\CC^d$ whose eigenvalues all have
distinct absolute values, any hyperplane section $V$ of $\Grass$
containing no eigenspace of $B$, and any $0\le m_1 < \cdots <
m_N$.
\end{proposition}

To deduce Proposition~\ref{p.noatoms} from this result, one can
use the same arguments as in Section~\ref{s.simple}, just replacing
the paragraph that contains \eqref{eq.disjointW} by the following one.

Applying Proposition~\ref{p.disjointness} with $B=A^q(p)$ and $V=W_\hp$
we conclude that the $W^n_\hp$ are $N$-wise disjoint:
$$
W_\hp^{m_1} \cap \cdots \cap W_\hp^{m_N} = \emptyset
 \quad\text{for all $1\le m_1 < \cdots < m_N$.}
$$
Fix $C\ge 1$ such $C\gamma_0>1$. By continuity, we have
$W_\eta^{m_1} \cap \cdots \cap W_\eta^{m_N} = \emptyset$ for all
$1\le m_1 < \cdots < m_N \le C N$ and every $\eta$ in a small
neighborhood of $\hp$ inside the local stable set. Then, for
$\hmu_p$-almost every $\eta$ in that neighborhood,
$$
\hm_\eta(\bigcup_{j=1}^{C N} W_\eta^j)
 \ge \frac 1 N \sum_{j=1}^{C N} \hm_\eta(W_\eta^j)
 = C \gamma_0 > 1.
$$
This is a contradiction, since $\hm_\eta$ is a probability. This
contradiction reduces the proof of Proposition~\ref{p.noatoms}
to proving Proposition~\ref{p.disjointness}.

\medskip


In the proof of Proposition~\ref{p.disjointness} we use the
following classical fact about Vandermonde type determinants (see
Mitchell~\cite{Mi81}). Given $N\ge 1$, $x=(x_1, \ldots,
x_N)\in\RR^N$, and $m=(m_1, \ldots, m_N)\in(\NN\cup\{0\})^N$, define
$$
\Delta_m(x) = \left|\begin{array}{ccc}
   x_1^{m_1} & \cdots & x_N^{m_1} \\
   \cdots & \cdots & \cdots     \\
   x_1^{m_N} & \cdots & x_N^{m_N}
   \end{array}\right|.
$$

\begin{proposition}\label{p.vandermonde}
Suppose $0\le m_1 < m_2 < \cdots < m_N$. Then
$$
\Delta_m(x) = \cP_m(x) \prod_{1\le i < j \le d} (x_j - x_i)
$$
where $\cP_m$ is a positive polynomial, in the sense that all its
monomials have positive coefficients. In particular, $\Delta_m(x)$
is different from zero whenever the $x_j$ are all positive and
distinct.
\end{proposition}

Notice that the contents of the proposition remains the same if one
replaces $B$ by its square. Indeed, it is trivial that the statement
for $B$ implies the one for $B^2$, and the converse is also easy to
check: if the $B^2$-iterates of any hyperplane section $V$ as in the
statement are $N$-wise disjoint then, using this fact both for $V$
and for $B(V)$, the $B$-iterates of any such hyperplane section $V$
are $2N$-wise disjoint. Thus, we may always assume the eigenvalues
of $B$ to be positive.

Let $\{\theta_1, \ldots, \theta_d\}$ be a basis of eigenvectors of
$B$, in decreasing order of the eigenvalues $b_1>\cdots> b_d>0$. Let
$V=\pi_v(\extlv{\CC^d}\cap H)$ be as in the statement, where $H$ is
the geometric hyperplane of $\extl{\CC^d}$ defined by some non-zero
$(d-\ell)$-vector $\upsilon$. Let us write
$$
\upsilon=\sum_{I} \upsilon(i_1, \ldots, i_\ell)\,
(\theta_{j_{\ell+1}}\wedge\cdots\wedge\theta_{j_d}),
$$
where the sum is over all sequences $I=(i_1, \ldots, i_\ell)$ with
$1\le i_1<\cdots<i_\ell\le d$, the $\upsilon(I)$ are scalars, and
$j_{\ell+1} < \cdots < j_d$ are the elements of $\{1, \ldots, d\}$
that are \emph{not} in $I$. The assumption that $V$ contains no
eigenspaces of $B$ implies that every $\upsilon(I)$ is non-zero:
otherwise, $\upsilon \wedge (\theta_{i_1}\wedge \cdots \wedge
\theta_{i_\ell})$ would vanish, that is, $\pi_v(\upsilon)$ would
have a non-trivial intersection with the subspace generated by
$\theta_{i_1}$, \ldots, $\theta_{i_\ell}$. Likewise, let us write
\begin{equation}\label{eq.scalars}
\omega = \sum_{I} \omega(i_1, \ldots, i_\ell) \,
(\theta_{i_1}\wedge\cdots\wedge\theta_{i_\ell}),
\end{equation}
where the $\omega(I)$ are scalars. Then $B^{-m}(H)=\{\omega : \omega
\wedge B^{-m}\upsilon =0\}$, and
$$
\omega \wedge B^{-m}\upsilon = \sum_I b_I^{-m} \sigma_I\, \omega(I)
\, \upsilon(I),
$$
where $b_I=b_{j_{\ell+1}}\cdots b_{j_d}>0$ and $\sigma_I =
\theta_{i_1}\wedge\cdots\wedge\theta_{i_\ell} \wedge
\theta_{j_{\ell+1}}\wedge\cdots\wedge\theta_{j_d}$ is either $\pm
1$.

Fix $N=\dim \extl{\CC^d}$ and then let $0\le m_1 < \cdots < m_N$. In
view of the previous paragraph, in order to prove that the
intersection of all the $B^{-m_u}(H)$ is empty it suffices to show
that there does not exist any non-zero $\omega\in \extlv{\CC^d}$
such that
\begin{equation}\label{eq.kernel}
\sum_I b_I^{-m_u}  \sigma_I \, \omega(I) \, \upsilon(I) = 0 \quad
\text{for all } u = 1, \ldots, N,
\end{equation}
that is, such that the vector $(\sigma_I\,\omega(I)\,\upsilon(I))_I$
is in the kernel of $X=\big(b_I^{m_u}\big)_{I,u}$. It is useful to
consider first the special case when the $b_I$ are all distinct (and
positive). Then, by Proposition~\ref{p.vandermonde}, the kernel of
$X$ is trivial. This means that \eqref{eq.kernel} implies
$\sigma_I\,\omega(I)\,\upsilon(I) = 0$ for every $I$. Since
$\sigma_I\,\upsilon(I)$ never vanishes, this means that
$\omega(I)=0$ for every $I$. This proves
Proposition~\ref{p.disjointness} in this case. Notice that this
argument applies to any element $\omega$ of $\extl{\CC^d}$, not only
$\ell$-vectors. Hence, it proves that, under this stronger
assumption, the relation \eqref{eq.kernel} has no non-zero solution
in the whole exterior power $\extl{\CC^d}$.

In general, when the products $b_I$ are not all distinct, condition
\eqref{eq.kernel} may hold on a subspace of $\extl{\CC^d}$ with
positive dimension. The main point in the proof of
Proposition~\ref{p.disjointness} is then to show that this subspace
intersects the set of $\ell$-vectors at the origin only. From
Proposition~\ref{p.vandermonde} we do get that the relation
\eqref{eq.kernel} implies
\begin{equation}\label{eq.kernel2}
\sum_{b_J = b_I} \sigma_J \, \omega(J)\,\upsilon(J)= 0 \quad\text{
for any admissible sequence $I$}
\end{equation}
(\emph{admissible} means that $1\le i_1 < \cdots < i_\ell\le d$),
where the sum is over all admissible sequences $J$ such that
$b_J=b_I$. So, what we really need to prove is

\begin{lemma}\label{l.kernel2}
If an $\ell$-vector $\omega=\omega_1\wedge\cdots\wedge\omega_\ell$
is a solution of \eqref{eq.kernel2} then $\omega(I)=0$ for every
admissible sequence $I=(i_1, \ldots, i_l)$.
\end{lemma}
\begin{proof}
Begin by noting that, for an $\ell$-vector
$\omega=\omega_1\wedge\cdots\wedge\omega_\ell$, the coefficients
$\omega(I)$ in \eqref{eq.scalars} may be expressed in terms of the
vectors $\omega_i$, as follows:
$$
 \omega(I) =
 \left|\begin{array}{ccc}
 \omega_1^{i_1} & \cdots & \omega_1^{i_\ell} \\
 \cdots & \cdots & \cdots \\
 \omega_\ell^{i_1} & \cdots & \omega_\ell^{i_\ell} \\
\end{array}\right|\,,
$$
where $\omega_j=(\omega_j^1, \ldots, \omega_j^d)$. For each $1\le
j\le d$, let $\omega^i=(\omega_1^i, \ldots, \omega_\ell^i)$ be a
column vector. Hence, $\omega(i_1, \ldots, i_\ell)\neq 0$ if and
only if the vectors $\omega^{i_1}$, \ldots, $\omega^{i_\ell}$, are
linearly independent. More generally, given any $1\le s\le \ell$
and $j_1, \ldots, j_s$, we write $\omega(j_1, \ldots, j_s)\neq 0$
to mean the vectors $\omega^{j_1}$, \ldots, $\omega^{j_s}$, are
linearly independent.

Consider first $I=(1, \ldots, \ell)$. Since we assume $b_1 > \cdots
> b_d$, we have $b_I > b_J$ for any admissible sequence $J\neq I$.
Thus, relation \eqref{eq.kernel2} reduces to
$\sigma_I\omega(I)\upsilon(I)=0$. Since $\sigma_I\upsilon(I)$ is
non-zero, that gives $\omega(I)=0$. Now the proof of
Lemma~\ref{l.kernel2} continues by induction: we consider any
admissible sequence $I$, and assume $\omega(J)=0$ for every
admissible sequence $J$ such that $b_J > b_I$. We use the following
simple observation:

\begin{lemma}\label{l.simpletool}
Suppose $\omega(j_1, \ldots, j_s, j, j_{s+1}) = 0$ and
$\omega(j_1, \ldots, j_s, j, j_{s+2})=0$, but $\omega(j_1, \ldots,
j_s, j) \neq 0$. Then $\omega(j_1, \ldots, j_s, j_{s+1}, j_{s+2})
= 0$.
\end{lemma}

\begin{proof}
The assumptions mean that both $\omega^{j_{s+1}}$ and
$\omega^{j_{s+2}}$ are linear combinations of $\{\omega^{j_1},
\ldots, \omega^{j_s}, \omega^j\}$, and so the set $\{\omega^{j_1},
\ldots, \omega^{j_s}, \omega^{j_{s+1}}, \omega^{j_{s+2}}\}$ is
contained in the $(s+1)$-dimensional subspace generated by
$\{\omega^{j_1}, \ldots, \omega^{j_s}, \omega^{j}\}$. This implies
that $\omega(j_1, \ldots, j_s, j_{s+1}, j_{s+2}) = 0$.
\end{proof}

\begin{lemma}\label{l.maintool}
If $\omega(I)\neq 0$ then we have $\omega(j_1, \ldots, j_s, j)=0$
for every $0\le s \le \ell-1$, every $j\notin\{i_1, \ldots,
i_\ell\}$, and every $\{j_1, \ldots, j_s\}\subset\{i_1, \ldots,
i_\ell\}$ that contains all $i_t<j$.
\end{lemma}

\begin{proof}
Consider first the case $\ell-s=1$. Then $(j_1, \ldots, j_s)$
misses exactly one element $i_t$ of $I$, and we have $j < i_t$.
Let $J$ be the admissible sequence obtained by ordering $(j_1,
\ldots, j_s, j)$. Notice that $b_J>b_I$, because $b_j>b_{i_t}$. By
induction, we get that $\omega(J)=0$, as claimed. Now the proof
proceeds by induction on $\ell-s$. Suppose $\ell-s\ge 2$ and let
$j_1, \ldots, j_s, j$ be as in the statement. Choose two different
elements $j_{s+1}$ and $j_{s+2}$ of $\{i_1, \ldots,
i_\ell\}\setminus\{j_1, \ldots, j_s\}$. By induction,
$$
\omega(j_1, \ldots, j_s, j, j_{s+1}) = 0 \quand \omega(j_1,
\ldots, j_s, j, j_{s+2}) = 0.
$$
Suppose $\omega(j_1, \ldots, j_s, j)\neq 0$. Then, we would be
able to use Lemma~\ref{l.simpletool} to conclude that
$$\omega(j_1, \ldots, j_s, j_{s+1}, j_{s+2}) = 0.$$ Since the $j_i$
are distinct elements of $\{i_1, \ldots, i_\ell\}$, that would
imply $\omega(i_1, \ldots, i_\ell)=0$, which would contradict the
hypothesis. This proves that $\omega(j_1, \ldots, j_s, j)= 0$, and
so the proof of Lemma~\ref{l.maintool} is complete.
\end{proof}

\begin{remark}
Notice that $s=0$ is compatible with the other assumptions only if
$i_1>1$. Then the lemma gives that $\omega(j)=0$ or, equivalently,
the column vector $\omega^j=0$, for every $1 \le j < i_1$. This
means that the $\ell$-vector $\omega$ really lives inside a lower
dimensional space, corresponding to coordinates $i_1$ through $d$
only. This case could be easily disposed of, just by assuming
Lemma~\ref{l.kernel2} has already been proved for dimensions
smaller than $d$.
\end{remark}

Let $\prec$ be the usual lexicographical order: $(j_1, \ldots,
j_r) \prec (i_1, \ldots, i_r)$ if and only there exists $0 \le s
\le r-1$ such that $j_1=i_1$, \ldots, $j_s=i_s$, and
$j_{s+1}<i_{s+1}$.

\begin{corollary}\label{c.souma}
If $\omega(I)\neq 0$ then $\omega(J)=0$ for every $J\prec I$.
\end{corollary}

\begin{proof}
Fix $0\le s\le\ell-1$ as in the definition of $J\prec I$, that is,
such that $j_1=i_1$, \ldots, $j_s=i_s$, and $j_{s+1}<i_{s+1}$. By
Lemma~\ref{l.maintool}, we have $\omega(j_1, \ldots, j_s,
j_{s+1})=0$. Consequently, $\omega(j_1, \ldots, j_\ell)=0$, as
claimed.
\end{proof}

Now the inductive step in the proof of Lemma~\ref{l.kernel2} is an
easy consequence. By Corollary~\ref{c.souma}, inside the class of
all sequences $J$ with $b_J=b_I$ there exists at most one $J$ such
that $\omega(J)\neq 0$. Then the relation~\eqref{eq.kernel2}
reduces to $\sigma_J\omega(J)\upsilon(J)=0$. Since
$\sigma_J\upsilon(J)$ never vanishes, this gives $\omega(J)=0$. In
other words, $\omega(J)=0$ for every $J$ such that $b_J=b_I$. This
finishes the proof of Lemma~\ref{l.kernel2}.
\end{proof}

The proof of Proposition~\ref{p.disjointness} is complete.


\end{document}